\documentclass[12pt]{article}

\def\endofps{EndOfTheIncludedPostscriptMagicCookie}
\chardef\other=12
\newwrite\psdumphandle 
\outer\def\psdump#1{\par\medbreak
  \immediate\openout\psdumphandle=#1
  \copytoblankline}
\def\copytoblankline{\begingroup\setupcopy\copypsline}
\def\setupcopy{\def\do##1{\catcode`##1=\other}\dospecials
  \catcode`\\=\other \obeylines}
{\obeylines \gdef\copypsline#1
  {\def\next{#1}%
  \ifx\next\endofps\let\next=\endgroup %
  \else\immediate\write\psdumphandle{\next} \let\next=\copypsline\fi\next}}
\outer\def\closepsdump{
  \immediate\closeout\psdumphandle}

\psdump{l1l3.eps}

/Mathdict 150 dict def
Mathdict begin
/Mlmarg		0 72 mul def
/Mrmarg		0 72 mul def
/Mbmarg		0 72 mul def
/Mtmarg		0 72 mul def
/Mwidth		4 72 mul def
/Mheight	4 72 mul def
/Mtransform	{  } bind def
/Mnodistort	true def
/Mfixwid	true	def
/Mfixdash	false def
/Mrot		0	def
/Mpstart {
MathPictureStart
} bind def
/Mpend {
MathPictureEnd
} bind def
/Mscale {
0 1 0 1
5 -1 roll
MathScale
} bind def
/ISOLatin1Encoding dup where
{ pop pop }
{
[
/.notdef /.notdef /.notdef /.notdef /.notdef /.notdef /.notdef /.notdef
/.notdef /.notdef /.notdef /.notdef /.notdef /.notdef /.notdef /.notdef
/.notdef /.notdef /.notdef /.notdef /.notdef /.notdef /.notdef /.notdef
/.notdef /.notdef /.notdef /.notdef /.notdef /.notdef /.notdef /.notdef
/space /exclam /quotedbl /numbersign /dollar /percent /ampersand /quoteright
/parenleft /parenright /asterisk /plus /comma /minus /period /slash
/zero /one /two /three /four /five /six /seven
/eight /nine /colon /semicolon /less /equal /greater /question
/at /A /B /C /D /E /F /G
/H /I /J /K /L /M /N /O
/P /Q /R /S /T /U /V /W
/X /Y /Z /bracketleft /backslash /bracketright /asciicircum /underscore
/quoteleft /a /b /c /d /e /f /g
/h /i /j /k /l /m /n /o
/p /q /r /s /t /u /v /w
/x /y /z /braceleft /bar /braceright /asciitilde /.notdef
/.notdef /.notdef /.notdef /.notdef /.notdef /.notdef /.notdef /.notdef
/.notdef /.notdef /.notdef /.notdef /.notdef /.notdef /.notdef /.notdef
/dotlessi /grave /acute /circumflex /tilde /macron /breve /dotaccent
/dieresis /.notdef /ring /cedilla /.notdef /hungarumlaut /ogonek /caron
/space /exclamdown /cent /sterling /currency /yen /brokenbar /section
/dieresis /copyright /ordfeminine /guillemotleft
/logicalnot /hyphen /registered /macron
/degree /plusminus /twosuperior /threesuperior
/acute /mu /paragraph /periodcentered
/cedilla /onesuperior /ordmasculine /guillemotright
/onequarter /onehalf /threequarters /questiondown
/Agrave /Aacute /Acircumflex /Atilde /Adieresis /Aring /AE /Ccedilla
/Egrave /Eacute /Ecircumflex /Edieresis /Igrave /Iacute /Icircumflex /Idieresis
/Eth /Ntilde /Ograve /Oacute /Ocircumflex /Otilde /Odieresis /multiply
/Oslash /Ugrave /Uacute /Ucircumflex /Udieresis /Yacute /Thorn /germandbls
/agrave /aacute /acircumflex /atilde /adieresis /aring /ae /ccedilla
/egrave /eacute /ecircumflex /edieresis /igrave /iacute /icircumflex /idieresis
/eth /ntilde /ograve /oacute /ocircumflex /otilde /odieresis /divide
/oslash /ugrave /uacute /ucircumflex /udieresis /yacute /thorn /ydieresis
] def
} ifelse
/MFontDict 50 dict def
/MStrCat
{
exch
dup length
2 index length add
string
dup 3 1 roll
copy
length
exch dup
4 2 roll exch
putinterval
} def
/MCreateEncoding
{
1 index
255 string cvs
(-) MStrCat
1 index MStrCat
cvn exch
(Encoding) MStrCat
cvn dup where
{
exch get
}
{
pop
StandardEncoding
} ifelse
3 1 roll
dup MFontDict exch known not
{
1 index findfont
dup length dict
begin
{1 index /FID ne
{def}
{pop pop}
ifelse} forall
/Encoding 3 index
def
currentdict
end
1 index exch definefont pop
MFontDict 1 index
null put
}
if
exch pop
exch pop
} def
/ISOLatin1 { (ISOLatin1) MCreateEncoding } def
/ISO8859 { (ISOLatin1) MCreateEncoding } def
/Mcopyfont {
dup
maxlength
dict
exch
{
1 index
/FID
eq
{
pop pop
}
{
2 index
3 1 roll
put
}
ifelse
}
forall
} def
/Plain	/Courier findfont Mcopyfont definefont pop
/Bold	/Courier-Bold findfont Mcopyfont definefont pop
/Italic /Courier-Oblique findfont Mcopyfont definefont pop
/MathPictureStart {
gsave
Mtransform
Mlmarg
Mbmarg
translate
Mwidth
Mlmarg Mrmarg add
sub
/Mwidth exch def
Mheight
Mbmarg Mtmarg add
sub
/Mheight exch def
/Mtmatrix
matrix currentmatrix
def
/Mgmatrix
matrix currentmatrix
def
} bind def
/MathPictureEnd {
grestore
} bind def
/MFill {
0 0 		moveto
Mwidth 0 	lineto
Mwidth Mheight 	lineto
0 Mheight 	lineto
fill
} bind def
/MPlotRegion {
3 index
Mwidth mul
2 index
Mheight mul
translate
exch sub
Mheight mul
/Mheight
exch def
exch sub
Mwidth mul
/Mwidth
exch def
} bind def
/MathSubStart {
Momatrix
Mgmatrix Mtmatrix
Mwidth Mheight
7 -2 roll
moveto
Mtmatrix setmatrix
currentpoint
Mgmatrix setmatrix
9 -2 roll
moveto
Mtmatrix setmatrix
currentpoint
2 copy translate
/Mtmatrix matrix currentmatrix def
3 -1 roll
exch sub
/Mheight exch def
sub
/Mwidth exch def
} bind def
/MathSubEnd {
/Mheight exch def
/Mwidth exch def
/Mtmatrix exch def
dup setmatrix
/Mgmatrix exch def
/Momatrix exch def
} bind def
/Mdot {
moveto
0 0 rlineto
stroke
} bind def
/Mtetra {
moveto
lineto
lineto
lineto
fill
} bind def
/Metetra {
moveto
lineto
lineto
lineto
closepath
gsave
fill
grestore
0 setgray
stroke
} bind def
/Mistroke {
flattenpath
0 0 0
{
4 2 roll
pop pop
}
{
4 -1 roll
2 index
sub dup mul
4 -1 roll
2 index
sub dup mul
add sqrt
4 -1 roll
add
3 1 roll
}
{
stop
}
{
stop
}
pathforall
pop pop
currentpoint
stroke
moveto
currentdash
3 -1 roll
add
setdash
} bind def
/Mfstroke {
stroke
currentdash
pop 0
setdash
} bind def
/Mrotsboxa {
gsave
dup
/Mrot
exch def
Mrotcheck
Mtmatrix
dup
setmatrix
7 1 roll
4 index
4 index
translate
rotate
3 index
-1 mul
3 index
-1 mul
translate
/Mtmatrix
matrix
currentmatrix
def
grestore
Msboxa
3  -1 roll
/Mtmatrix
exch def
/Mrot
0 def
} bind def
/Msboxa {
newpath
5 -1 roll
Mvboxa
pop
Mboxout
6 -1 roll
5 -1 roll
4 -1 roll
Msboxa1
5 -3 roll
Msboxa1
Mboxrot
[
7 -2 roll
2 copy
[
3 1 roll
10 -1 roll
9 -1 roll
]
6 1 roll
5 -2 roll
]
} bind def
/Msboxa1 {
sub
2 div
dup
2 index
1 add
mul
3 -1 roll
-1 add
3 -1 roll
mul
} bind def
/Mvboxa	{
Mfixwid
{
Mvboxa1
}
{
dup
Mwidthcal
0 exch
{
add
}
forall
exch
Mvboxa1
4 index
7 -1 roll
add
4 -1 roll
pop
3 1 roll
}
ifelse
} bind def
/Mvboxa1 {
gsave
newpath
[ true
3 -1 roll
{
Mbbox
5 -1 roll
{
0
5 1 roll
}
{
7 -1 roll
exch sub
(m) stringwidth pop
.3 mul
sub
7 1 roll
6 -1 roll
4 -1 roll
Mmin
3 -1 roll
5 index
add
5 -1 roll
4 -1 roll
Mmax
4 -1 roll
}
ifelse
false
}
forall
{ stop } if
counttomark
1 add
4 roll
]
grestore
} bind def
/Mbbox {
1 dict begin
0 0 moveto
/temp (T) def
{	gsave
currentpoint newpath moveto
temp 0 3 -1 roll put temp
false charpath flattenpath currentpoint
pathbbox
grestore moveto lineto moveto} forall
pathbbox
newpath
end
} bind def
/Mmin {
2 copy
gt
{ exch } if
pop
} bind def
/Mmax {
2 copy
lt
{ exch } if
pop
} bind def
/Mrotshowa {
dup
/Mrot
exch def
Mrotcheck
Mtmatrix
dup
setmatrix
7 1 roll
4 index
4 index
translate
rotate
3 index
-1 mul
3 index
-1 mul
translate
/Mtmatrix
matrix
currentmatrix
def
Mgmatrix setmatrix
Mshowa
/Mtmatrix
exch def
/Mrot 0 def
} bind def
/Mshowa {
4 -2 roll
moveto
2 index
Mtmatrix setmatrix
Mvboxa
7 1 roll
Mboxout
6 -1 roll
5 -1 roll
4 -1 roll
Mshowa1
4 1 roll
Mshowa1
rmoveto
currentpoint
Mfixwid
{
Mshowax
}
{
Mshoway
}
ifelse
pop pop pop pop
Mgmatrix setmatrix
} bind def
/Mshowax {	
0 1
4 index length
-1 add
{
2 index
4 index
2 index
get
3 index
add
moveto
4 index
exch get
Mfixdash
{
Mfixdashp		
}
if
show
} for
} bind def
/Mfixdashp {
dup
length
1
gt
1 index
true exch
{
45
eq
and
} forall
and
{
gsave
(--)		
stringwidth pop
(-)
stringwidth pop
sub
2 div
0 rmoveto
dup
length
1 sub
{
(-)
show
}
repeat
grestore
}
if
} bind def	
/Mshoway {
3 index
Mwidthcal
5 1 roll
0 1
4 index length
-1 add
{
2 index
4 index
2 index
get
3 index
add
moveto
4 index
exch get
[
6 index
aload
length
2 add
-1 roll
{
pop
Strform
stringwidth
pop
neg 
exch
add
0 rmoveto
}
exch
kshow
cleartomark
} for
pop	
} bind def
/Mwidthcal {
[
exch
{
Mwidthcal1
}
forall
]
[
exch
dup
Maxlen
-1 add
0 1
3 -1 roll
{
[
exch
2 index
{		 
1 index	
Mget
exch
}
forall		
pop
Maxget
exch
}
for
pop
]
Mreva
} bind def
/Mreva	{
[
exch
aload
length
-1 1
{1 roll}
for
]
} bind def
/Mget	{
1 index
length
-1 add
1 index
ge
{
get
}
{
pop pop
0
}
ifelse
} bind def
/Maxlen	{
[
exch
{
length
}
forall
Maxget
} bind def
/Maxget	{
counttomark
-1 add
1 1
3 -1 roll
{
pop
Mmax
}
for
exch
pop
} bind def
/Mwidthcal1 {
[
exch
{
Strform
stringwidth
pop
}
forall
]
} bind def
/Strform {
/tem (x) def
tem 0
3 -1 roll
put
tem
} bind def
/Mshowa1 {
2 copy
add
4 1 roll
sub
mul
sub
-2 div
} bind def
/MathScale {
Mwidth
Mheight
Mlp
translate
scale
/yscale exch def
/ybias exch def
/xscale exch def
/xbias exch def
/Momatrix
xscale yscale matrix scale
xbias ybias matrix translate
matrix concatmatrix def
/Mgmatrix
matrix currentmatrix
def
} bind def
/Mlp {
3 copy
Mlpfirst
{
Mnodistort
{
Mmin
dup
} if
4 index
2 index
2 index
Mlprun
11 index
11 -1 roll
10 -4 roll
Mlp1
8 index
9 -5 roll
Mlp1
4 -1 roll
and
{ exit } if
3 -1 roll
pop pop
} loop
exch
3 1 roll
7 -3 roll
pop pop pop
} bind def
/Mlpfirst {
3 -1 roll
dup length
2 copy
-2 add
get
aload
pop pop pop
4 -2 roll
-1 add
get
aload
pop pop pop
6 -1 roll
3 -1 roll
5 -1 roll
sub
div
4 1 roll
exch sub
div
} bind def
/Mlprun {
2 copy
4 index
0 get
dup
4 1 roll
Mlprun1
3 copy
8 -2 roll
9 -1 roll
{
3 copy
Mlprun1
3 copy
11 -3 roll
/gt Mlpminmax
8 3 roll
11 -3 roll
/lt Mlpminmax
8 3 roll
} forall
pop pop pop pop
3 1 roll
pop pop
aload pop
5 -1 roll
aload pop
exch
6 -1 roll
Mlprun2
8 2 roll
4 -1 roll
Mlprun2
6 2 roll
3 -1 roll
Mlprun2
4 2 roll
exch
Mlprun2
6 2 roll
} bind def
/Mlprun1 {
aload pop
exch
6 -1 roll
5 -1 roll
mul add
4 -2 roll
mul
3 -1 roll
add
} bind def
/Mlprun2 {
2 copy
add 2 div
3 1 roll
exch sub
} bind def
/Mlpminmax {
cvx
2 index
6 index
2 index
exec
{
7 -3 roll
4 -1 roll
} if
1 index
5 index
3 -1 roll
exec
{
4 1 roll
pop
5 -1 roll
aload
pop pop
4 -1 roll
aload pop
[
8 -2 roll
pop
5 -2 roll
pop
6 -2 roll
pop
5 -1 roll
]
4 1 roll
pop
}
{
pop pop pop
} ifelse
} bind def
/Mlp1 {
5 index
3 index	sub
5 index
2 index mul
1 index
le
1 index
0 le
or
dup
not
{
1 index
3 index	div
.99999 mul
8 -1 roll
pop
7 1 roll
}
if
8 -1 roll
2 div
7 -2 roll
pop sub
5 index
6 -3 roll
pop pop
mul sub
exch
} bind def
/intop 0 def
/inrht 0 def
/inflag 0 def
/outflag 0 def
/xadrht 0 def
/xadlft 0 def
/yadtop 0 def
/yadbot 0 def
/Minner {
outflag
1
eq
{
/outflag 0 def
/intop 0 def
/inrht 0 def
} if		
5 index
gsave
Mtmatrix setmatrix
Mvboxa pop
grestore
3 -1 roll
pop
dup
intop
gt
{
/intop
exch def
}
{ pop }
ifelse
dup
inrht
gt
{
/inrht
exch def
}
{ pop }
ifelse
pop
/inflag
1 def
} bind def
/Mouter {
/xadrht 0 def
/xadlft 0 def
/yadtop 0 def
/yadbot 0 def
inflag
1 eq
{
dup
0 lt
{
dup
intop
mul
neg
/yadtop
exch def
} if
dup
0 gt
{
dup
intop
mul
/yadbot
exch def
}
if
pop
dup
0 lt
{
dup
inrht
mul
neg
/xadrht
exch def
} if
dup
0 gt
{
dup
inrht
mul
/xadlft
exch def
} if
pop
/outflag 1 def
}
{ pop pop}
ifelse
/inflag 0 def
/inrht 0 def
/intop 0 def
} bind def	
/Mboxout {
outflag
1
eq
{
4 -1
roll
xadlft
leadjust
add
sub
4 1 roll
3 -1
roll
yadbot
leadjust
add
sub
3 1
roll
exch
xadrht
leadjust
add
add
exch
yadtop
leadjust
add
add
/outflag 0 def
/xadlft 0 def
/yadbot 0 def
/xadrht 0 def
/yadtop 0 def
} if
} bind def
/leadjust {
(m) stringwidth pop
.5 mul
} bind def
/Mrotcheck {
dup
90
eq
{
yadbot
/yadbot
xadrht
def	
/xadrht
yadtop
def
/yadtop
xadlft
def
/xadlft
exch
def
}
if
dup
cos
1 index
sin
Checkaux
dup
cos
1 index
sin neg
exch
Checkaux
3 1 roll
pop pop
} bind def
/Checkaux {
4 index
exch
4 index
mul
3 1 roll
mul add
4 1 roll
} bind def
/Mboxrot {
Mrot
90 eq
{
brotaux
4 2
roll
}
if
Mrot
180 eq
{
4 2
roll
brotaux
4 2
roll
brotaux
}	
if
Mrot
270 eq
{
4 2
roll
brotaux
}
if
} bind def
/brotaux {
neg
exch
neg
} bind def
/Mabsproc {
0
matrix defaultmatrix
dtransform idtransform
dup mul exch
dup mul
add sqrt
} bind def
/Mabswid {
Mabsproc
setlinewidth	
} bind def
/Mabsdash {
exch
[
exch
{
Mabsproc
}
forall
]
exch
setdash
} bind def
/MBeginOrig { Momatrix concat} bind def
/MEndOrig { Mgmatrix setmatrix} bind def
/sampledsound where
{ pop}
{ /sampledsound {
exch
pop
exch
5 1 roll
mul
4 idiv
mul
2 idiv
exch pop
exch
/Mtempproc exch def
{ Mtempproc pop}
repeat
} bind def
} ifelse
/g { setgray} bind def
/k { setcmykcolor} bind def
/m { moveto} bind def
/p { gsave} bind def
/r { setrgbcolor} bind def
/w { setlinewidth} bind def
/C { curveto} bind def
/F { fill} bind def
/L { lineto} bind def
/P { grestore} bind def
/s { stroke} bind def
/setcmykcolor where
{ pop}
{ /setcmykcolor {
4 1
roll
[
4 1 
roll
]
{
1 index
sub
1
sub neg
dup
0
lt
{
pop
0
}
if
dup
1
gt
{
pop
1
}
if
exch
} forall
pop
setrgbcolor
} bind def
} ifelse
/Mcharproc
{
currentfile
(x)
readhexstring
pop
0 get
exch
div
} bind def
/Mshadeproc
{
dup
3 1
roll
{
dup
Mcharproc
3 1
roll
} repeat
1 eq
{
setgray
}
{
3 eq
{
setrgbcolor
}
{
setcmykcolor
} ifelse
} ifelse
} bind def
/Mrectproc
{
3 index
2 index
moveto
2 index
3 -1
roll
lineto
dup
3 1
roll
lineto
lineto
fill
} bind def
/Mcolorimage
{
7 1
roll
pop
pop
matrix
invertmatrix
concat
2 exch exp
1 sub
3 1 roll
1 1
2 index
{
1 1
4 index
{
dup
1 sub
exch
2 index
dup
1 sub
exch
7 index
9 index
Mshadeproc
Mrectproc
} for
pop
} for
pop pop pop pop
} bind def
/Mimage
{
pop
matrix
invertmatrix
concat
2 exch exp
1 sub
3 1 roll
1 1
2 index
{
1 1
4 index
{
dup
1 sub
exch
2 index
dup
1 sub
exch
7 index
Mcharproc
setgray
Mrectproc
} for
pop
} for
pop pop pop
} bind def

MathPictureStart
/Courier findfont 10  scalefont  setfont
-0.928571 0.952381 0 6.18034e-05 [
[(1)] .02381 0 0 2 Msboxa
[(1.2)] .21429 0 0 2 Msboxa
[(1.4)] .40476 0 0 2 Msboxa
[(1.6)] .59524 0 0 2 Msboxa
[(1.8)] .78571 0 0 2 Msboxa
[(2)] .97619 0 0 2 Msboxa
[(2000)] .01131 .12361 1 0 Msboxa
[(4000)] .01131 .24721 1 0 Msboxa
[(6000)] .01131 .37082 1 0 Msboxa
[(8000)] .01131 .49443 1 0 Msboxa
[(10000)] .01131 .61803 1 0 Msboxa
[ -0.001 -0.001 0 0 ]
[ 1.001 .61903 0 0 ]
] MathScale
1 setlinecap
1 setlinejoin
newpath
[ ] 0 setdash
0 g
p
p
.002 w
.02381 0 m
.02381 .00625 L
s
P
[(1)] .02381 0 0 2 Mshowa
p
.002 w
.21429 0 m
.21429 .00625 L
s
P
[(1.2)] .21429 0 0 2 Mshowa
p
.002 w
.40476 0 m
.40476 .00625 L
s
P
[(1.4)] .40476 0 0 2 Mshowa
p
.002 w
.59524 0 m
.59524 .00625 L
s
P
[(1.6)] .59524 0 0 2 Mshowa
p
.002 w
.78571 0 m
.78571 .00625 L
s
P
[(1.8)] .78571 0 0 2 Mshowa
p
.002 w
.97619 0 m
.97619 .00625 L
s
P
[(2)] .97619 0 0 2 Mshowa
p
.001 w
.0619 0 m
.0619 .00375 L
s
P
p
.001 w
.1 0 m
.1 .00375 L
s
P
p
.001 w
.1381 0 m
.1381 .00375 L
s
P
p
.001 w
.17619 0 m
.17619 .00375 L
s
P
p
.001 w
.25238 0 m
.25238 .00375 L
s
P
p
.001 w
.29048 0 m
.29048 .00375 L
s
P
p
.001 w
.32857 0 m
.32857 .00375 L
s
P
p
.001 w
.36667 0 m
.36667 .00375 L
s
P
p
.001 w
.44286 0 m
.44286 .00375 L
s
P
p
.001 w
.48095 0 m
.48095 .00375 L
s
P
p
.001 w
.51905 0 m
.51905 .00375 L
s
P
p
.001 w
.55714 0 m
.55714 .00375 L
s
P
p
.001 w
.63333 0 m
.63333 .00375 L
s
P
p
.001 w
.67143 0 m
.67143 .00375 L
s
P
p
.001 w
.70952 0 m
.70952 .00375 L
s
P
p
.001 w
.74762 0 m
.74762 .00375 L
s
P
p
.001 w
.82381 0 m
.82381 .00375 L
s
P
p
.001 w
.8619 0 m
.8619 .00375 L
s
P
p
.001 w
.9 0 m
.9 .00375 L
s
P
p
.001 w
.9381 0 m
.9381 .00375 L
s
P
p
.002 w
0 0 m
1 0 L
s
P
p
.002 w
.02381 .12361 m
.03006 .12361 L
s
P
[(2000)] .01131 .12361 1 0 Mshowa
p
.002 w
.02381 .24721 m
.03006 .24721 L
s
P
[(4000)] .01131 .24721 1 0 Mshowa
p
.002 w
.02381 .37082 m
.03006 .37082 L
s
P
[(6000)] .01131 .37082 1 0 Mshowa
p
.002 w
.02381 .49443 m
.03006 .49443 L
s
P
[(8000)] .01131 .49443 1 0 Mshowa
p
.002 w
.02381 .61803 m
.03006 .61803 L
s
P
[(10000)] .01131 .61803 1 0 Mshowa
p
.001 w
.02381 .02472 m
.02756 .02472 L
s
P
p
.001 w
.02381 .04944 m
.02756 .04944 L
s
P
p
.001 w
.02381 .07416 m
.02756 .07416 L
s
P
p
.001 w
.02381 .09889 m
.02756 .09889 L
s
P
p
.001 w
.02381 .14833 m
.02756 .14833 L
s
P
p
.001 w
.02381 .17305 m
.02756 .17305 L
s
P
p
.001 w
.02381 .19777 m
.02756 .19777 L
s
P
p
.001 w
.02381 .22249 m
.02756 .22249 L
s
P
p
.001 w
.02381 .27193 m
.02756 .27193 L
s
P
p
.001 w
.02381 .29666 m
.02756 .29666 L
s
P
p
.001 w
.02381 .32138 m
.02756 .32138 L
s
P
p
.001 w
.02381 .3461 m
.02756 .3461 L
s
P
p
.001 w
.02381 .39554 m
.02756 .39554 L
s
P
p
.001 w
.02381 .42026 m
.02756 .42026 L
s
P
p
.001 w
.02381 .44498 m
.02756 .44498 L
s
P
p
.001 w
.02381 .46971 m
.02756 .46971 L
s
P
p
.001 w
.02381 .51915 m
.02756 .51915 L
s
P
p
.001 w
.02381 .54387 m
.02756 .54387 L
s
P
p
.001 w
.02381 .56859 m
.02756 .56859 L
s
P
p
.001 w
.02381 .59331 m
.02756 .59331 L
s
P
p
.002 w
.02381 0 m
.02381 .61803 L
s
P
P
0 0 m
1 0 L
1 .61803 L
0 .61803 L
closepath
clip
newpath
p
p
[ .01 .01 ] 0 setdash
.5 g
p
.002 w
s
s
s
s
s
s
s
s
s
s
s
s
.49675 .61803 m
.5 .61282 L
s
.5 .61282 m
.53968 .55717 L
.57937 .50832 L
.61905 .46527 L
.65873 .42721 L
.69841 .39343 L
.7381 .36338 L
.77778 .33655 L
.81746 .31253 L
.85714 .29097 L
.89683 .27156 L
.91667 .26258 L
.92659 .25826 L
.93155 .25614 L
.93651 .25405 L
.93775 .25353 L
.93899 .25302 L
.94023 .2525 L
.94147 .25235 L
.94395 .25229 L
.94643 .25224 L
.95635 .25204 L
.97619 .25165 L
s
P
P
p
[ .01 .01 ] 0 setdash
.5 g
p
.002 w
s
s
s
s
.17833 .61803 m
.18254 .60921 L
s
.18254 .60921 m
.22222 .5392 L
.2619 .47992 L
.30159 .42942 L
.34127 .38616 L
.35119 .37633 L
.36111 .36685 L
.36607 .36224 L
.37103 .35772 L
.37227 .3566 L
.37351 .35573 L
.37599 .35502 L
.38095 .35361 L
.42063 .34314 L
.46032 .33394 L
.5 .32583 L
.53968 .31864 L
.57937 .31224 L
.61905 .30653 L
.65873 .30141 L
.69841 .29681 L
.7381 .29266 L
.77778 .28891 L
.81746 .2855 L
.85714 .2824 L
.89683 .27958 L
.93651 .277 L
.97619 .27463 L
s
P
P
p
.5 g
p
.002 w
.02381 .57 m
.06349 .52676 L
.10317 .49082 L
.14286 .46072 L
.18254 .43531 L
.22222 .41371 L
.2619 .39524 L
.30159 .37934 L
.34127 .36558 L
.35119 .36244 L
.36111 .35939 L
.36607 .35791 L
.37103 .35645 L
.37227 .35609 L
.37351 .35549 L
.37599 .35328 L
.38095 .34891 L
.42063 .31668 L
.46032 .28866 L
.48016 .27602 L
.48512 .27299 L
.49008 .27001 L
.49132 .26927 L
.49256 .26853 L
.4938 .2678 L
.49504 .26754 L
.49628 .26748 L
.49752 .26741 L
.5 .26728 L
.53968 .26525 L
.57937 .2634 L
.61905 .26173 L
.65873 .2602 L
.69841 .2588 L
.7381 .25752 L
.77778 .25634 L
.81746 .25525 L
.85714 .25424 L
.89683 .25331 L
.91667 .25287 L
.92659 .25266 L
.93155 .25255 L
.93403 .2525 L
.93651 .25245 L
.93775 .25242 L
.93899 .2524 L
.94023 .25237 L
.94147 .25198 L
.94395 .25096 L
.94643 .24994 L
Mistroke
.95635 .24594 L
.97619 .23821 L
Mfstroke
P
P
p
.5 g
p
.002 w
.02381 .32476 m
.06349 .28473 L
.10317 .25165 L
.14286 .22409 L
.18254 .20096 L
.22222 .18141 L
.2619 .16478 L
.30159 .15055 L
.34127 .13831 L
.38095 .12772 L
.42063 .11851 L
.46032 .11047 L
.5 .10341 L
.53968 .09719 L
.57937 .09169 L
.61905 .08681 L
.65873 .08246 L
.69841 .07857 L
.7381 .07509 L
.77778 .07195 L
.81746 .06912 L
.85714 .06656 L
.89683 .06424 L
.93651 .06213 L
.97619 .0602 L
s
P
P
p
p
.004 w
.02381 .32476 m
.06349 .3155 L
.10317 .3076 L
.14286 .30081 L
.18254 .29494 L
.22222 .28983 L
.2619 .28534 L
.30159 .2814 L
.34127 .2779 L
.38095 .27479 L
.42063 .27201 L
.44048 .27073 L
.46032 .26952 L
.48016 .26837 L
.48512 .26809 L
.4876 .26795 L
.49008 .26782 L
.49132 .26775 L
.49256 .26768 L
.4938 .26761 L
.49504 .26707 L
.49628 .26635 L
.49752 .26563 L
.5 .26419 L
.53968 .24274 L
.57937 .22386 L
.61905 .20718 L
.65873 .19239 L
.69841 .17924 L
.7381 .1675 L
.77778 .15699 L
.81746 .14756 L
.85714 .13907 L
.89683 .1314 L
.93651 .12447 L
.97619 .11818 L
s
P
P
p
p
.004 w
.02381 .07952 m
.06349 .07348 L
.10317 .06843 L
.14286 .06419 L
.18254 .06059 L
.22222 .05752 L
.2619 .05489 L
.30159 .05261 L
.34127 .05063 L
.38095 .0489 L
.42063 .04738 L
.46032 .04604 L
.5 .04485 L
.53968 .0438 L
.57937 .04286 L
.61905 .04201 L
.65873 .04125 L
.69841 .04057 L
.7381 .03995 L
.77778 .03938 L
.81746 .03887 L
.85714 .0384 L
.89683 .03798 L
.93651 .03758 L
.97619 .03722 L
s
P
P
P
MathPictureEnd
end
EndOfTheIncludedPostscriptMagicCookie

\closepsdump

\psdump{l1n1.eps}

/Mathdict 150 dict def
Mathdict begin
/Mlmarg		0 72 mul def
/Mrmarg		0 72 mul def
/Mbmarg		0 72 mul def
/Mtmarg		0 72 mul def
/Mwidth		4 72 mul def
/Mheight	4 72 mul def
/Mtransform	{  } bind def
/Mnodistort	true def
/Mfixwid	true	def
/Mfixdash	false def
/Mrot		0	def
/Mpstart {
MathPictureStart
} bind def
/Mpend {
MathPictureEnd
} bind def
/Mscale {
0 1 0 1
5 -1 roll
MathScale
} bind def
/ISOLatin1Encoding dup where
{ pop pop }
{
[
/.notdef /.notdef /.notdef /.notdef /.notdef /.notdef /.notdef /.notdef
/.notdef /.notdef /.notdef /.notdef /.notdef /.notdef /.notdef /.notdef
/.notdef /.notdef /.notdef /.notdef /.notdef /.notdef /.notdef /.notdef
/.notdef /.notdef /.notdef /.notdef /.notdef /.notdef /.notdef /.notdef
/space /exclam /quotedbl /numbersign /dollar /percent /ampersand /quoteright
/parenleft /parenright /asterisk /plus /comma /minus /period /slash
/zero /one /two /three /four /five /six /seven
/eight /nine /colon /semicolon /less /equal /greater /question
/at /A /B /C /D /E /F /G
/H /I /J /K /L /M /N /O
/P /Q /R /S /T /U /V /W
/X /Y /Z /bracketleft /backslash /bracketright /asciicircum /underscore
/quoteleft /a /b /c /d /e /f /g
/h /i /j /k /l /m /n /o
/p /q /r /s /t /u /v /w
/x /y /z /braceleft /bar /braceright /asciitilde /.notdef
/.notdef /.notdef /.notdef /.notdef /.notdef /.notdef /.notdef /.notdef
/.notdef /.notdef /.notdef /.notdef /.notdef /.notdef /.notdef /.notdef
/dotlessi /grave /acute /circumflex /tilde /macron /breve /dotaccent
/dieresis /.notdef /ring /cedilla /.notdef /hungarumlaut /ogonek /caron
/space /exclamdown /cent /sterling /currency /yen /brokenbar /section
/dieresis /copyright /ordfeminine /guillemotleft
/logicalnot /hyphen /registered /macron
/degree /plusminus /twosuperior /threesuperior
/acute /mu /paragraph /periodcentered
/cedilla /onesuperior /ordmasculine /guillemotright
/onequarter /onehalf /threequarters /questiondown
/Agrave /Aacute /Acircumflex /Atilde /Adieresis /Aring /AE /Ccedilla
/Egrave /Eacute /Ecircumflex /Edieresis /Igrave /Iacute /Icircumflex /Idieresis
/Eth /Ntilde /Ograve /Oacute /Ocircumflex /Otilde /Odieresis /multiply
/Oslash /Ugrave /Uacute /Ucircumflex /Udieresis /Yacute /Thorn /germandbls
/agrave /aacute /acircumflex /atilde /adieresis /aring /ae /ccedilla
/egrave /eacute /ecircumflex /edieresis /igrave /iacute /icircumflex /idieresis
/eth /ntilde /ograve /oacute /ocircumflex /otilde /odieresis /divide
/oslash /ugrave /uacute /ucircumflex /udieresis /yacute /thorn /ydieresis
] def
} ifelse
/MFontDict 50 dict def
/MStrCat
{
exch
dup length
2 index length add
string
dup 3 1 roll
copy
length
exch dup
4 2 roll exch
putinterval
} def
/MCreateEncoding
{
1 index
255 string cvs
(-) MStrCat
1 index MStrCat
cvn exch
(Encoding) MStrCat
cvn dup where
{
exch get
}
{
pop
StandardEncoding
} ifelse
3 1 roll
dup MFontDict exch known not
{
1 index findfont
dup length dict
begin
{1 index /FID ne
{def}
{pop pop}
ifelse} forall
/Encoding 3 index
def
currentdict
end
1 index exch definefont pop
MFontDict 1 index
null put
}
if
exch pop
exch pop
} def
/ISOLatin1 { (ISOLatin1) MCreateEncoding } def
/ISO8859 { (ISOLatin1) MCreateEncoding } def
/Mcopyfont {
dup
maxlength
dict
exch
{
1 index
/FID
eq
{
pop pop
}
{
2 index
3 1 roll
put
}
ifelse
}
forall
} def
/Plain	/Courier findfont Mcopyfont definefont pop
/Bold	/Courier-Bold findfont Mcopyfont definefont pop
/Italic /Courier-Oblique findfont Mcopyfont definefont pop
/MathPictureStart {
gsave
Mtransform
Mlmarg
Mbmarg
translate
Mwidth
Mlmarg Mrmarg add
sub
/Mwidth exch def
Mheight
Mbmarg Mtmarg add
sub
/Mheight exch def
/Mtmatrix
matrix currentmatrix
def
/Mgmatrix
matrix currentmatrix
def
} bind def
/MathPictureEnd {
grestore
} bind def
/MFill {
0 0 		moveto
Mwidth 0 	lineto
Mwidth Mheight 	lineto
0 Mheight 	lineto
fill
} bind def
/MPlotRegion {
3 index
Mwidth mul
2 index
Mheight mul
translate
exch sub
Mheight mul
/Mheight
exch def
exch sub
Mwidth mul
/Mwidth
exch def
} bind def
/MathSubStart {
Momatrix
Mgmatrix Mtmatrix
Mwidth Mheight
7 -2 roll
moveto
Mtmatrix setmatrix
currentpoint
Mgmatrix setmatrix
9 -2 roll
moveto
Mtmatrix setmatrix
currentpoint
2 copy translate
/Mtmatrix matrix currentmatrix def
3 -1 roll
exch sub
/Mheight exch def
sub
/Mwidth exch def
} bind def
/MathSubEnd {
/Mheight exch def
/Mwidth exch def
/Mtmatrix exch def
dup setmatrix
/Mgmatrix exch def
/Momatrix exch def
} bind def
/Mdot {
moveto
0 0 rlineto
stroke
} bind def
/Mtetra {
moveto
lineto
lineto
lineto
fill
} bind def
/Metetra {
moveto
lineto
lineto
lineto
closepath
gsave
fill
grestore
0 setgray
stroke
} bind def
/Mistroke {
flattenpath
0 0 0
{
4 2 roll
pop pop
}
{
4 -1 roll
2 index
sub dup mul
4 -1 roll
2 index
sub dup mul
add sqrt
4 -1 roll
add
3 1 roll
}
{
stop
}
{
stop
}
pathforall
pop pop
currentpoint
stroke
moveto
currentdash
3 -1 roll
add
setdash
} bind def
/Mfstroke {
stroke
currentdash
pop 0
setdash
} bind def
/Mrotsboxa {
gsave
dup
/Mrot
exch def
Mrotcheck
Mtmatrix
dup
setmatrix
7 1 roll
4 index
4 index
translate
rotate
3 index
-1 mul
3 index
-1 mul
translate
/Mtmatrix
matrix
currentmatrix
def
grestore
Msboxa
3  -1 roll
/Mtmatrix
exch def
/Mrot
0 def
} bind def
/Msboxa {
newpath
5 -1 roll
Mvboxa
pop
Mboxout
6 -1 roll
5 -1 roll
4 -1 roll
Msboxa1
5 -3 roll
Msboxa1
Mboxrot
[
7 -2 roll
2 copy
[
3 1 roll
10 -1 roll
9 -1 roll
]
6 1 roll
5 -2 roll
]
} bind def
/Msboxa1 {
sub
2 div
dup
2 index
1 add
mul
3 -1 roll
-1 add
3 -1 roll
mul
} bind def
/Mvboxa	{
Mfixwid
{
Mvboxa1
}
{
dup
Mwidthcal
0 exch
{
add
}
forall
exch
Mvboxa1
4 index
7 -1 roll
add
4 -1 roll
pop
3 1 roll
}
ifelse
} bind def
/Mvboxa1 {
gsave
newpath
[ true
3 -1 roll
{
Mbbox
5 -1 roll
{
0
5 1 roll
}
{
7 -1 roll
exch sub
(m) stringwidth pop
.3 mul
sub
7 1 roll
6 -1 roll
4 -1 roll
Mmin
3 -1 roll
5 index
add
5 -1 roll
4 -1 roll
Mmax
4 -1 roll
}
ifelse
false
}
forall
{ stop } if
counttomark
1 add
4 roll
]
grestore
} bind def
/Mbbox {
1 dict begin
0 0 moveto
/temp (T) def
{	gsave
currentpoint newpath moveto
temp 0 3 -1 roll put temp
false charpath flattenpath currentpoint
pathbbox
grestore moveto lineto moveto} forall
pathbbox
newpath
end
} bind def
/Mmin {
2 copy
gt
{ exch } if
pop
} bind def
/Mmax {
2 copy
lt
{ exch } if
pop
} bind def
/Mrotshowa {
dup
/Mrot
exch def
Mrotcheck
Mtmatrix
dup
setmatrix
7 1 roll
4 index
4 index
translate
rotate
3 index
-1 mul
3 index
-1 mul
translate
/Mtmatrix
matrix
currentmatrix
def
Mgmatrix setmatrix
Mshowa
/Mtmatrix
exch def
/Mrot 0 def
} bind def
/Mshowa {
4 -2 roll
moveto
2 index
Mtmatrix setmatrix
Mvboxa
7 1 roll
Mboxout
6 -1 roll
5 -1 roll
4 -1 roll
Mshowa1
4 1 roll
Mshowa1
rmoveto
currentpoint
Mfixwid
{
Mshowax
}
{
Mshoway
}
ifelse
pop pop pop pop
Mgmatrix setmatrix
} bind def
/Mshowax {	
0 1
4 index length
-1 add
{
2 index
4 index
2 index
get
3 index
add
moveto
4 index
exch get
Mfixdash
{
Mfixdashp		
}
if
show
} for
} bind def
/Mfixdashp {
dup
length
1
gt
1 index
true exch
{
45
eq
and
} forall
and
{
gsave
(--)		
stringwidth pop
(-)
stringwidth pop
sub
2 div
0 rmoveto
dup
length
1 sub
{
(-)
show
}
repeat
grestore
}
if
} bind def	
/Mshoway {
3 index
Mwidthcal
5 1 roll
0 1
4 index length
-1 add
{
2 index
4 index
2 index
get
3 index
add
moveto
4 index
exch get
[
6 index
aload
length
2 add
-1 roll
{
pop
Strform
stringwidth
pop
neg 
exch
add
0 rmoveto
}
exch
kshow
cleartomark
} for
pop	
} bind def
/Mwidthcal {
[
exch
{
Mwidthcal1
}
forall
]
[
exch
dup
Maxlen
-1 add
0 1
3 -1 roll
{
[
exch
2 index
{		 
1 index	
Mget
exch
}
forall		
pop
Maxget
exch
}
for
pop
]
Mreva
} bind def
/Mreva	{
[
exch
aload
length
-1 1
{1 roll}
for
]
} bind def
/Mget	{
1 index
length
-1 add
1 index
ge
{
get
}
{
pop pop
0
}
ifelse
} bind def
/Maxlen	{
[
exch
{
length
}
forall
Maxget
} bind def
/Maxget	{
counttomark
-1 add
1 1
3 -1 roll
{
pop
Mmax
}
for
exch
pop
} bind def
/Mwidthcal1 {
[
exch
{
Strform
stringwidth
pop
}
forall
]
} bind def
/Strform {
/tem (x) def
tem 0
3 -1 roll
put
tem
} bind def
/Mshowa1 {
2 copy
add
4 1 roll
sub
mul
sub
-2 div
} bind def
/MathScale {
Mwidth
Mheight
Mlp
translate
scale
/yscale exch def
/ybias exch def
/xscale exch def
/xbias exch def
/Momatrix
xscale yscale matrix scale
xbias ybias matrix translate
matrix concatmatrix def
/Mgmatrix
matrix currentmatrix
def
} bind def
/Mlp {
3 copy
Mlpfirst
{
Mnodistort
{
Mmin
dup
} if
4 index
2 index
2 index
Mlprun
11 index
11 -1 roll
10 -4 roll
Mlp1
8 index
9 -5 roll
Mlp1
4 -1 roll
and
{ exit } if
3 -1 roll
pop pop
} loop
exch
3 1 roll
7 -3 roll
pop pop pop
} bind def
/Mlpfirst {
3 -1 roll
dup length
2 copy
-2 add
get
aload
pop pop pop
4 -2 roll
-1 add
get
aload
pop pop pop
6 -1 roll
3 -1 roll
5 -1 roll
sub
div
4 1 roll
exch sub
div
} bind def
/Mlprun {
2 copy
4 index
0 get
dup
4 1 roll
Mlprun1
3 copy
8 -2 roll
9 -1 roll
{
3 copy
Mlprun1
3 copy
11 -3 roll
/gt Mlpminmax
8 3 roll
11 -3 roll
/lt Mlpminmax
8 3 roll
} forall
pop pop pop pop
3 1 roll
pop pop
aload pop
5 -1 roll
aload pop
exch
6 -1 roll
Mlprun2
8 2 roll
4 -1 roll
Mlprun2
6 2 roll
3 -1 roll
Mlprun2
4 2 roll
exch
Mlprun2
6 2 roll
} bind def
/Mlprun1 {
aload pop
exch
6 -1 roll
5 -1 roll
mul add
4 -2 roll
mul
3 -1 roll
add
} bind def
/Mlprun2 {
2 copy
add 2 div
3 1 roll
exch sub
} bind def
/Mlpminmax {
cvx
2 index
6 index
2 index
exec
{
7 -3 roll
4 -1 roll
} if
1 index
5 index
3 -1 roll
exec
{
4 1 roll
pop
5 -1 roll
aload
pop pop
4 -1 roll
aload pop
[
8 -2 roll
pop
5 -2 roll
pop
6 -2 roll
pop
5 -1 roll
]
4 1 roll
pop
}
{
pop pop pop
} ifelse
} bind def
/Mlp1 {
5 index
3 index	sub
5 index
2 index mul
1 index
le
1 index
0 le
or
dup
not
{
1 index
3 index	div
.99999 mul
8 -1 roll
pop
7 1 roll
}
if
8 -1 roll
2 div
7 -2 roll
pop sub
5 index
6 -3 roll
pop pop
mul sub
exch
} bind def
/intop 0 def
/inrht 0 def
/inflag 0 def
/outflag 0 def
/xadrht 0 def
/xadlft 0 def
/yadtop 0 def
/yadbot 0 def
/Minner {
outflag
1
eq
{
/outflag 0 def
/intop 0 def
/inrht 0 def
} if		
5 index
gsave
Mtmatrix setmatrix
Mvboxa pop
grestore
3 -1 roll
pop
dup
intop
gt
{
/intop
exch def
}
{ pop }
ifelse
dup
inrht
gt
{
/inrht
exch def
}
{ pop }
ifelse
pop
/inflag
1 def
} bind def
/Mouter {
/xadrht 0 def
/xadlft 0 def
/yadtop 0 def
/yadbot 0 def
inflag
1 eq
{
dup
0 lt
{
dup
intop
mul
neg
/yadtop
exch def
} if
dup
0 gt
{
dup
intop
mul
/yadbot
exch def
}
if
pop
dup
0 lt
{
dup
inrht
mul
neg
/xadrht
exch def
} if
dup
0 gt
{
dup
inrht
mul
/xadlft
exch def
} if
pop
/outflag 1 def
}
{ pop pop}
ifelse
/inflag 0 def
/inrht 0 def
/intop 0 def
} bind def	
/Mboxout {
outflag
1
eq
{
4 -1
roll
xadlft
leadjust
add
sub
4 1 roll
3 -1
roll
yadbot
leadjust
add
sub
3 1
roll
exch
xadrht
leadjust
add
add
exch
yadtop
leadjust
add
add
/outflag 0 def
/xadlft 0 def
/yadbot 0 def
/xadrht 0 def
/yadtop 0 def
} if
} bind def
/leadjust {
(m) stringwidth pop
.5 mul
} bind def
/Mrotcheck {
dup
90
eq
{
yadbot
/yadbot
xadrht
def	
/xadrht
yadtop
def
/yadtop
xadlft
def
/xadlft
exch
def
}
if
dup
cos
1 index
sin
Checkaux
dup
cos
1 index
sin neg
exch
Checkaux
3 1 roll
pop pop
} bind def
/Checkaux {
4 index
exch
4 index
mul
3 1 roll
mul add
4 1 roll
} bind def
/Mboxrot {
Mrot
90 eq
{
brotaux
4 2
roll
}
if
Mrot
180 eq
{
4 2
roll
brotaux
4 2
roll
brotaux
}	
if
Mrot
270 eq
{
4 2
roll
brotaux
}
if
} bind def
/brotaux {
neg
exch
neg
} bind def
/Mabsproc {
0
matrix defaultmatrix
dtransform idtransform
dup mul exch
dup mul
add sqrt
} bind def
/Mabswid {
Mabsproc
setlinewidth	
} bind def
/Mabsdash {
exch
[
exch
{
Mabsproc
}
forall
]
exch
setdash
} bind def
/MBeginOrig { Momatrix concat} bind def
/MEndOrig { Mgmatrix setmatrix} bind def
/sampledsound where
{ pop}
{ /sampledsound {
exch
pop
exch
5 1 roll
mul
4 idiv
mul
2 idiv
exch pop
exch
/Mtempproc exch def
{ Mtempproc pop}
repeat
} bind def
} ifelse
/g { setgray} bind def
/k { setcmykcolor} bind def
/m { moveto} bind def
/p { gsave} bind def
/r { setrgbcolor} bind def
/w { setlinewidth} bind def
/C { curveto} bind def
/F { fill} bind def
/L { lineto} bind def
/P { grestore} bind def
/s { stroke} bind def
/setcmykcolor where
{ pop}
{ /setcmykcolor {
4 1
roll
[
4 1 
roll
]
{
1 index
sub
1
sub neg
dup
0
lt
{
pop
0
}
if
dup
1
gt
{
pop
1
}
if
exch
} forall
pop
setrgbcolor
} bind def
} ifelse
/Mcharproc
{
currentfile
(x)
readhexstring
pop
0 get
exch
div
} bind def
/Mshadeproc
{
dup
3 1
roll
{
dup
Mcharproc
3 1
roll
} repeat
1 eq
{
setgray
}
{
3 eq
{
setrgbcolor
}
{
setcmykcolor
} ifelse
} ifelse
} bind def
/Mrectproc
{
3 index
2 index
moveto
2 index
3 -1
roll
lineto
dup
3 1
roll
lineto
lineto
fill
} bind def
/Mcolorimage
{
7 1
roll
pop
pop
matrix
invertmatrix
concat
2 exch exp
1 sub
3 1 roll
1 1
2 index
{
1 1
4 index
{
dup
1 sub
exch
2 index
dup
1 sub
exch
7 index
9 index
Mshadeproc
Mrectproc
} for
pop
} for
pop pop pop pop
} bind def
/Mimage
{
pop
matrix
invertmatrix
concat
2 exch exp
1 sub
3 1 roll
1 1
2 index
{
1 1
4 index
{
dup
1 sub
exch
2 index
dup
1 sub
exch
7 index
Mcharproc
setgray
Mrectproc
} for
pop
} for
pop pop pop
} bind def

MathPictureStart
/Courier findfont 10  scalefont  setfont
-1.88095 1.90476 -1.05 0.0015 [
[(1)] .02381 0 0 2 Msboxa
[(1.1)] .21429 0 0 2 Msboxa
[(1.2)] .40476 0 0 2 Msboxa
[(1.3)] .59524 0 0 2 Msboxa
[(1.4)] .78571 0 0 2 Msboxa
[(1.5)] .97619 0 0 2 Msboxa
[(800)] .01131 .15 1 0 Msboxa
[(900)] .01131 .3 1 0 Msboxa
[(1000)] .01131 .45 1 0 Msboxa
[(1100)] .01131 .6 1 0 Msboxa
[(1200)] .01131 .75 1 0 Msboxa
[(1300)] .01131 .9 1 0 Msboxa
[ -0.001 -0.001 0 0 ]
[ 1.001 .901 0 0 ]
] MathScale
1 setlinecap
1 setlinejoin
newpath
[ ] 0 setdash
0 g
p
p
.002 w
.02381 0 m
.02381 .00625 L
s
P
[(1)] .02381 0 0 2 Mshowa
p
.002 w
.21429 0 m
.21429 .00625 L
s
P
[(1.1)] .21429 0 0 2 Mshowa
p
.002 w
.40476 0 m
.40476 .00625 L
s
P
[(1.2)] .40476 0 0 2 Mshowa
p
.002 w
.59524 0 m
.59524 .00625 L
s
P
[(1.3)] .59524 0 0 2 Mshowa
p
.002 w
.78571 0 m
.78571 .00625 L
s
P
[(1.4)] .78571 0 0 2 Mshowa
p
.002 w
.97619 0 m
.97619 .00625 L
s
P
[(1.5)] .97619 0 0 2 Mshowa
p
.001 w
.0619 0 m
.0619 .00375 L
s
P
p
.001 w
.1 0 m
.1 .00375 L
s
P
p
.001 w
.1381 0 m
.1381 .00375 L
s
P
p
.001 w
.17619 0 m
.17619 .00375 L
s
P
p
.001 w
.25238 0 m
.25238 .00375 L
s
P
p
.001 w
.29048 0 m
.29048 .00375 L
s
P
p
.001 w
.32857 0 m
.32857 .00375 L
s
P
p
.001 w
.36667 0 m
.36667 .00375 L
s
P
p
.001 w
.44286 0 m
.44286 .00375 L
s
P
p
.001 w
.48095 0 m
.48095 .00375 L
s
P
p
.001 w
.51905 0 m
.51905 .00375 L
s
P
p
.001 w
.55714 0 m
.55714 .00375 L
s
P
p
.001 w
.63333 0 m
.63333 .00375 L
s
P
p
.001 w
.67143 0 m
.67143 .00375 L
s
P
p
.001 w
.70952 0 m
.70952 .00375 L
s
P
p
.001 w
.74762 0 m
.74762 .00375 L
s
P
p
.001 w
.82381 0 m
.82381 .00375 L
s
P
p
.001 w
.8619 0 m
.8619 .00375 L
s
P
p
.001 w
.9 0 m
.9 .00375 L
s
P
p
.001 w
.9381 0 m
.9381 .00375 L
s
P
p
.002 w
0 0 m
1 0 L
s
P
p
.002 w
.02381 .15 m
.03006 .15 L
s
P
[(800)] .01131 .15 1 0 Mshowa
p
.002 w
.02381 .3 m
.03006 .3 L
s
P
[(900)] .01131 .3 1 0 Mshowa
p
.002 w
.02381 .45 m
.03006 .45 L
s
P
[(1000)] .01131 .45 1 0 Mshowa
p
.002 w
.02381 .6 m
.03006 .6 L
s
P
[(1100)] .01131 .6 1 0 Mshowa
p
.002 w
.02381 .75 m
.03006 .75 L
s
P
[(1200)] .01131 .75 1 0 Mshowa
p
.002 w
.02381 .9 m
.03006 .9 L
s
P
[(1300)] .01131 .9 1 0 Mshowa
p
.001 w
.02381 .03 m
.02756 .03 L
s
P
p
.001 w
.02381 .06 m
.02756 .06 L
s
P
p
.001 w
.02381 .09 m
.02756 .09 L
s
P
p
.001 w
.02381 .12 m
.02756 .12 L
s
P
p
.001 w
.02381 .18 m
.02756 .18 L
s
P
p
.001 w
.02381 .21 m
.02756 .21 L
s
P
p
.001 w
.02381 .24 m
.02756 .24 L
s
P
p
.001 w
.02381 .27 m
.02756 .27 L
s
P
p
.001 w
.02381 .33 m
.02756 .33 L
s
P
p
.001 w
.02381 .36 m
.02756 .36 L
s
P
p
.001 w
.02381 .39 m
.02756 .39 L
s
P
p
.001 w
.02381 .42 m
.02756 .42 L
s
P
p
.001 w
.02381 .48 m
.02756 .48 L
s
P
p
.001 w
.02381 .51 m
.02756 .51 L
s
P
p
.001 w
.02381 .54 m
.02756 .54 L
s
P
p
.001 w
.02381 .57 m
.02756 .57 L
s
P
p
.001 w
.02381 .63 m
.02756 .63 L
s
P
p
.001 w
.02381 .66 m
.02756 .66 L
s
P
p
.001 w
.02381 .69 m
.02756 .69 L
s
P
p
.001 w
.02381 .72 m
.02756 .72 L
s
P
p
.001 w
.02381 .78 m
.02756 .78 L
s
P
p
.001 w
.02381 .81 m
.02756 .81 L
s
P
p
.001 w
.02381 .84 m
.02756 .84 L
s
P
p
.001 w
.02381 .87 m
.02756 .87 L
s
P
p
.002 w
.02381 0 m
.02381 .9 L
s
P
P
0 0 m
1 0 L
1 .9 L
0 .9 L
closepath
clip
newpath
p
p
p
.001 w
.02381 .88 m
.06349 .80325 L
.10317 .73329 L
.14286 .66937 L
.18254 .61086 L
.22222 .5572 L
.2619 .50788 L
.30159 .46248 L
.34127 .4206 L
.38095 .38191 L
.42063 .34611 L
.46032 .31293 L
.5 .28213 L
.53968 .2535 L
.57937 .22684 L
.61905 .20199 L
.65873 .17879 L
.69841 .1571 L
.7381 .13681 L
.77778 .11779 L
.81746 .09995 L
.85714 .0832 L
.89683 .06745 L
.93651 .05262 L
.97619 .03865 L
s
P
P
p
p
.001 w
.02381 .89389 m
.06349 .81657 L
.10317 .74608 L
.14286 .68166 L
.18254 .62267 L
.22222 .56855 L
.2619 .5188 L
.30159 .47298 L
.34127 .43072 L
.38095 .39165 L
.42063 .3555 L
.46032 .32198 L
.5 .29086 L
.53968 .26191 L
.57937 .23496 L
.61905 .20983 L
.65873 .18636 L
.69841 .16442 L
.7381 .14388 L
.77778 .12463 L
.81746 .10657 L
.85714 .0896 L
.89683 .07364 L
.93651 .05862 L
.97619 .04446 L
s
P
P
P
MathPictureEnd
end
EndOfTheIncludedPostscriptMagicCookie

\closepsdump

\psdump{negpart.eps}

/Mathdict 150 dict def
Mathdict begin
/Mlmarg		0 72 mul def
/Mrmarg		0 72 mul def
/Mbmarg		0 72 mul def
/Mtmarg		0 72 mul def
/Mwidth		4 72 mul def
/Mheight	4 72 mul def
/Mtransform	{  } bind def
/Mnodistort	true def
/Mfixwid	true	def
/Mfixdash	false def
/Mrot		0	def
/Mpstart {
MathPictureStart
} bind def
/Mpend {
MathPictureEnd
} bind def
/Mscale {
0 1 0 1
5 -1 roll
MathScale
} bind def
/ISOLatin1Encoding dup where
{ pop pop }
{
[
/.notdef /.notdef /.notdef /.notdef /.notdef /.notdef /.notdef /.notdef
/.notdef /.notdef /.notdef /.notdef /.notdef /.notdef /.notdef /.notdef
/.notdef /.notdef /.notdef /.notdef /.notdef /.notdef /.notdef /.notdef
/.notdef /.notdef /.notdef /.notdef /.notdef /.notdef /.notdef /.notdef
/space /exclam /quotedbl /numbersign /dollar /percent /ampersand /quoteright
/parenleft /parenright /asterisk /plus /comma /minus /period /slash
/zero /one /two /three /four /five /six /seven
/eight /nine /colon /semicolon /less /equal /greater /question
/at /A /B /C /D /E /F /G
/H /I /J /K /L /M /N /O
/P /Q /R /S /T /U /V /W
/X /Y /Z /bracketleft /backslash /bracketright /asciicircum /underscore
/quoteleft /a /b /c /d /e /f /g
/h /i /j /k /l /m /n /o
/p /q /r /s /t /u /v /w
/x /y /z /braceleft /bar /braceright /asciitilde /.notdef
/.notdef /.notdef /.notdef /.notdef /.notdef /.notdef /.notdef /.notdef
/.notdef /.notdef /.notdef /.notdef /.notdef /.notdef /.notdef /.notdef
/dotlessi /grave /acute /circumflex /tilde /macron /breve /dotaccent
/dieresis /.notdef /ring /cedilla /.notdef /hungarumlaut /ogonek /caron
/space /exclamdown /cent /sterling /currency /yen /brokenbar /section
/dieresis /copyright /ordfeminine /guillemotleft
/logicalnot /hyphen /registered /macron
/degree /plusminus /twosuperior /threesuperior
/acute /mu /paragraph /periodcentered
/cedilla /onesuperior /ordmasculine /guillemotright
/onequarter /onehalf /threequarters /questiondown
/Agrave /Aacute /Acircumflex /Atilde /Adieresis /Aring /AE /Ccedilla
/Egrave /Eacute /Ecircumflex /Edieresis /Igrave /Iacute /Icircumflex /Idieresis
/Eth /Ntilde /Ograve /Oacute /Ocircumflex /Otilde /Odieresis /multiply
/Oslash /Ugrave /Uacute /Ucircumflex /Udieresis /Yacute /Thorn /germandbls
/agrave /aacute /acircumflex /atilde /adieresis /aring /ae /ccedilla
/egrave /eacute /ecircumflex /edieresis /igrave /iacute /icircumflex /idieresis
/eth /ntilde /ograve /oacute /ocircumflex /otilde /odieresis /divide
/oslash /ugrave /uacute /ucircumflex /udieresis /yacute /thorn /ydieresis
] def
} ifelse
/MFontDict 50 dict def
/MStrCat
{
exch
dup length
2 index length add
string
dup 3 1 roll
copy
length
exch dup
4 2 roll exch
putinterval
} def
/MCreateEncoding
{
1 index
255 string cvs
(-) MStrCat
1 index MStrCat
cvn exch
(Encoding) MStrCat
cvn dup where
{
exch get
}
{
pop
StandardEncoding
} ifelse
3 1 roll
dup MFontDict exch known not
{
1 index findfont
dup length dict
begin
{1 index /FID ne
{def}
{pop pop}
ifelse} forall
/Encoding 3 index
def
currentdict
end
1 index exch definefont pop
MFontDict 1 index
null put
}
if
exch pop
exch pop
} def
/ISOLatin1 { (ISOLatin1) MCreateEncoding } def
/ISO8859 { (ISOLatin1) MCreateEncoding } def
/Mcopyfont {
dup
maxlength
dict
exch
{
1 index
/FID
eq
{
pop pop
}
{
2 index
3 1 roll
put
}
ifelse
}
forall
} def
/Plain	/Courier findfont Mcopyfont definefont pop
/Bold	/Courier-Bold findfont Mcopyfont definefont pop
/Italic /Courier-Oblique findfont Mcopyfont definefont pop
/MathPictureStart {
gsave
Mtransform
Mlmarg
Mbmarg
translate
Mwidth
Mlmarg Mrmarg add
sub
/Mwidth exch def
Mheight
Mbmarg Mtmarg add
sub
/Mheight exch def
/Mtmatrix
matrix currentmatrix
def
/Mgmatrix
matrix currentmatrix
def
} bind def
/MathPictureEnd {
grestore
} bind def
/MFill {
0 0 		moveto
Mwidth 0 	lineto
Mwidth Mheight 	lineto
0 Mheight 	lineto
fill
} bind def
/MPlotRegion {
3 index
Mwidth mul
2 index
Mheight mul
translate
exch sub
Mheight mul
/Mheight
exch def
exch sub
Mwidth mul
/Mwidth
exch def
} bind def
/MathSubStart {
Momatrix
Mgmatrix Mtmatrix
Mwidth Mheight
7 -2 roll
moveto
Mtmatrix setmatrix
currentpoint
Mgmatrix setmatrix
9 -2 roll
moveto
Mtmatrix setmatrix
currentpoint
2 copy translate
/Mtmatrix matrix currentmatrix def
3 -1 roll
exch sub
/Mheight exch def
sub
/Mwidth exch def
} bind def
/MathSubEnd {
/Mheight exch def
/Mwidth exch def
/Mtmatrix exch def
dup setmatrix
/Mgmatrix exch def
/Momatrix exch def
} bind def
/Mdot {
moveto
0 0 rlineto
stroke
} bind def
/Mtetra {
moveto
lineto
lineto
lineto
fill
} bind def
/Metetra {
moveto
lineto
lineto
lineto
closepath
gsave
fill
grestore
0 setgray
stroke
} bind def
/Mistroke {
flattenpath
0 0 0
{
4 2 roll
pop pop
}
{
4 -1 roll
2 index
sub dup mul
4 -1 roll
2 index
sub dup mul
add sqrt
4 -1 roll
add
3 1 roll
}
{
stop
}
{
stop
}
pathforall
pop pop
currentpoint
stroke
moveto
currentdash
3 -1 roll
add
setdash
} bind def
/Mfstroke {
stroke
currentdash
pop 0
setdash
} bind def
/Mrotsboxa {
gsave
dup
/Mrot
exch def
Mrotcheck
Mtmatrix
dup
setmatrix
7 1 roll
4 index
4 index
translate
rotate
3 index
-1 mul
3 index
-1 mul
translate
/Mtmatrix
matrix
currentmatrix
def
grestore
Msboxa
3  -1 roll
/Mtmatrix
exch def
/Mrot
0 def
} bind def
/Msboxa {
newpath
5 -1 roll
Mvboxa
pop
Mboxout
6 -1 roll
5 -1 roll
4 -1 roll
Msboxa1
5 -3 roll
Msboxa1
Mboxrot
[
7 -2 roll
2 copy
[
3 1 roll
10 -1 roll
9 -1 roll
]
6 1 roll
5 -2 roll
]
} bind def
/Msboxa1 {
sub
2 div
dup
2 index
1 add
mul
3 -1 roll
-1 add
3 -1 roll
mul
} bind def
/Mvboxa	{
Mfixwid
{
Mvboxa1
}
{
dup
Mwidthcal
0 exch
{
add
}
forall
exch
Mvboxa1
4 index
7 -1 roll
add
4 -1 roll
pop
3 1 roll
}
ifelse
} bind def
/Mvboxa1 {
gsave
newpath
[ true
3 -1 roll
{
Mbbox
5 -1 roll
{
0
5 1 roll
}
{
7 -1 roll
exch sub
(m) stringwidth pop
.3 mul
sub
7 1 roll
6 -1 roll
4 -1 roll
Mmin
3 -1 roll
5 index
add
5 -1 roll
4 -1 roll
Mmax
4 -1 roll
}
ifelse
false
}
forall
{ stop } if
counttomark
1 add
4 roll
]
grestore
} bind def
/Mbbox {
1 dict begin
0 0 moveto
/temp (T) def
{	gsave
currentpoint newpath moveto
temp 0 3 -1 roll put temp
false charpath flattenpath currentpoint
pathbbox
grestore moveto lineto moveto} forall
pathbbox
newpath
end
} bind def
/Mmin {
2 copy
gt
{ exch } if
pop
} bind def
/Mmax {
2 copy
lt
{ exch } if
pop
} bind def
/Mrotshowa {
dup
/Mrot
exch def
Mrotcheck
Mtmatrix
dup
setmatrix
7 1 roll
4 index
4 index
translate
rotate
3 index
-1 mul
3 index
-1 mul
translate
/Mtmatrix
matrix
currentmatrix
def
Mgmatrix setmatrix
Mshowa
/Mtmatrix
exch def
/Mrot 0 def
} bind def
/Mshowa {
4 -2 roll
moveto
2 index
Mtmatrix setmatrix
Mvboxa
7 1 roll
Mboxout
6 -1 roll
5 -1 roll
4 -1 roll
Mshowa1
4 1 roll
Mshowa1
rmoveto
currentpoint
Mfixwid
{
Mshowax
}
{
Mshoway
}
ifelse
pop pop pop pop
Mgmatrix setmatrix
} bind def
/Mshowax {	
0 1
4 index length
-1 add
{
2 index
4 index
2 index
get
3 index
add
moveto
4 index
exch get
Mfixdash
{
Mfixdashp		
}
if
show
} for
} bind def
/Mfixdashp {
dup
length
1
gt
1 index
true exch
{
45
eq
and
} forall
and
{
gsave
(--)		
stringwidth pop
(-)
stringwidth pop
sub
2 div
0 rmoveto
dup
length
1 sub
{
(-)
show
}
repeat
grestore
}
if
} bind def	
/Mshoway {
3 index
Mwidthcal
5 1 roll
0 1
4 index length
-1 add
{
2 index
4 index
2 index
get
3 index
add
moveto
4 index
exch get
[
6 index
aload
length
2 add
-1 roll
{
pop
Strform
stringwidth
pop
neg 
exch
add
0 rmoveto
}
exch
kshow
cleartomark
} for
pop	
} bind def
/Mwidthcal {
[
exch
{
Mwidthcal1
}
forall
]
[
exch
dup
Maxlen
-1 add
0 1
3 -1 roll
{
[
exch
2 index
{		 
1 index	
Mget
exch
}
forall		
pop
Maxget
exch
}
for
pop
]
Mreva
} bind def
/Mreva	{
[
exch
aload
length
-1 1
{1 roll}
for
]
} bind def
/Mget	{
1 index
length
-1 add
1 index
ge
{
get
}
{
pop pop
0
}
ifelse
} bind def
/Maxlen	{
[
exch
{
length
}
forall
Maxget
} bind def
/Maxget	{
counttomark
-1 add
1 1
3 -1 roll
{
pop
Mmax
}
for
exch
pop
} bind def
/Mwidthcal1 {
[
exch
{
Strform
stringwidth
pop
}
forall
]
} bind def
/Strform {
/tem (x) def
tem 0
3 -1 roll
put
tem
} bind def
/Mshowa1 {
2 copy
add
4 1 roll
sub
mul
sub
-2 div
} bind def
/MathScale {
Mwidth
Mheight
Mlp
translate
scale
/yscale exch def
/ybias exch def
/xscale exch def
/xbias exch def
/Momatrix
xscale yscale matrix scale
xbias ybias matrix translate
matrix concatmatrix def
/Mgmatrix
matrix currentmatrix
def
} bind def
/Mlp {
3 copy
Mlpfirst
{
Mnodistort
{
Mmin
dup
} if
4 index
2 index
2 index
Mlprun
11 index
11 -1 roll
10 -4 roll
Mlp1
8 index
9 -5 roll
Mlp1
4 -1 roll
and
{ exit } if
3 -1 roll
pop pop
} loop
exch
3 1 roll
7 -3 roll
pop pop pop
} bind def
/Mlpfirst {
3 -1 roll
dup length
2 copy
-2 add
get
aload
pop pop pop
4 -2 roll
-1 add
get
aload
pop pop pop
6 -1 roll
3 -1 roll
5 -1 roll
sub
div
4 1 roll
exch sub
div
} bind def
/Mlprun {
2 copy
4 index
0 get
dup
4 1 roll
Mlprun1
3 copy
8 -2 roll
9 -1 roll
{
3 copy
Mlprun1
3 copy
11 -3 roll
/gt Mlpminmax
8 3 roll
11 -3 roll
/lt Mlpminmax
8 3 roll
} forall
pop pop pop pop
3 1 roll
pop pop
aload pop
5 -1 roll
aload pop
exch
6 -1 roll
Mlprun2
8 2 roll
4 -1 roll
Mlprun2
6 2 roll
3 -1 roll
Mlprun2
4 2 roll
exch
Mlprun2
6 2 roll
} bind def
/Mlprun1 {
aload pop
exch
6 -1 roll
5 -1 roll
mul add
4 -2 roll
mul
3 -1 roll
add
} bind def
/Mlprun2 {
2 copy
add 2 div
3 1 roll
exch sub
} bind def
/Mlpminmax {
cvx
2 index
6 index
2 index
exec
{
7 -3 roll
4 -1 roll
} if
1 index
5 index
3 -1 roll
exec
{
4 1 roll
pop
5 -1 roll
aload
pop pop
4 -1 roll
aload pop
[
8 -2 roll
pop
5 -2 roll
pop
6 -2 roll
pop
5 -1 roll
]
4 1 roll
pop
}
{
pop pop pop
} ifelse
} bind def
/Mlp1 {
5 index
3 index	sub
5 index
2 index mul
1 index
le
1 index
0 le
or
dup
not
{
1 index
3 index	div
.99999 mul
8 -1 roll
pop
7 1 roll
}
if
8 -1 roll
2 div
7 -2 roll
pop sub
5 index
6 -3 roll
pop pop
mul sub
exch
} bind def
/intop 0 def
/inrht 0 def
/inflag 0 def
/outflag 0 def
/xadrht 0 def
/xadlft 0 def
/yadtop 0 def
/yadbot 0 def
/Minner {
outflag
1
eq
{
/outflag 0 def
/intop 0 def
/inrht 0 def
} if		
5 index
gsave
Mtmatrix setmatrix
Mvboxa pop
grestore
3 -1 roll
pop
dup
intop
gt
{
/intop
exch def
}
{ pop }
ifelse
dup
inrht
gt
{
/inrht
exch def
}
{ pop }
ifelse
pop
/inflag
1 def
} bind def
/Mouter {
/xadrht 0 def
/xadlft 0 def
/yadtop 0 def
/yadbot 0 def
inflag
1 eq
{
dup
0 lt
{
dup
intop
mul
neg
/yadtop
exch def
} if
dup
0 gt
{
dup
intop
mul
/yadbot
exch def
}
if
pop
dup
0 lt
{
dup
inrht
mul
neg
/xadrht
exch def
} if
dup
0 gt
{
dup
inrht
mul
/xadlft
exch def
} if
pop
/outflag 1 def
}
{ pop pop}
ifelse
/inflag 0 def
/inrht 0 def
/intop 0 def
} bind def	
/Mboxout {
outflag
1
eq
{
4 -1
roll
xadlft
leadjust
add
sub
4 1 roll
3 -1
roll
yadbot
leadjust
add
sub
3 1
roll
exch
xadrht
leadjust
add
add
exch
yadtop
leadjust
add
add
/outflag 0 def
/xadlft 0 def
/yadbot 0 def
/xadrht 0 def
/yadtop 0 def
} if
} bind def
/leadjust {
(m) stringwidth pop
.5 mul
} bind def
/Mrotcheck {
dup
90
eq
{
yadbot
/yadbot
xadrht
def	
/xadrht
yadtop
def
/yadtop
xadlft
def
/xadlft
exch
def
}
if
dup
cos
1 index
sin
Checkaux
dup
cos
1 index
sin neg
exch
Checkaux
3 1 roll
pop pop
} bind def
/Checkaux {
4 index
exch
4 index
mul
3 1 roll
mul add
4 1 roll
} bind def
/Mboxrot {
Mrot
90 eq
{
brotaux
4 2
roll
}
if
Mrot
180 eq
{
4 2
roll
brotaux
4 2
roll
brotaux
}	
if
Mrot
270 eq
{
4 2
roll
brotaux
}
if
} bind def
/brotaux {
neg
exch
neg
} bind def
/Mabsproc {
0
matrix defaultmatrix
dtransform idtransform
dup mul exch
dup mul
add sqrt
} bind def
/Mabswid {
Mabsproc
setlinewidth	
} bind def
/Mabsdash {
exch
[
exch
{
Mabsproc
}
forall
]
exch
setdash
} bind def
/MBeginOrig { Momatrix concat} bind def
/MEndOrig { Mgmatrix setmatrix} bind def
/sampledsound where
{ pop}
{ /sampledsound {
exch
pop
exch
5 1 roll
mul
4 idiv
mul
2 idiv
exch pop
exch
/Mtempproc exch def
{ Mtempproc pop}
repeat
} bind def
} ifelse
/g { setgray} bind def
/k { setcmykcolor} bind def
/m { moveto} bind def
/p { gsave} bind def
/r { setrgbcolor} bind def
/w { setlinewidth} bind def
/C { curveto} bind def
/F { fill} bind def
/L { lineto} bind def
/P { grestore} bind def
/s { stroke} bind def
/setcmykcolor where
{ pop}
{ /setcmykcolor {
4 1
roll
[
4 1 
roll
]
{
1 index
sub
1
sub neg
dup
0
lt
{
pop
0
}
if
dup
1
gt
{
pop
1
}
if
exch
} forall
pop
setrgbcolor
} bind def
} ifelse
/Mcharproc
{
currentfile
(x)
readhexstring
pop
0 get
exch
div
} bind def
/Mshadeproc
{
dup
3 1
roll
{
dup
Mcharproc
3 1
roll
} repeat
1 eq
{
setgray
}
{
3 eq
{
setrgbcolor
}
{
setcmykcolor
} ifelse
} ifelse
} bind def
/Mrectproc
{
3 index
2 index
moveto
2 index
3 -1
roll
lineto
dup
3 1
roll
lineto
lineto
fill
} bind def
/Mcolorimage
{
7 1
roll
pop
pop
matrix
invertmatrix
concat
2 exch exp
1 sub
3 1 roll
1 1
2 index
{
1 1
4 index
{
dup
1 sub
exch
2 index
dup
1 sub
exch
7 index
9 index
Mshadeproc
Mrectproc
} for
pop
} for
pop pop pop pop
} bind def
/Mimage
{
pop
matrix
invertmatrix
concat
2 exch exp
1 sub
3 1 roll
1 1
2 index
{
1 1
4 index
{
dup
1 sub
exch
2 index
dup
1 sub
exch
7 index
Mcharproc
setgray
Mrectproc
} for
pop
} for
pop pop pop
} bind def

MathPictureStart
/Courier findfont 10  scalefont  setfont
0.0238095 0.00952381 0 12.3607 [
[(0)] .02381 0 0 2 Msboxa
[(20)] .21429 0 0 2 Msboxa
[(40)] .40476 0 0 2 Msboxa
[(60)] .59524 0 0 2 Msboxa
[(80)] .78571 0 0 2 Msboxa
[(100)] .97619 0 0 2 Msboxa
[(0.01)] .01131 .12361 1 0 Msboxa
[(0.02)] .01131 .24721 1 0 Msboxa
[(0.03)] .01131 .37082 1 0 Msboxa
[(0.04)] .01131 .49443 1 0 Msboxa
[(0.05)] .01131 .61803 1 0 Msboxa
[ -0.001 -0.001 0 0 ]
[ 1.001 .61903 0 0 ]
] MathScale
1 setlinecap
1 setlinejoin
newpath
[ ] 0 setdash
0 g
p
p
.002 w
.02381 0 m
.02381 .00625 L
s
P
[(0)] .02381 0 0 2 Mshowa
p
.002 w
.21429 0 m
.21429 .00625 L
s
P
[(20)] .21429 0 0 2 Mshowa
p
.002 w
.40476 0 m
.40476 .00625 L
s
P
[(40)] .40476 0 0 2 Mshowa
p
.002 w
.59524 0 m
.59524 .00625 L
s
P
[(60)] .59524 0 0 2 Mshowa
p
.002 w
.78571 0 m
.78571 .00625 L
s
P
[(80)] .78571 0 0 2 Mshowa
p
.002 w
.97619 0 m
.97619 .00625 L
s
P
[(100)] .97619 0 0 2 Mshowa
p
.001 w
.0619 0 m
.0619 .00375 L
s
P
p
.001 w
.1 0 m
.1 .00375 L
s
P
p
.001 w
.1381 0 m
.1381 .00375 L
s
P
p
.001 w
.17619 0 m
.17619 .00375 L
s
P
p
.001 w
.25238 0 m
.25238 .00375 L
s
P
p
.001 w
.29048 0 m
.29048 .00375 L
s
P
p
.001 w
.32857 0 m
.32857 .00375 L
s
P
p
.001 w
.36667 0 m
.36667 .00375 L
s
P
p
.001 w
.44286 0 m
.44286 .00375 L
s
P
p
.001 w
.48095 0 m
.48095 .00375 L
s
P
p
.001 w
.51905 0 m
.51905 .00375 L
s
P
p
.001 w
.55714 0 m
.55714 .00375 L
s
P
p
.001 w
.63333 0 m
.63333 .00375 L
s
P
p
.001 w
.67143 0 m
.67143 .00375 L
s
P
p
.001 w
.70952 0 m
.70952 .00375 L
s
P
p
.001 w
.74762 0 m
.74762 .00375 L
s
P
p
.001 w
.82381 0 m
.82381 .00375 L
s
P
p
.001 w
.8619 0 m
.8619 .00375 L
s
P
p
.001 w
.9 0 m
.9 .00375 L
s
P
p
.001 w
.9381 0 m
.9381 .00375 L
s
P
p
.002 w
0 0 m
1 0 L
s
P
p
.002 w
.02381 .12361 m
.03006 .12361 L
s
P
[(0.01)] .01131 .12361 1 0 Mshowa
p
.002 w
.02381 .24721 m
.03006 .24721 L
s
P
[(0.02)] .01131 .24721 1 0 Mshowa
p
.002 w
.02381 .37082 m
.03006 .37082 L
s
P
[(0.03)] .01131 .37082 1 0 Mshowa
p
.002 w
.02381 .49443 m
.03006 .49443 L
s
P
[(0.04)] .01131 .49443 1 0 Mshowa
p
.002 w
.02381 .61803 m
.03006 .61803 L
s
P
[(0.05)] .01131 .61803 1 0 Mshowa
p
.001 w
.02381 .02472 m
.02756 .02472 L
s
P
p
.001 w
.02381 .04944 m
.02756 .04944 L
s
P
p
.001 w
.02381 .07416 m
.02756 .07416 L
s
P
p
.001 w
.02381 .09889 m
.02756 .09889 L
s
P
p
.001 w
.02381 .14833 m
.02756 .14833 L
s
P
p
.001 w
.02381 .17305 m
.02756 .17305 L
s
P
p
.001 w
.02381 .19777 m
.02756 .19777 L
s
P
p
.001 w
.02381 .22249 m
.02756 .22249 L
s
P
p
.001 w
.02381 .27193 m
.02756 .27193 L
s
P
p
.001 w
.02381 .29666 m
.02756 .29666 L
s
P
p
.001 w
.02381 .32138 m
.02756 .32138 L
s
P
p
.001 w
.02381 .3461 m
.02756 .3461 L
s
P
p
.001 w
.02381 .39554 m
.02756 .39554 L
s
P
p
.001 w
.02381 .42026 m
.02756 .42026 L
s
P
p
.001 w
.02381 .44498 m
.02756 .44498 L
s
P
p
.001 w
.02381 .46971 m
.02756 .46971 L
s
P
p
.001 w
.02381 .51915 m
.02756 .51915 L
s
P
p
.001 w
.02381 .54387 m
.02756 .54387 L
s
P
p
.001 w
.02381 .56859 m
.02756 .56859 L
s
P
p
.001 w
.02381 .59331 m
.02756 .59331 L
s
P
p
.002 w
.02381 0 m
.02381 .61803 L
s
P
P
0 0 m
1 0 L
1 .61803 L
0 .61803 L
closepath
clip
newpath
p
p
.004 w
.03333 .59775 m
.03571 .48858 L
.0369 .44429 L
.0375 .4243 L
.0378 .4148 L
.0381 .40846 L
.03869 .41895 L
.03929 .4288 L
.04048 .44666 L
.04286 .47573 L
.04524 .49718 L
.04643 .5055 L
.04762 .51244 L
.04881 .51816 L
.05 .52278 L
.05119 .52643 L
.05238 .52923 L
.05298 .53034 L
.05357 .53128 L
.05417 .53205 L
.05476 .53266 L
.05536 .53313 L
.05565 .53332 L
.05595 .53347 L
.05625 .53359 L
.05655 .53367 L
.05685 .53373 L
.05714 .53376 L
.05744 .53376 L
.05774 .53374 L
.05804 .53368 L
.05833 .53361 L
.05893 .53338 L
.05952 .53307 L
.06071 .53218 L
.0619 .53101 L
.06429 .52792 L
.06667 .52406 L
.07143 .5148 L
.08095 .49347 L
.09048 .47169 L
.1 .45111 L
.10952 .4322 L
.11905 .415 L
.12857 .39941 L
.1381 .38526 L
.14762 .37239 L
.15714 .36063 L
.16667 .34985 L
.17619 .33994 L
Mistroke
.18571 .33079 L
.19524 .32231 L
.20476 .31443 L
.21429 .30709 L
.22381 .30022 L
.23333 .29379 L
.24286 .28774 L
.25238 .28205 L
.2619 .27667 L
.27143 .27158 L
.28095 .26677 L
.29048 .26219 L
.3 .25784 L
.30952 .2537 L
.31905 .24975 L
.32857 .24597 L
.3381 .24236 L
.34762 .23891 L
.35714 .23559 L
.36667 .23241 L
.37619 .22936 L
.38571 .22642 L
.39524 .22359 L
.40476 .22086 L
.41429 .21823 L
.42381 .21569 L
.43333 .21324 L
.44286 .21087 L
.45238 .20857 L
.4619 .20635 L
.47143 .2042 L
.48095 .20211 L
.49048 .20009 L
.5 .19812 L
.50952 .19621 L
.51905 .19436 L
.52857 .19256 L
.5381 .1908 L
.54762 .18909 L
.55714 .18743 L
.56667 .18581 L
.57619 .18424 L
.58571 .1827 L
.59524 .1812 L
.60476 .17973 L
.61429 .1783 L
.62381 .17691 L
.63333 .17554 L
.64286 .17421 L
.65238 .17291 L
Mistroke
.6619 .17163 L
.67143 .17038 L
.68095 .16916 L
.69048 .16797 L
.7 .1668 L
.70952 .16566 L
.71905 .16453 L
.72857 .16343 L
.7381 .16236 L
.74762 .1613 L
.75714 .16026 L
.76667 .15925 L
.77619 .15825 L
.78571 .15727 L
.79524 .15631 L
.80476 .15536 L
.81429 .15444 L
.82381 .15353 L
.83333 .15263 L
.84286 .15175 L
.85238 .15089 L
.8619 .15004 L
.87143 .1492 L
.88095 .14838 L
.89048 .14757 L
.9 .14678 L
.90952 .146 L
.91905 .14523 L
.92857 .14447 L
.9381 .14372 L
.94762 .14299 L
.95714 .14226 L
.96667 .14155 L
.97619 .14085 L
Mfstroke
P
P
MathPictureEnd
end
EndOfTheIncludedPostscriptMagicCookie

\closepsdump

\psdump{shadr1234.eps}

/Mathdict 150 dict def
Mathdict begin
/Mlmarg		0 72 mul def
/Mrmarg		0 72 mul def
/Mbmarg		0 72 mul def
/Mtmarg		0 72 mul def
/Mwidth		4 72 mul def
/Mheight	4 72 mul def
/Mtransform	{  } bind def
/Mnodistort	true def
/Mfixwid	true	def
/Mfixdash	false def
/Mrot		0	def
/Mpstart {
MathPictureStart
} bind def
/Mpend {
MathPictureEnd
} bind def
/Mscale {
0 1 0 1
5 -1 roll
MathScale
} bind def
/ISOLatin1Encoding dup where
{ pop pop }
{
[
/.notdef /.notdef /.notdef /.notdef /.notdef /.notdef /.notdef /.notdef
/.notdef /.notdef /.notdef /.notdef /.notdef /.notdef /.notdef /.notdef
/.notdef /.notdef /.notdef /.notdef /.notdef /.notdef /.notdef /.notdef
/.notdef /.notdef /.notdef /.notdef /.notdef /.notdef /.notdef /.notdef
/space /exclam /quotedbl /numbersign /dollar /percent /ampersand /quoteright
/parenleft /parenright /asterisk /plus /comma /minus /period /slash
/zero /one /two /three /four /five /six /seven
/eight /nine /colon /semicolon /less /equal /greater /question
/at /A /B /C /D /E /F /G
/H /I /J /K /L /M /N /O
/P /Q /R /S /T /U /V /W
/X /Y /Z /bracketleft /backslash /bracketright /asciicircum /underscore
/quoteleft /a /b /c /d /e /f /g
/h /i /j /k /l /m /n /o
/p /q /r /s /t /u /v /w
/x /y /z /braceleft /bar /braceright /asciitilde /.notdef
/.notdef /.notdef /.notdef /.notdef /.notdef /.notdef /.notdef /.notdef
/.notdef /.notdef /.notdef /.notdef /.notdef /.notdef /.notdef /.notdef
/dotlessi /grave /acute /circumflex /tilde /macron /breve /dotaccent
/dieresis /.notdef /ring /cedilla /.notdef /hungarumlaut /ogonek /caron
/space /exclamdown /cent /sterling /currency /yen /brokenbar /section
/dieresis /copyright /ordfeminine /guillemotleft
/logicalnot /hyphen /registered /macron
/degree /plusminus /twosuperior /threesuperior
/acute /mu /paragraph /periodcentered
/cedilla /onesuperior /ordmasculine /guillemotright
/onequarter /onehalf /threequarters /questiondown
/Agrave /Aacute /Acircumflex /Atilde /Adieresis /Aring /AE /Ccedilla
/Egrave /Eacute /Ecircumflex /Edieresis /Igrave /Iacute /Icircumflex /Idieresis
/Eth /Ntilde /Ograve /Oacute /Ocircumflex /Otilde /Odieresis /multiply
/Oslash /Ugrave /Uacute /Ucircumflex /Udieresis /Yacute /Thorn /germandbls
/agrave /aacute /acircumflex /atilde /adieresis /aring /ae /ccedilla
/egrave /eacute /ecircumflex /edieresis /igrave /iacute /icircumflex /idieresis
/eth /ntilde /ograve /oacute /ocircumflex /otilde /odieresis /divide
/oslash /ugrave /uacute /ucircumflex /udieresis /yacute /thorn /ydieresis
] def
} ifelse
/MFontDict 50 dict def
/MStrCat
{
exch
dup length
2 index length add
string
dup 3 1 roll
copy
length
exch dup
4 2 roll exch
putinterval
} def
/MCreateEncoding
{
1 index
255 string cvs
(-) MStrCat
1 index MStrCat
cvn exch
(Encoding) MStrCat
cvn dup where
{
exch get
}
{
pop
StandardEncoding
} ifelse
3 1 roll
dup MFontDict exch known not
{
1 index findfont
dup length dict
begin
{1 index /FID ne
{def}
{pop pop}
ifelse} forall
/Encoding 3 index
def
currentdict
end
1 index exch definefont pop
MFontDict 1 index
null put
}
if
exch pop
exch pop
} def
/ISOLatin1 { (ISOLatin1) MCreateEncoding } def
/ISO8859 { (ISOLatin1) MCreateEncoding } def
/Mcopyfont {
dup
maxlength
dict
exch
{
1 index
/FID
eq
{
pop pop
}
{
2 index
3 1 roll
put
}
ifelse
}
forall
} def
/Plain	/Courier findfont Mcopyfont definefont pop
/Bold	/Courier-Bold findfont Mcopyfont definefont pop
/Italic /Courier-Oblique findfont Mcopyfont definefont pop
/MathPictureStart {
gsave
Mtransform
Mlmarg
Mbmarg
translate
Mwidth
Mlmarg Mrmarg add
sub
/Mwidth exch def
Mheight
Mbmarg Mtmarg add
sub
/Mheight exch def
/Mtmatrix
matrix currentmatrix
def
/Mgmatrix
matrix currentmatrix
def
} bind def
/MathPictureEnd {
grestore
} bind def
/MFill {
0 0 		moveto
Mwidth 0 	lineto
Mwidth Mheight 	lineto
0 Mheight 	lineto
fill
} bind def
/MPlotRegion {
3 index
Mwidth mul
2 index
Mheight mul
translate
exch sub
Mheight mul
/Mheight
exch def
exch sub
Mwidth mul
/Mwidth
exch def
} bind def
/MathSubStart {
Momatrix
Mgmatrix Mtmatrix
Mwidth Mheight
7 -2 roll
moveto
Mtmatrix setmatrix
currentpoint
Mgmatrix setmatrix
9 -2 roll
moveto
Mtmatrix setmatrix
currentpoint
2 copy translate
/Mtmatrix matrix currentmatrix def
3 -1 roll
exch sub
/Mheight exch def
sub
/Mwidth exch def
} bind def
/MathSubEnd {
/Mheight exch def
/Mwidth exch def
/Mtmatrix exch def
dup setmatrix
/Mgmatrix exch def
/Momatrix exch def
} bind def
/Mdot {
moveto
0 0 rlineto
stroke
} bind def
/Mtetra {
moveto
lineto
lineto
lineto
fill
} bind def
/Metetra {
moveto
lineto
lineto
lineto
closepath
gsave
fill
grestore
0 setgray
stroke
} bind def
/Mistroke {
flattenpath
0 0 0
{
4 2 roll
pop pop
}
{
4 -1 roll
2 index
sub dup mul
4 -1 roll
2 index
sub dup mul
add sqrt
4 -1 roll
add
3 1 roll
}
{
stop
}
{
stop
}
pathforall
pop pop
currentpoint
stroke
moveto
currentdash
3 -1 roll
add
setdash
} bind def
/Mfstroke {
stroke
currentdash
pop 0
setdash
} bind def
/Mrotsboxa {
gsave
dup
/Mrot
exch def
Mrotcheck
Mtmatrix
dup
setmatrix
7 1 roll
4 index
4 index
translate
rotate
3 index
-1 mul
3 index
-1 mul
translate
/Mtmatrix
matrix
currentmatrix
def
grestore
Msboxa
3  -1 roll
/Mtmatrix
exch def
/Mrot
0 def
} bind def
/Msboxa {
newpath
5 -1 roll
Mvboxa
pop
Mboxout
6 -1 roll
5 -1 roll
4 -1 roll
Msboxa1
5 -3 roll
Msboxa1
Mboxrot
[
7 -2 roll
2 copy
[
3 1 roll
10 -1 roll
9 -1 roll
]
6 1 roll
5 -2 roll
]
} bind def
/Msboxa1 {
sub
2 div
dup
2 index
1 add
mul
3 -1 roll
-1 add
3 -1 roll
mul
} bind def
/Mvboxa	{
Mfixwid
{
Mvboxa1
}
{
dup
Mwidthcal
0 exch
{
add
}
forall
exch
Mvboxa1
4 index
7 -1 roll
add
4 -1 roll
pop
3 1 roll
}
ifelse
} bind def
/Mvboxa1 {
gsave
newpath
[ true
3 -1 roll
{
Mbbox
5 -1 roll
{
0
5 1 roll
}
{
7 -1 roll
exch sub
(m) stringwidth pop
.3 mul
sub
7 1 roll
6 -1 roll
4 -1 roll
Mmin
3 -1 roll
5 index
add
5 -1 roll
4 -1 roll
Mmax
4 -1 roll
}
ifelse
false
}
forall
{ stop } if
counttomark
1 add
4 roll
]
grestore
} bind def
/Mbbox {
1 dict begin
0 0 moveto
/temp (T) def
{	gsave
currentpoint newpath moveto
temp 0 3 -1 roll put temp
false charpath flattenpath currentpoint
pathbbox
grestore moveto lineto moveto} forall
pathbbox
newpath
end
} bind def
/Mmin {
2 copy
gt
{ exch } if
pop
} bind def
/Mmax {
2 copy
lt
{ exch } if
pop
} bind def
/Mrotshowa {
dup
/Mrot
exch def
Mrotcheck
Mtmatrix
dup
setmatrix
7 1 roll
4 index
4 index
translate
rotate
3 index
-1 mul
3 index
-1 mul
translate
/Mtmatrix
matrix
currentmatrix
def
Mgmatrix setmatrix
Mshowa
/Mtmatrix
exch def
/Mrot 0 def
} bind def
/Mshowa {
4 -2 roll
moveto
2 index
Mtmatrix setmatrix
Mvboxa
7 1 roll
Mboxout
6 -1 roll
5 -1 roll
4 -1 roll
Mshowa1
4 1 roll
Mshowa1
rmoveto
currentpoint
Mfixwid
{
Mshowax
}
{
Mshoway
}
ifelse
pop pop pop pop
Mgmatrix setmatrix
} bind def
/Mshowax {	
0 1
4 index length
-1 add
{
2 index
4 index
2 index
get
3 index
add
moveto
4 index
exch get
Mfixdash
{
Mfixdashp		
}
if
show
} for
} bind def
/Mfixdashp {
dup
length
1
gt
1 index
true exch
{
45
eq
and
} forall
and
{
gsave
(--)		
stringwidth pop
(-)
stringwidth pop
sub
2 div
0 rmoveto
dup
length
1 sub
{
(-)
show
}
repeat
grestore
}
if
} bind def	
/Mshoway {
3 index
Mwidthcal
5 1 roll
0 1
4 index length
-1 add
{
2 index
4 index
2 index
get
3 index
add
moveto
4 index
exch get
[
6 index
aload
length
2 add
-1 roll
{
pop
Strform
stringwidth
pop
neg 
exch
add
0 rmoveto
}
exch
kshow
cleartomark
} for
pop	
} bind def
/Mwidthcal {
[
exch
{
Mwidthcal1
}
forall
]
[
exch
dup
Maxlen
-1 add
0 1
3 -1 roll
{
[
exch
2 index
{		 
1 index	
Mget
exch
}
forall		
pop
Maxget
exch
}
for
pop
]
Mreva
} bind def
/Mreva	{
[
exch
aload
length
-1 1
{1 roll}
for
]
} bind def
/Mget	{
1 index
length
-1 add
1 index
ge
{
get
}
{
pop pop
0
}
ifelse
} bind def
/Maxlen	{
[
exch
{
length
}
forall
Maxget
} bind def
/Maxget	{
counttomark
-1 add
1 1
3 -1 roll
{
pop
Mmax
}
for
exch
pop
} bind def
/Mwidthcal1 {
[
exch
{
Strform
stringwidth
pop
}
forall
]
} bind def
/Strform {
/tem (x) def
tem 0
3 -1 roll
put
tem
} bind def
/Mshowa1 {
2 copy
add
4 1 roll
sub
mul
sub
-2 div
} bind def
/MathScale {
Mwidth
Mheight
Mlp
translate
scale
/yscale exch def
/ybias exch def
/xscale exch def
/xbias exch def
/Momatrix
xscale yscale matrix scale
xbias ybias matrix translate
matrix concatmatrix def
/Mgmatrix
matrix currentmatrix
def
} bind def
/Mlp {
3 copy
Mlpfirst
{
Mnodistort
{
Mmin
dup
} if
4 index
2 index
2 index
Mlprun
11 index
11 -1 roll
10 -4 roll
Mlp1
8 index
9 -5 roll
Mlp1
4 -1 roll
and
{ exit } if
3 -1 roll
pop pop
} loop
exch
3 1 roll
7 -3 roll
pop pop pop
} bind def
/Mlpfirst {
3 -1 roll
dup length
2 copy
-2 add
get
aload
pop pop pop
4 -2 roll
-1 add
get
aload
pop pop pop
6 -1 roll
3 -1 roll
5 -1 roll
sub
div
4 1 roll
exch sub
div
} bind def
/Mlprun {
2 copy
4 index
0 get
dup
4 1 roll
Mlprun1
3 copy
8 -2 roll
9 -1 roll
{
3 copy
Mlprun1
3 copy
11 -3 roll
/gt Mlpminmax
8 3 roll
11 -3 roll
/lt Mlpminmax
8 3 roll
} forall
pop pop pop pop
3 1 roll
pop pop
aload pop
5 -1 roll
aload pop
exch
6 -1 roll
Mlprun2
8 2 roll
4 -1 roll
Mlprun2
6 2 roll
3 -1 roll
Mlprun2
4 2 roll
exch
Mlprun2
6 2 roll
} bind def
/Mlprun1 {
aload pop
exch
6 -1 roll
5 -1 roll
mul add
4 -2 roll
mul
3 -1 roll
add
} bind def
/Mlprun2 {
2 copy
add 2 div
3 1 roll
exch sub
} bind def
/Mlpminmax {
cvx
2 index
6 index
2 index
exec
{
7 -3 roll
4 -1 roll
} if
1 index
5 index
3 -1 roll
exec
{
4 1 roll
pop
5 -1 roll
aload
pop pop
4 -1 roll
aload pop
[
8 -2 roll
pop
5 -2 roll
pop
6 -2 roll
pop
5 -1 roll
]
4 1 roll
pop
}
{
pop pop pop
} ifelse
} bind def
/Mlp1 {
5 index
3 index	sub
5 index
2 index mul
1 index
le
1 index
0 le
or
dup
not
{
1 index
3 index	div
.99999 mul
8 -1 roll
pop
7 1 roll
}
if
8 -1 roll
2 div
7 -2 roll
pop sub
5 index
6 -3 roll
pop pop
mul sub
exch
} bind def
/intop 0 def
/inrht 0 def
/inflag 0 def
/outflag 0 def
/xadrht 0 def
/xadlft 0 def
/yadtop 0 def
/yadbot 0 def
/Minner {
outflag
1
eq
{
/outflag 0 def
/intop 0 def
/inrht 0 def
} if		
5 index
gsave
Mtmatrix setmatrix
Mvboxa pop
grestore
3 -1 roll
pop
dup
intop
gt
{
/intop
exch def
}
{ pop }
ifelse
dup
inrht
gt
{
/inrht
exch def
}
{ pop }
ifelse
pop
/inflag
1 def
} bind def
/Mouter {
/xadrht 0 def
/xadlft 0 def
/yadtop 0 def
/yadbot 0 def
inflag
1 eq
{
dup
0 lt
{
dup
intop
mul
neg
/yadtop
exch def
} if
dup
0 gt
{
dup
intop
mul
/yadbot
exch def
}
if
pop
dup
0 lt
{
dup
inrht
mul
neg
/xadrht
exch def
} if
dup
0 gt
{
dup
inrht
mul
/xadlft
exch def
} if
pop
/outflag 1 def
}
{ pop pop}
ifelse
/inflag 0 def
/inrht 0 def
/intop 0 def
} bind def	
/Mboxout {
outflag
1
eq
{
4 -1
roll
xadlft
leadjust
add
sub
4 1 roll
3 -1
roll
yadbot
leadjust
add
sub
3 1
roll
exch
xadrht
leadjust
add
add
exch
yadtop
leadjust
add
add
/outflag 0 def
/xadlft 0 def
/yadbot 0 def
/xadrht 0 def
/yadtop 0 def
} if
} bind def
/leadjust {
(m) stringwidth pop
.5 mul
} bind def
/Mrotcheck {
dup
90
eq
{
yadbot
/yadbot
xadrht
def	
/xadrht
yadtop
def
/yadtop
xadlft
def
/xadlft
exch
def
}
if
dup
cos
1 index
sin
Checkaux
dup
cos
1 index
sin neg
exch
Checkaux
3 1 roll
pop pop
} bind def
/Checkaux {
4 index
exch
4 index
mul
3 1 roll
mul add
4 1 roll
} bind def
/Mboxrot {
Mrot
90 eq
{
brotaux
4 2
roll
}
if
Mrot
180 eq
{
4 2
roll
brotaux
4 2
roll
brotaux
}	
if
Mrot
270 eq
{
4 2
roll
brotaux
}
if
} bind def
/brotaux {
neg
exch
neg
} bind def
/Mabsproc {
0
matrix defaultmatrix
dtransform idtransform
dup mul exch
dup mul
add sqrt
} bind def
/Mabswid {
Mabsproc
setlinewidth	
} bind def
/Mabsdash {
exch
[
exch
{
Mabsproc
}
forall
]
exch
setdash
} bind def
/MBeginOrig { Momatrix concat} bind def
/MEndOrig { Mgmatrix setmatrix} bind def
/sampledsound where
{ pop}
{ /sampledsound {
exch
pop
exch
5 1 roll
mul
4 idiv
mul
2 idiv
exch pop
exch
/Mtempproc exch def
{ Mtempproc pop}
repeat
} bind def
} ifelse
/g { setgray} bind def
/k { setcmykcolor} bind def
/m { moveto} bind def
/p { gsave} bind def
/r { setrgbcolor} bind def
/w { setlinewidth} bind def
/C { curveto} bind def
/F { fill} bind def
/L { lineto} bind def
/P { grestore} bind def
/s { stroke} bind def
/setcmykcolor where
{ pop}
{ /setcmykcolor {
4 1
roll
[
4 1 
roll
]
{
1 index
sub
1
sub neg
dup
0
lt
{
pop
0
}
if
dup
1
gt
{
pop
1
}
if
exch
} forall
pop
setrgbcolor
} bind def
} ifelse
/Mcharproc
{
currentfile
(x)
readhexstring
pop
0 get
exch
div
} bind def
/Mshadeproc
{
dup
3 1
roll
{
dup
Mcharproc
3 1
roll
} repeat
1 eq
{
setgray
}
{
3 eq
{
setrgbcolor
}
{
setcmykcolor
} ifelse
} ifelse
} bind def
/Mrectproc
{
3 index
2 index
moveto
2 index
3 -1
roll
lineto
dup
3 1
roll
lineto
lineto
fill
} bind def
/Mcolorimage
{
7 1
roll
pop
pop
matrix
invertmatrix
concat
2 exch exp
1 sub
3 1 roll
1 1
2 index
{
1 1
4 index
{
dup
1 sub
exch
2 index
dup
1 sub
exch
7 index
9 index
Mshadeproc
Mrectproc
} for
pop
} for
pop pop pop pop
} bind def
/Mimage
{
pop
matrix
invertmatrix
concat
2 exch exp
1 sub
3 1 roll
1 1
2 index
{
1 1
4 index
{
dup
1 sub
exch
2 index
dup
1 sub
exch
7 index
Mcharproc
setgray
Mrectproc
} for
pop
} for
pop pop pop
} bind def

MathPictureStart
/Courier findfont 10  scalefont  setfont
0.735744 0.00240447 0.309017 3.98729e-06 [
[(      2)(-30 Pi)] .02381 .30902 0 2 Msboxa
[(      2)(-20 Pi)] .26112 .30902 0 2 Msboxa
[(      2)(-10 Pi)] .49843 .30902 0 2 Msboxa
[(50)] .85597 .30902 0 2 Msboxa
[(100)] .97619 .30902 0 2 Msboxa
[(-75000)] .72324 .00997 1 0 Msboxa
[(-50000)] .72324 .10965 1 0 Msboxa
[(-25000)] .72324 .20933 1 0 Msboxa
[(25000)] .72324 .4087 1 0 Msboxa
[(50000)] .72324 .50838 1 0 Msboxa
[(75000)] .72324 .60806 1 0 Msboxa
[ -0.001 -0.001 0 0 ]
[ 1.001 .61903 0 0 ]
] MathScale
1 setlinecap
1 setlinejoin
newpath
[ ] 0 setdash
0 g
p
p
.002 w
.02381 .30902 m
.02381 .31527 L
s
P
[(      2)(-30 Pi)] .02381 .30902 0 2 Mshowa
p
.002 w
.26112 .30902 m
.26112 .31527 L
s
P
[(      2)(-20 Pi)] .26112 .30902 0 2 Mshowa
p
.002 w
.49843 .30902 m
.49843 .31527 L
s
P
[(      2)(-10 Pi)] .49843 .30902 0 2 Mshowa
p
.002 w
.85597 .30902 m
.85597 .31527 L
s
P
[(50)] .85597 .30902 0 2 Mshowa
p
.002 w
.97619 .30902 m
.97619 .31527 L
s
P
[(100)] .97619 .30902 0 2 Mshowa
p
.001 w
.04754 .30902 m
.04754 .31277 L
s
P
p
.001 w
.07127 .30902 m
.07127 .31277 L
s
P
p
.001 w
.095 .30902 m
.095 .31277 L
s
P
p
.001 w
.11873 .30902 m
.11873 .31277 L
s
P
p
.001 w
.14247 .30902 m
.14247 .31277 L
s
P
p
.001 w
.1662 .30902 m
.1662 .31277 L
s
P
p
.001 w
.18993 .30902 m
.18993 .31277 L
s
P
p
.001 w
.21366 .30902 m
.21366 .31277 L
s
P
p
.001 w
.23739 .30902 m
.23739 .31277 L
s
P
p
.001 w
.26112 .30902 m
.26112 .31277 L
s
P
p
.001 w
.28485 .30902 m
.28485 .31277 L
s
P
p
.001 w
.30858 .30902 m
.30858 .31277 L
s
P
p
.001 w
.33231 .30902 m
.33231 .31277 L
s
P
p
.001 w
.35605 .30902 m
.35605 .31277 L
s
P
p
.001 w
.37978 .30902 m
.37978 .31277 L
s
P
p
.001 w
.40351 .30902 m
.40351 .31277 L
s
P
p
.001 w
.42724 .30902 m
.42724 .31277 L
s
P
p
.001 w
.45097 .30902 m
.45097 .31277 L
s
P
p
.001 w
.4747 .30902 m
.4747 .31277 L
s
P
p
.001 w
.49843 .30902 m
.49843 .31277 L
s
P
p
.001 w
.52216 .30902 m
.52216 .31277 L
s
P
p
.001 w
.54589 .30902 m
.54589 .31277 L
s
P
p
.001 w
.56963 .30902 m
.56963 .31277 L
s
P
p
.001 w
.59336 .30902 m
.59336 .31277 L
s
P
p
.001 w
.61709 .30902 m
.61709 .31277 L
s
P
p
.001 w
.64082 .30902 m
.64082 .31277 L
s
P
p
.001 w
.66455 .30902 m
.66455 .31277 L
s
P
p
.001 w
.68828 .30902 m
.68828 .31277 L
s
P
p
.001 w
.71201 .30902 m
.71201 .31277 L
s
P
p
.001 w
.75979 .30902 m
.75979 .31277 L
s
P
p
.001 w
.78383 .30902 m
.78383 .31277 L
s
P
p
.001 w
.80788 .30902 m
.80788 .31277 L
s
P
p
.001 w
.83192 .30902 m
.83192 .31277 L
s
P
p
.001 w
.85597 .30902 m
.85597 .31277 L
s
P
p
.001 w
.88001 .30902 m
.88001 .31277 L
s
P
p
.001 w
.90406 .30902 m
.90406 .31277 L
s
P
p
.001 w
.9281 .30902 m
.9281 .31277 L
s
P
p
.001 w
.95215 .30902 m
.95215 .31277 L
s
P
p
.002 w
0 .30902 m
1 .30902 L
s
P
p
.002 w
.73574 .00997 m
.74199 .00997 L
s
P
[(-75000)] .72324 .00997 1 0 Mshowa
p
.002 w
.73574 .10965 m
.74199 .10965 L
s
P
[(-50000)] .72324 .10965 1 0 Mshowa
p
.002 w
.73574 .20933 m
.74199 .20933 L
s
P
[(-25000)] .72324 .20933 1 0 Mshowa
p
.002 w
.73574 .4087 m
.74199 .4087 L
s
P
[(25000)] .72324 .4087 1 0 Mshowa
p
.002 w
.73574 .50838 m
.74199 .50838 L
s
P
[(50000)] .72324 .50838 1 0 Mshowa
p
.002 w
.73574 .60806 m
.74199 .60806 L
s
P
[(75000)] .72324 .60806 1 0 Mshowa
p
.001 w
.73574 .02991 m
.73949 .02991 L
s
P
p
.001 w
.73574 .04984 m
.73949 .04984 L
s
P
p
.001 w
.73574 .06978 m
.73949 .06978 L
s
P
p
.001 w
.73574 .08972 m
.73949 .08972 L
s
P
p
.001 w
.73574 .10965 m
.73949 .10965 L
s
P
p
.001 w
.73574 .12959 m
.73949 .12959 L
s
P
p
.001 w
.73574 .14953 m
.73949 .14953 L
s
P
p
.001 w
.73574 .16946 m
.73949 .16946 L
s
P
p
.001 w
.73574 .1894 m
.73949 .1894 L
s
P
p
.001 w
.73574 .20933 m
.73949 .20933 L
s
P
p
.001 w
.73574 .22927 m
.73949 .22927 L
s
P
p
.001 w
.73574 .24921 m
.73949 .24921 L
s
P
p
.001 w
.73574 .26914 m
.73949 .26914 L
s
P
p
.001 w
.73574 .28908 m
.73949 .28908 L
s
P
p
.001 w
.73574 .32895 m
.73949 .32895 L
s
P
p
.001 w
.73574 .34889 m
.73949 .34889 L
s
P
p
.001 w
.73574 .36883 m
.73949 .36883 L
s
P
p
.001 w
.73574 .38876 m
.73949 .38876 L
s
P
p
.001 w
.73574 .4087 m
.73949 .4087 L
s
P
p
.001 w
.73574 .42864 m
.73949 .42864 L
s
P
p
.001 w
.73574 .44857 m
.73949 .44857 L
s
P
p
.001 w
.73574 .46851 m
.73949 .46851 L
s
P
p
.001 w
.73574 .48844 m
.73949 .48844 L
s
P
p
.001 w
.73574 .50838 m
.73949 .50838 L
s
P
p
.001 w
.73574 .52832 m
.73949 .52832 L
s
P
p
.001 w
.73574 .54825 m
.73949 .54825 L
s
P
p
.001 w
.73574 .56819 m
.73949 .56819 L
s
P
p
.001 w
.73574 .58813 m
.73949 .58813 L
s
P
p
.002 w
.73574 0 m
.73574 .61803 L
s
P
P
0 0 m
1 0 L
1 .61803 L
0 .61803 L
closepath
clip
newpath
p
p
.5 g
p
.002 w
.02381 .12568 m
.03343 .13063 L
.03824 .13307 L
.04305 .13549 L
.04545 .13669 L
.04666 .13729 L
.04696 .13744 L
.04726 .13759 L
.04786 .13804 L
.04906 .13918 L
.05026 .14032 L
.05267 .14259 L
.06229 .15162 L
.07191 .16049 L
.08153 .16912 L
.09115 .17745 L
.10077 .18537 L
.11039 .19282 L
.12001 .19972 L
.12963 .20604 L
.13925 .21175 L
.14887 .21689 L
.15849 .22151 L
.16811 .22566 L
.17773 .22942 L
.18735 .23285 L
.19697 .23601 L
.20659 .23895 L
.21621 .24171 L
.22583 .24432 L
.23545 .24682 L
.24507 .24923 L
.24988 .2504 L
.25469 .25156 L
.25709 .25213 L
.2583 .25242 L
.2595 .2527 L
.2601 .25285 L
.2604 .25292 L
.2607 .25299 L
.261 .25306 L
.2613 .25321 L
.2619 .2536 L
.26431 .25517 L
.27393 .26142 L
.28355 .2675 L
.29317 .27325 L
.30279 .27843 L
.31241 .28283 L
.32203 .28634 L
Mistroke
.33165 .28903 L
.34127 .29106 L
.35089 .29264 L
.36051 .29391 L
.37013 .29496 L
.37975 .29587 L
.38937 .29667 L
.39899 .29741 L
.40861 .29809 L
.41823 .29873 L
.42304 .29904 L
.42424 .29912 L
.42544 .29919 L
.42605 .29923 L
.42665 .29927 L
.42695 .29929 L
.42725 .29931 L
.42785 .29957 L
.43747 .30379 L
.44709 .30823 L
.45671 .31286 L
.46633 .31764 L
.47595 .32256 L
.48557 .32759 L
.49519 .3327 L
.50481 .33788 L
.51443 .34312 L
.52405 .34841 L
.53367 .35374 L
.54329 .35909 L
.55291 .36448 L
.56253 .36988 L
.57215 .3753 L
.58177 .38073 L
.59139 .38617 L
.60101 .39162 L
.61063 .39708 L
.62025 .40254 L
.62987 .40801 L
.63949 .41348 L
.64911 .41896 L
.65873 .42444 L
.66835 .42991 L
.67797 .43539 L
.68759 .44087 L
.69721 .44635 L
.70683 .45183 L
.71645 .4573 L
.72607 .46278 L
.73569 .46825 L
Mistroke
.74531 .47373 L
.75493 .4792 L
.76455 .48467 L
.77417 .49014 L
.78379 .49561 L
.79341 .50107 L
.80303 .50654 L
.81265 .512 L
.82227 .51746 L
.83189 .52292 L
.84151 .52837 L
.85113 .53383 L
.86075 .53928 L
.87037 .54473 L
.87999 .55018 L
.88961 .55562 L
.89923 .56107 L
.90885 .56651 L
.91847 .57195 L
.92809 .57739 L
.93771 .58282 L
.94733 .58826 L
.95695 .59369 L
.96657 .59912 L
.97619 .60455 L
Mfstroke
P
P
p
p
.002 w
.02381 .11511 m
.03343 .12429 L
.04305 .13346 L
.04545 .13575 L
.04666 .13689 L
.04726 .13747 L
.04756 .13774 L
.04786 .13789 L
.04846 .13819 L
.04906 .13849 L
.05026 .13909 L
.05267 .14028 L
.06229 .14501 L
.07191 .14969 L
.08153 .15435 L
.09115 .15899 L
.10077 .16362 L
.11039 .16827 L
.12001 .17295 L
.12963 .17767 L
.13925 .18244 L
.14887 .18729 L
.15849 .19222 L
.16811 .19726 L
.17773 .20241 L
.18735 .2077 L
.19697 .21314 L
.20659 .21873 L
.21621 .22448 L
.22583 .23039 L
.23545 .23644 L
.24507 .24262 L
.25469 .24888 L
.25709 .25045 L
.2595 .25203 L
.2601 .25242 L
.2607 .25281 L
.261 .25301 L
.2613 .25313 L
.2619 .25327 L
.26311 .25355 L
.26431 .25384 L
.27393 .25607 L
.28355 .25828 L
.29317 .26047 L
.30279 .26266 L
.31241 .26487 L
.32203 .2671 L
.33165 .26937 L
.34127 .2717 L
Mistroke
.35089 .27412 L
.36051 .27664 L
.37013 .27929 L
.37975 .28211 L
.38937 .28513 L
.39899 .28837 L
.40861 .29185 L
.41823 .29558 L
.42304 .29754 L
.42544 .29855 L
.42605 .2988 L
.42665 .29906 L
.42695 .29918 L
.42725 .29931 L
.42785 .29935 L
.43025 .2995 L
.43266 .29965 L
.43747 .29994 L
.44709 .30053 L
.45671 .30111 L
.46633 .30169 L
.47595 .30229 L
.48557 .30291 L
.49519 .30357 L
.50481 .30428 L
.51443 .30509 L
.52405 .30602 L
.53367 .30716 L
.54329 .30858 L
.55291 .31034 L
.56253 .31247 L
.57215 .31492 L
.58177 .3176 L
.59139 .32045 L
.60101 .3234 L
.61063 .32643 L
.62025 .3295 L
.62987 .3326 L
.63949 .33573 L
.64911 .33887 L
.65873 .34202 L
.66835 .34517 L
.67797 .34833 L
.68759 .35149 L
.69721 .35465 L
.70683 .35781 L
.71645 .36097 L
.72607 .36413 L
.73569 .36728 L
.74531 .37044 L
Mistroke
.75493 .37359 L
.76455 .37674 L
.77417 .37989 L
.78379 .38303 L
.79341 .38618 L
.80303 .38931 L
.81265 .39245 L
.82227 .39559 L
.83189 .39872 L
.84151 .40185 L
.85113 .40497 L
.86075 .4081 L
.87037 .41122 L
.87999 .41434 L
.88961 .41745 L
.89923 .42057 L
.90885 .42368 L
.91847 .42679 L
.92809 .4299 L
.93771 .433 L
.94733 .4361 L
.95695 .4392 L
.96657 .4423 L
.97619 .4454 L
Mfstroke
P
P
p
.5 g
p
.002 w
.02381 .0097 m
.03343 .0181 L
.04305 .02631 L
.05267 .03434 L
.06229 .0422 L
.07191 .04989 L
.08153 .05742 L
.09115 .06482 L
.10077 .07207 L
.11039 .07921 L
.12001 .08638 L
.12963 .0945 L
.13925 .10248 L
.14887 .11032 L
.15849 .118 L
.16811 .12552 L
.17773 .13285 L
.18735 .13999 L
.19697 .14691 L
.20659 .15361 L
.21621 .16008 L
.22583 .16631 L
.23545 .1723 L
.24507 .17806 L
.25469 .18359 L
.26431 .18891 L
.27393 .19402 L
.28355 .19895 L
.29317 .20371 L
.30279 .20831 L
.31241 .21277 L
.32203 .21711 L
.33165 .22134 L
.34127 .22656 L
.35089 .23174 L
.36051 .23677 L
.37013 .24163 L
.37975 .24629 L
.38937 .25073 L
.39899 .25491 L
.40861 .25882 L
.41823 .26245 L
.42785 .2658 L
.43747 .2689 L
.44709 .27175 L
.45671 .2744 L
.46633 .27687 L
.47595 .27919 L
.48557 .28139 L
.49519 .28348 L
Mistroke
.50481 .28623 L
.51443 .28925 L
.52405 .29211 L
.53367 .29475 L
.54329 .29709 L
.55291 .29907 L
.56253 .30067 L
.57215 .30193 L
.58177 .30294 L
.59139 .30378 L
.60101 .3045 L
.60582 .30482 L
.61063 .30513 L
.61304 .30528 L
.61424 .30535 L
.61544 .30542 L
.61604 .30546 L
.61634 .30548 L
.61664 .3055 L
.61694 .30551 L
.61724 .30555 L
.61754 .3056 L
.61785 .30565 L
.62025 .30604 L
.62987 .3076 L
.63949 .30915 L
.64911 .31069 L
.65873 .31222 L
.66835 .31374 L
.67797 .31525 L
.68759 .31676 L
.69721 .31825 L
.70683 .31974 L
.71645 .32123 L
.72607 .3227 L
.73569 .32417 L
.74531 .32564 L
.75493 .3271 L
.76455 .32856 L
.77417 .33001 L
.78379 .33146 L
.79341 .33291 L
.80303 .33435 L
.81265 .33579 L
.82227 .33722 L
.83189 .33866 L
.84151 .34009 L
.85113 .34151 L
.86075 .34294 L
.87037 .34436 L
Mistroke
.87999 .34578 L
.88961 .34719 L
.89923 .34861 L
.90885 .35002 L
.91847 .35143 L
.92809 .35284 L
.93771 .35424 L
.94733 .35565 L
.95695 .35705 L
.96657 .35845 L
.97619 .35985 L
Mfstroke
P
P
p
p
.002 w
.02381 0 m
.03343 .00899 L
.04305 .01791 L
.05267 .02676 L
.06229 .03554 L
.07191 .04423 L
.08153 .05285 L
.09115 .06138 L
.10077 .06982 L
.11039 .07816 L
.12001 .08622 L
.12963 .09313 L
.13925 .09994 L
.14887 .10666 L
.15849 .11329 L
.16811 .11984 L
.17773 .12631 L
.18735 .13271 L
.19697 .13904 L
.20659 .14531 L
.21621 .15152 L
.22583 .15766 L
.23545 .16375 L
.24507 .16979 L
.25469 .17576 L
.26431 .18168 L
.27393 .18754 L
.28355 .19334 L
.29317 .19908 L
.30279 .20475 L
.31241 .21034 L
.32203 .21584 L
.33165 .22126 L
.34127 .22547 L
.35089 .2295 L
.36051 .23346 L
.37013 .23734 L
.37975 .24116 L
.38937 .24492 L
.39899 .24862 L
.40861 .25227 L
.41823 .25588 L
.42785 .25944 L
.43747 .26296 L
.44709 .26643 L
.45671 .26987 L
.46633 .27325 L
.47595 .27659 L
.48557 .27988 L
.49519 .28309 L
Mistroke
.50481 .28548 L
.51443 .28741 L
.52405 .28928 L
.53367 .29109 L
.54329 .29287 L
.55291 .2946 L
.56253 .29631 L
.57215 .29798 L
.58177 .29963 L
.59139 .30126 L
.60101 .30287 L
.61063 .30446 L
.61304 .30486 L
.61544 .30525 L
.61604 .30535 L
.61664 .30545 L
.61694 .3055 L
.61724 .30553 L
.61754 .30555 L
.61785 .30557 L
.62025 .30571 L
.62987 .30624 L
.63949 .30674 L
.64911 .30722 L
.65873 .30768 L
.66835 .30813 L
.67797 .30856 L
.68759 .30899 L
.69721 .3094 L
.70683 .30981 L
.71645 .31022 L
.72607 .31062 L
.73569 .31101 L
.74531 .3114 L
.75493 .31179 L
.76455 .31217 L
.77417 .31255 L
.78379 .31293 L
.79341 .31331 L
.80303 .31368 L
.81265 .31406 L
.82227 .31443 L
.83189 .3148 L
.84151 .31517 L
.85113 .31553 L
.86075 .3159 L
.87037 .31626 L
.87999 .31663 L
.88961 .31699 L
.89923 .31735 L
Mistroke
.90885 .31771 L
.91847 .31807 L
.92809 .31843 L
.93771 .31879 L
.94733 .31914 L
.95695 .3195 L
.96657 .31985 L
.97619 .32021 L
Mfstroke
P
P
P
MathPictureEnd
end
EndOfTheIncludedPostscriptMagicCookie

\closepsdump

\usepackage[dvips]{epsfig}
\usepackage{amstex,ifthen,amssymb,alltt,calc,epic}

\allowdisplaybreaks

\makeatletter
\renewcommand{\@seccntformat}[1]{\csname the#1\endcsname.\hspace{1em}}
\renewcommand\section{\@startsection {section}{1}{\z@}
                                   {-3.5ex \@plus -1ex \@minus -.2ex}
                                   {2.3ex \@plus.2ex}
                                   {\normalfont\normalsize\scshape\centering}}
\renewcommand\subsection{\@startsection{subsection}{2}{\z@}
                                     {-3.25ex\@plus -1ex \@minus -.2ex}
                                     {1.5ex \@plus .2ex}
                                     {\normalfont\normalsize\scshape\centering}}
\renewcommand{\l@section}[2]{
  \ifnum \c@tocdepth >\z@
    \addpenalty{\@secpenalty}
    \addvspace{\p@}
    \setlength\@tempdima{1.5em}
    \begingroup
      \parindent \z@ \rightskip \@pnumwidth
      \parfillskip -\@pnumwidth
      \leavevmode \scshape
      \advance\leftskip\@tempdima
      \hskip -\leftskip
      #1\nobreak\hfil \nobreak\hbox to\@pnumwidth{\hss #2}\par
    \endgroup
  \fi}
\def\@maketitle{
  \newpage
  \null
  \vskip 2em
  \begin{center}
    {\sffamily \large \@title \par}
    \vskip 1.5em
      {\normalsize
      \lineskip .5em
        {\scshape \normalsize \@author}}
  \end{center}
  \par
  \vskip 1.5em}
\makeatother

\newcommand\mathematica{\textsl{Mathe\-matica}}
\newcommand\real{\mathbb R}

\renewcommand\natural{\mathbb N}
\newcommand\re{\operatorname{Re}}
\newcommand\dom{\operatorname{Dom}}
\newcommand\dir{|_{\operatorname{DIR}}} 
\newcommand\ip[2]{\langle #1,#2\rangle}
\newcommand\bip[2]{\big\langle #1,#2\big\rangle}

\newcommand\dup{\mathrm{d}}
\newcommand\finbox{~\hfill$\Box$}
\newcommand\earlybox{\tag*{$\Box$}}%
\newcommand\lowearlybox{\tag*{$\frac[0pt]{\displaystyle{\phantom{\big|}}}{%
  \displaystyle{\Box}}$}}%

\newcommand\fgap{\,}
\newcommand\cbytwo{(c/2)}
\newcommand\oper{H(h,\alpha)}
\newcommand\bihar{\Delta^2\dir}
\newcommand\fmtwo{\|f_1^-\|\makebox[0cm][l]{$_2$}}
\newcommand\ftwo{\|f_1^{\phantom-}\|\makebox[0cm][l]{$_2$}}
\newcommand\fminfty{\|f_1^-\|\makebox[0cm][l]{$_\infty$}}
\newcommand\finfty{\|f_1^{\phantom-}\|\makebox[0cm][l]{$_\infty$}}
\newcommand\slref[1]{\textsl{\ref{#1}}}%
\renewcommand\eqref[1]{(\ref{#1})}
\newcommand\sleqref[1]{\textsl{\eqref{#1}}}%

\newcommand\afigwidth0%
\newcommand\afigheight0%
\newcommand\bfigwidth0%
\newcommand\bfigheight0%
\newlength\aeqspace%
\newlength\beqspace%
\newlength\border%
\newlength\tempeqspace%
\newcounter{tempactr}%
\newcounter{tempbctr}%
\newcounter{tempcctr}%
1%

\newcommand{\fig}[6]{%
  \setcounter{tempactr}{#6}%
  \setcounter{tempbctr}{200+\border/\unitlength}%
  \setcounter{tempcctr}{#6+\border/\unitlength}%
  \setlength\tempeqspace{#5}%
  \vskip\tempeqspace%
  \vskip\border%
  \centering%
  \begin{picture}(200,\value{tempactr})%
    \put(0,0){\epsfig{figure=#1.eps,%
                      width=200\unitlength}}%
    #3%
  \end{picture}%
  \vskip\border%
  \vskip\tempeqspace%
  \vskip\abovecaptionskip%
  \refstepcounter{figure}#2%
  \makebox[0pt][c]{\parbox[t]{1.5\textwidth}{\centering%
    \ifthenelse{\equal{#4}{}}{Figure \thefigure}{Figure \thefigure: #4}}}%
  \vskip\belowcaptionskip}%

\newcommand{\tab}[6]{%
  \setcounter{tempactr}{#6}%
  \setcounter{tempbctr}{200+\border/\unitlength}%
  \setcounter{tempcctr}{#6+\border/\unitlength}%
  \setlength\tempeqspace{#5}%
  \vskip\tempeqspace%
  \vskip\border%
  \centering%
  {\small%
  \begin{picture}(200,\value{tempactr})%
    \put(0,0){\begin{tabular}[b]{#2}%
      #3%
    \end{tabular}}%
  \end{picture}}%
  \vskip\border%
  \vskip\tempeqspace%
  \vskip\abovecaptionskip%
  \refstepcounter{table}#1%
  \makebox[0pt][c]{\parbox[t]{1.5\textwidth}{\centering%
    \ifthenelse{\equal{#4}{}}{Table \thetable}{Table \thetable: #4}}}%
  \vskip\belowcaptionskip}%

\newcommand{\dbledgram}[3]{{\setlength\border{0pt}
  \renewcommand\baselinestretch1\normalsize%
  \setlength\aeqspace{0pt}%
  \setlength\beqspace{0pt}%
  #1%
  \ifthenelse{\bfigheight > \afigheight}{%
    \setlength\aeqspace{\textwidth/200*(\bfigheight-\afigheight)}}{%
    \setlength\beqspace{\textwidth/200*(\afigheight-\bfigheight)}}%
  \begin{figure}[t]%
    \par%
    \noindent%
    \hfill%
    \setlength\unitlength{\textwidth/100*\afigwidth/200}%
    \begin{minipage}[t]{\textwidth/100*\afigwidth+2\border}%
      #2\aeqspace{200*\afigheight/\afigwidth}%
    \end{minipage}%
    \hfill%
    \setlength\unitlength{\textwidth/100*\bfigwidth/200}%
    \begin{minipage}[t]{\textwidth/100*\bfigwidth+2\border}%
      #3\beqspace{200*\bfigheight/\bfigwidth}%
    \end{minipage}%
    \hfill%
    \null%
  \end{figure}}}%

\newcommand{\sngldiagm}[2]{{%
  \setlength\border{0pt}%
  \renewcommand\baselinestretch1\normalsize%
  #1%
  \begin{figure}[t]%
    \par%
    \noindent%
    \hfill%
    \setlength\unitlength{\textwidth/100*\afigwidth/200}%
    \begin{minipage}[t]{\textwidth/100*\afigwidth+2\border}%
      #2{0pt}{200*\afigheight/\afigwidth}%
    \end{minipage}%
    \hfill%
    \null%
  \end{figure}}}%

\newcounter{theorem}%
\newenvironment{thm}[1]{%
    \refstepcounter{theorem}%
    \textbf{#1 \thetheorem:}{}%
    \begin{itshape}}%
  {\end{itshape}}%

\newenvironment{proof}{\textbf{Proof:}{}}{}%

\newenvironment{note}{%
  \refstepcounter{theorem}%
  \textbf{Note \thetheorem:}{}}%
  {}%

\begin{document}

\title{Asymptotic First Eigenvalue Estimates For\\
  The Biharmonic Operator On A Rectangle}
\author{M. P. Owen\\[1em]
  \slshape{\small Department of Mathematics, King's College London, Strand,
    London, WC2R 2LS}}
\maketitle

\begin{abstract}
  We find an asymptotic expression for the first eigenvalue of the biharmonic
  operator on a long thin rectangle. This is done by finding lower and upper
  bounds which become increasingly accurate with increasing length. The lower
  bound is found by algebraic manipulation of the operator, and the upper bound
  is found by minimising the quadratic form for the operator over a test space
  consisting of separable functions. These bounds can be used to show that the
  negative part of the groundstate is small.
\end{abstract}%
\section{Introduction}%
There is a considerable literature on studies of the biharmonic operator acting 
in $L^2(\Omega)$ for particular regions $\Omega\subseteq\real^2$, especially the
square, disk and punctured disk. In this paper we study the operator acting 
in $L^2(R_h)$, where $R_h$ is the rectangle $[0,h]\times[0,1]$. Difficulties in 
studying the biharmonic operator arise because the eigenvalue problem
\begin{equation}\label{eqn:eprob}
  \Delta^2f=\mu f,\qquad f=\frac{\partial f}{\partial n}=0\textup{ on }
  \partial R_h
\end{equation}
is not exactly soluble. The boundary conditions in this problem are 
called Dirichlet boundary conditions or clamped plate boundary conditions. 
Numerical analysts have succeeded in proving a number of interesting results 
about the groundstate of the biharmonic operator and the corresponding 
eigenvalue for the square and some rectangles, however there are very few 
previous results concerning the $h$-dependence of spectral quantities. One
such result is by Behnke and Mertins~\cite{BehnMert}. See 
note~\ref{not:evalveer}. 

For the unit square the best current enclosure
\[ \mu_1=1294.9339_{40}^{88} \]
for the first eigenvalue is due to C. Wieners~\cite{Wien3} using the 
Tempel-Lehmann-G\"orisch method to obtain the lower bound, and minimisation of 
the quadratic form of the operator on a certain space of test functions for the 
upper bound. The enclosure is guaranteed by interval arithmetic programming. 
The accuracy with which this value has been computed has increased with 
computing power over the last sixty years. In 1937 Weinstein~\cite{Wein} 
introduced a method which theoretically enabled him to calculate a limiting 
sequence of lower bounds for $\mu_1$, although his hand calculations
\[ 1294.956\leq\mu_1\leq1302.360 \]
were slightly inaccurate. Results of intermediate accuracy have been given by
Aronszajn~\cite{Aron} (from a citation in~\cite[page 
74]{Fich}), Bazley, Fox and Stadter~\cite{BazlFoxStad}, De~Vito, Fichera,
Fusciardi and Sch{\"a}rf~\cite{DeviFichFuscScha} and many others.

The biharmonic operator, which we shall denote by $\bihar$, is defined as the
non-negative self-adjoint operator associated with the closed quadratic form
\begin{equation}
  Q(f)=
  \begin{cases}
    {\displaystyle\int_{R_h}|\Delta f|^2},&\textup{if }f\in W^{2,2}_0(R_h);\\
    \infty,&\textup{otherwise}.
  \end{cases}
\end{equation}
See~\cite[theorem 4.4.2]{Davi5} for details. Using the Rayleigh-Ritz
formula~\cite[section 4.5]{Davi5}, the first eigenvalue of the biharmonic
operator is given by the expression
\begin{equation}\label{eqn:asymu1}
  \mu_1(h)=\inf\{Q(f):f\in L^2(R_h),\|f\|_2=1\}.
\end{equation}
We give formulae for lower and upper bounds $\lambda_1(h),\nu_1(h)$ for
$\mu_1(h)$ and use their asymptotic expressions to prove that
\[ \mu_1(h)=c^4+2dc^2\pi^2h^{-2}+O(h^{-3}) \]
as $h\rightarrow\infty$, where $c\approx4.73004$ is the first positive
solution of the transcendental equation
\begin{gather*}
  \refstepcounter{equation} \cosh c\cos c=1, \\
  \intertext{and\hfill(\theequation)}
  d=\frac{2\tanh c\tan c-c\tanh c-c\tan c}{c\tanh c-c\tan c}\approx0.54988.
\end{gather*}%

Elementary algebraic manipulation is used to find the lower bound;
no benefit is derived for large $h$ by using more involved methods such as
Weinstein's truncation of operators. The upper bound is found, as in most other
papers, by minimisation of the quadratic form of the operator over a certain
test space of functions. Observing that eigenfunctions of the biharmonic
operator are close to being separable functions, we choose our test space to 
consist of all such functions. This simplicity of approach allows us to find the
asymptotic formulae. The lower and upper bounds we find are also useful for 
small values of $h$. They are within $0.72\%$ of the actual eigenvalue for all 
$h\in[1,\infty)$. See figure~\ref{fig:l1n1} and table~\ref{tab:lowupp}.

The advantage of our results lies in the fact that they are valid for all 
values of $h$, and also in the limit as $h\rightarrow\infty$. They are also 
simple to compute, without the need for finite element algorithms. The natural 
traded disadvantage is that for a particular rectangle with comparable side 
lengths, where $h$ is close to $1$, numerical analysts are able to use powerful 
computers using finite element methods to compute eigenvalues to a far higher 
degree of accuracy.

A fundamental issue in the study of fourth order operators is the fact that
the corresponding semigroup is not positivity preserving. This is exhibited
by non-positivity of the groundstate of the biharmonic operator for certain
domains $\Omega$, a feature first noticed by Bauer and Reiss
(1972)~\cite{BaueReis} for the square. A rigorous proof of this fact
has been given by Wieners (1995)~\cite{Wien2} who finds a function which
is slightly negative very near the corners of the square and pointwise
extremely close to the actual groundstate. Kozlov, Kondrat'ev, and Maz'ya
(1990)~\cite{KozlKondMazy} had in fact already obtained a more
informative result. They managed to show that if the region has an internal
angle of less than $146.30^\circ$ then the groundstate oscillates infinitely
often in sign as one approaches the corner, although there are no explicit
data concerning where the first oscillation occurs.

Using the bounds $\lambda_1$, $\nu_1$ and also a lower bound $\lambda_3$ on the
third eigenvalue, we prove the bound
\begin{equation}\label{eqn:l2bound}
  \frac{\fmtwo}{\ftwo}_{\phantom2}\leq\frac{(\nu_1-\lambda_1)^{1/2}}{(
  \lambda_3-\nu_1)^{1/2}}=O(h^{-1/2})
\end{equation}
as $h\rightarrow\infty$, on the size of the negative part $f_1^-$ of the
groundstate $f_1$ of $\bihar$. See figure~\ref{fig:posit} and
table~\ref{tab:posit} for a plot and values of this function for small values
of $h$. In particular we see that for the case of the square,
\[ \frac{\fmtwo}{\ftwo}_{\phantom2}\leq0.0484. \]
By `negative part of the groundstate' we mean the negative part when the
eigenfunction is positive in the centre of the square.

It is possible to use the $L^2$ bound and a Sobolev embedding theorem to
imply that
\[ \frac{\fminfty}{\ftwo}_{\phantom\infty}\leq\frac{{(\nu_1-\lambda_1)}
   ^{1/2}\lambda_3^{1/4}}{2(\lambda_3-\nu_1)^{1/2}}=O(h
   ^{-1/2}) \]
as $h\rightarrow\infty$, and so the negative part of the eigenfunction $f_1$,
already small for $h=1$, vanishes asymptotically as $h\rightarrow\infty$ in
the $L^2$ and $L^\infty$ senses. However since
\[ \frac{\|f_1\|\makebox[0cm][l]{$_\infty$}}{\|f_1\|\makebox[0cm][l]{$_2$}}
   _{\phantom\infty}=O(h^{-1/2}), \]
only the $L^2$ result~\eqref{eqn:l2bound} is of particular value. See
note~\ref{note} for details.%
\section{A Fourth Order Operator in one Dimension}%
Our approach to analysis of the biharmonic operator acting in $L^2(R_h)$
involves attempting to separate variables. Eigenfunctions of $\bihar$ are not
separable functions, a fact confirmed by oscillations at the corners, so it
is remarkable that we obtain such good estimates for the first eigenvalue.
The method is successful because the eigenfunctions are close to being
separable. The following sections will rely heavily upon spectral analysis of
the self-adjoint operator
\begin{equation}
  H(h,\alpha)=\frac{\dup^4\phantom x}{\dup x^4}-2\alpha\frac{\dup^2\phantom x}{\dup x^2}
\end{equation}
acting in $L^2([0,h])$, with quadratic form domain $W^{2,2}_0([0,h])$. Since
$\oper$ is bounded below and has compact resolvent we may order the
eigenvalues as an increasing list.

Let $\sigma(h,\alpha,n):(0,\infty)\times\real\times\natural\rightarrow\real$ be
the function which associates to each $\alpha$ the $n$-th eigenvalue
of the operator $\oper$ acting in $L^2([0,h])$, and let
$\rho_n(\alpha)=\sigma(1,\alpha,n)$.

\begin{thm}{Theorem}\label{thm:analysis}
  \begin{enumerate}
    \item
      The first eigenvalue $\rho_1$ of $H(1,\alpha)$ is an increasing and 
      concave function of $\alpha$. The functions $\sigma$ and $\rho$ are 
      related by the equation
      \begin{equation}\label{eqn:rhosig}
        \sigma(h,\alpha,n)=\rho_n(h^2\alpha)h^{-4}.
      \end{equation}
      For $\alpha>0$, $\rho_1$ is analytic\textup{\/;} 
    \item
      For $\alpha>0$ the Green's function of $H(h,\alpha)$ is positive. It
      follows that for such $\alpha$ the first eigenvalue is of multiplicity
      one and the groundstate is positive\textup{\/;}
    \item
      For $\alpha>0$ let $f$ be the $n$-th eigenfunction of $H(h,\alpha)$ and
      have unit $L^2$ norm. Then
      \begin{equation}\label{eqn:norm_f'}
        \|f'\|_2^2=\frac12\rho_n'(h^2\alpha)h^{-2}\mathrm;
      \end{equation}
    \item
      For $\alpha>0$ let $\beta>\alpha$ be the $n$-th solution of the
      transcendental equation
      \begin{equation}\label{eqn:trans}
        \cosh\sqrt{\beta+\alpha}\cos\sqrt{\beta-\alpha}-\frac\alpha{\sqrt{
        \beta^2-\alpha^2}}\sinh\sqrt{\beta+\alpha}\sin\sqrt{\beta-\alpha}=1.
      \end{equation}
      Then $\rho_n(\alpha)=\beta^2-\alpha^2$\textup{\/;}
    \item
      The following are asymptotic expansions of $\rho_n$ at $0$ and 
      $\infty$\textup{\/:}
      \begin{alignat}{2}\label{eqn:forms}
        \rho_n(\alpha)&=2n^2\pi^2\alpha+4\sqrt2n^2\pi^2\alpha^{1/2}+O(1)
        &&\text{as }\alpha\rightarrow\infty\mathrm;\notag\\
        \rho_n'(\alpha)&=2n^2\pi^2+2\sqrt2n^2\pi^2\alpha^{-1/2}+O(\alpha
        ^{-1})\qquad&&\text{as }\alpha\rightarrow\infty\mathrm;\notag \\
        \rho_n(\alpha)&=c_n^4+2d_nc_n^2\alpha+O(\alpha^2)&&\text{as }\alpha
        \rightarrow0\mathrm;\notag\\
        \rho_n'(\alpha)&=2d_nc_n^2+O(\alpha)&&\text{as }\alpha\rightarrow0
        \mathrm;
      \end{alignat}
      where $c_n$ is the $n$-th positive solution of the equation
      $\cosh c\cos c=1$, and
      \[ d_n=\frac{2\tanh c_n\tan c_n-c_n\tanh c_n-c_n\tan c_n}{c_n\tanh c_n-
         c_n\tan c_n}\mathrm; \]
    \item
      The function $\rho_1'$ is monotone decreasing and
      \begin{equation}\label{eqn:dr_ineq}
        2\pi^2\leq\rho_1'(\alpha)\leq2d_1c_1^2
      \end{equation}
      for all $\alpha>0$.
  \end{enumerate}%
\end{thm}%
\sngldiagm{%
  \renewcommand\afigwidth{65}%
  \renewcommand\afigheight{39}}{%
  \fig{shadr1234}{\label{fig:rho}}{%
  \put(196,64){$\rho_1$}%
  \put(196,70){$\rho_2$}%
  \put(196,86){$\rho_3$}%
  \put(196,116){$\rho_4$}%
  \put(0,60){$\alpha$}%
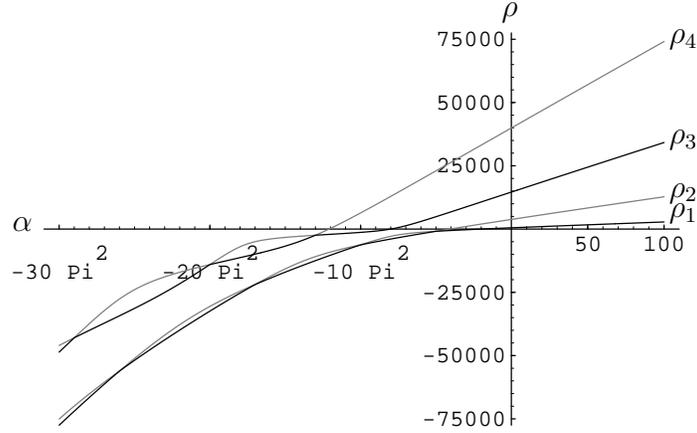
  \put(146,124){$\rho$}}{Plot of the first four eigenvalues of $H(1,\alpha)$}}%
Figure~\ref{fig:rho} shows the first four eigenvalues of $H(1,\alpha)$.
It was plotted by \mathematica{}, using the implicit formula~\eqref{eqn:trans} 
in theorem~\ref{thm:analysis} (iv) for positive $\rho$, with a similar formula
for negative $\rho$.

\begin{note}
  The portion of the graph for negative $\alpha$ is not relevant in this paper
  but has been included to demonstrate that even simple fourth order operators 
  have completely different eigenvalue behaviour to the second order theory. 
  Coincidence of $\rho_{2n-1}$ and $\rho_{2n}$ occurs for
  \begin{equation}
    (\alpha,\rho)=(-(m^2+n^2)\pi^2,-(m^2-n^2)^2\pi^4),
  \end{equation}
  where $m-n\in\natural$.\finbox
\end{note}

\textbf{Proof of theorem~\ref{thm:analysis} (i):}
  The function $\rho_1$ is increasing because the perturbing operator is 
  positive. Using the Rayleigh-Ritz formula,
  \[ \rho_1(\alpha)=\inf\{\ip{H(1,\alpha)g}g_1:\|g\|_2=1,g\in\dom(H(1,\alpha))
     \}. \]
  Suppose that $\alpha=\lambda\alpha_1+(1-\lambda)\alpha_0$ where
  $0<\lambda<1$. Then
  \[ \ip{H(1,\alpha)g}g_1=\lambda\ip{H(1,\alpha_1)g}g_1+(1-\lambda)\ip{H(1,
     \alpha_0)g}g_1. \]
  Let $\epsilon>0$ and choose $g\in\dom(H(1,\alpha))$ such that
  $\|g\|_2=1$ and $\ip{H(1,\alpha)g}g_1<\rho_1(\alpha)+\epsilon$. Then
  \begin{align*}
    \rho_1(\alpha)+\epsilon&>\lambda\ip{H(1,\alpha_1)g}g_1+(1-\lambda)\ip{H(1,
    \alpha_0)g}g_1\\
    &\geq\lambda\rho_1(\alpha_1)+(1-\lambda)\rho_1(\alpha_0).
  \end{align*}
  Since this holds for all positive $\epsilon$, $\rho_1(\alpha)$ is concave.

  For fixed $\alpha$ let $f_n\in L^2([0,h])$ be the $n$-th eigenfunction of
  $\oper$. Define $g_n\in L^2([0,1])$ by $g_n(y)=f_n(hy)$. By the chain rule,
  \begin{equation}\label{eqn:relat}
    H(1,h^2\alpha)g_n=\frac{\dup^4g_n}{\dup y^4}-2h^2\alpha\frac{\dup^2g_n}{\dup
    y^2}=h^4\sigma(h,\alpha,n)g_n,
  \end{equation}
  so $g_n$ is an eigenfunction of $H(1,h^2\alpha)$ with eigenvalue
  $h^4\sigma(h,\alpha,n)$. It follows that
  $h^4\sigma(h,\alpha,n)\geq\rho_n(h^2\alpha)$. By a similar reverse argument
  we obtain equality.

  The family $H(1,\alpha)$ of differential operators indexed by $\alpha$ is a
  holomorphic family of type B and so $\rho_n$ are analytic except where the
  eigenvalues swap. See~\cite[chapter VII \S 4]{Kato}. We shall see in
  part (ii) that swapping does not occur when $\alpha$ is positive.\finbox

\textbf{Proof of theorem~\ref{thm:analysis} (ii):}
  For $a>0$ define $G:[0,1]^2\rightarrow\real$ by
  \begin{equation}
    G(x,y)=
    \begin{cases}
      \frac1ck(1-x,y),&y\leq x\\
      \frac1ck(x,1-y),&x<y
    \end{cases}
  \end{equation}
  where
  \begin{equation}
    c=a^3(2(1-\cosh a)+a\sinh a)
  \end{equation}
  and
  \begin{align}
    k(x,y)=&\phantom-(\sinh a-a)(\cosh ax-1)(\cosh ay-1)\notag\\
    &-(\cosh a-1)(\cosh ax-1)(\sinh ay-ay)\notag\\
    &-(\cosh a-1)(\sinh ax-ax)(\cosh ay-1)\notag\\
    &+\sinh a(\sinh ax-ax)(\sinh ay-ay).
  \end{align}
  For $g\in L^2([0,1])$ define $f$ by
  \[ f(x)=\int_0^1G(x,y)g(y)\dup y. \]
  Putting
  \[  k^{(r,s)}(x,y)=\frac{\partial^{r+s}k}{\partial x^r\partial y^s}(x,y), \]
  the identities
  \begin{alignat}{3}
    k&(1-x,x)-&k&(x,1-x)\equiv0&\qquad k^{(4,0)}(x,y)&\equiv a^2k^{(2,0)}(x,y
    )\notag\\
    -k^{(1,0)}&(1-x,x)-&k^{(1,0)}&(x,1-x)\equiv0&&\notag\\
    k^{(2,0)}&(1-x,x)-&k^{(2,0)}&(x,1-x)\equiv0&k(0,y)&\equiv0\notag\\
    -k^{(3,0)}&(1-x,x)-&k^{(3,0)}&(x,1-x)\equiv c&k^{(1,0)}(0,y)&\equiv0
  \end{alignat}
  imply that
  \[ f(0)=f'(0)=f(1)=f'(1)=0 \]
  and
  \[ \frac{\dup^4f}{\dup x^4}-a^2\frac{\dup^2f}{\dup x^2}=g, \]
  so $G$ is the Green's function of the operator
  \[ \frac{\dup^4\phantom x}{\dup x^4}-a^2\frac{\dup^2\phantom x}{\dup x^2} \]
  acting in $L^2([0,1])$ with Dirichlet boundary conditions. For positivity of 
  $G(x,y)$ we need to to prove that $k(x,y)$ is positive in the triangular 
  region 
  \begin{equation}
    T=\{(x,y):x,y>0\text{ and }x+y<1\},
  \end{equation}
  In order to do this, we introduce the function
  \begin{equation}
    \phi(x)=\tanh^{-1}\left(\frac{\sinh x-x}{\cosh x-1}\right).
  \end{equation}
  The reader should verify that
  \[ \lim_{x\rightarrow0^+}\phi'(x)=\frac13 \]
  and
  \[ \phi''(x)=\int_0^{x/2}4\sinh 2y[\sinh^2y\tanh y-y^3]\dup y\bigg/\left(\int
     _0^xy(\cosh y-1)\dup y\right)^2>0 \]
  for $x>0$. These features imply that $\phi$ is convex and increasing for 
  positive $x$. Let $(x,y)\in T$. Then
  \begin{equation}\label{eqn:phi_ineq}
    \phi(a) > \phi(ax)+\phi(ay).
  \end{equation}
  Hence using inequality~\eqref{eqn:phi_ineq} and a two angle $\tanh$
  identity,
  \begin{align*}
    \lefteqn{\frac{\sinh a-a}{\cosh a-1}+\frac{\sinh a-a}{\cosh a-1}\fgap\frac{
    \sinh ax-ax}{\cosh ax-1}\fgap\frac{\sinh ay-ay}{\cosh ay-1}}\qquad & \\
    &=\tanh\{\phi(a)\}[1+\tanh\{\phi(ax)\}\tanh\{\phi(ay)\}]\\
    &>\tanh\{\phi(ax)+\phi(ay)\}[1+\tanh\{\phi(ax)\}\tanh\{\phi(ay)\}]\\
    &=\tanh\{\phi(ax)\}+\tanh\{\phi(ay)\}\\
    &=\frac{\sinh ax-ax}{\cosh ax-1}+\frac{\sinh ay-ay}{\cosh ay-1}.
  \end{align*}
  Thus
  \begin{align*}
    k(x,y)=&\quad\big[(\sinh a-a)(\cosh ax-1)(\cosh ay-1) \\
    &-(\cosh a-1)(\cosh ax-1)(\sinh ay-ay) \\
    &-(\cosh a-1)(\sinh ax-ax)(\cosh ay-1) \\
    &+(\sinh a-a)(\sinh ax-ax)(\sinh ay-ay)\big] \\
    &+a(\sinh ax-ax)(\sinh ay-ay)>0
  \end{align*}
  as required.

  It follows that the Green's function of $H(h,\alpha)$ for all $h,\alpha>0$
  is positive because of the relationship~\eqref{eqn:relat} between $H(h,
  \alpha)$ and $H(1,h^2\alpha)$ established in part (i).

  See~\cite[theorem XIII.44]{ReedSimo4} for a proof that the groundstate 
  of $H(h,\alpha)$ is positive and the associated eigenvalue is multiplicity
  one.\finbox

\textbf{Proof of theorem~\ref{thm:analysis} (iii):}
  For $\alpha>0$ let $f=f_\alpha\in W_0^{2,2}([0,h])\cap C^\infty$ satisfy
  \begin{equation}\label{eqn:wibble}
    \frac{\dup^4f}{\dup x^4}-2\alpha\frac{\dup^2f}{\dup x^2}=\sigma(h,\alpha,n)f
  \end{equation}
  and have unit $L^2$-norm. Then $f_\alpha$ is a critical point of the
  functional
  \begin{equation}
    \mathcal E_\alpha(\psi)=\bip{\frac{\dup^2\psi}{\dup x^2}}{\frac{\dup^2\psi}
    {\dup x^2}}-2\alpha\bip{\frac{\dup^2\psi}{\dup x^2}}\psi
  \end{equation}
  in the sense that if $\psi(t)$ is $C^1$ with respect to $t$,
  $\|\psi(t)\|=1$ for all $t$ and $\psi(0)=f_\alpha$ then
  \begin{equation*}
    \left.\frac{\dup\phantom t}{\dup t}\mathcal E_\alpha(\psi(t))\right|_{t=0}
     =0.
  \end{equation*}
  Letting $\psi(t)=f_{\alpha+t}$ and using the relationship~\eqref{eqn:rhosig} 
  between $\rho$ and $\sigma$,
  \begin{align*}
    \rho_n'(h^2\alpha)h^{-2}&=\frac{\dup\phantom\alpha}{\dup\alpha}\sigma(h,
    \alpha,n)\\
    &=\frac{\dup\phantom\alpha}{\dup\alpha}\mathcal E_\alpha(f_\alpha)\\
    &=\left.\frac{\dup\phantom t}{\dup t}\mathcal E_t(f_\alpha)\right|_{t=
    \alpha}+\left.\frac{\dup\phantom t}{\dup t}\mathcal E_\alpha(\psi(t))\right
    |_{t=0}\\
    &=-2\bip{\frac{\dup^2f}{\dup x^2}}f_h+0.\lowearlybox
  \end{align*}

\textbf{Proof of theorem~\ref{thm:analysis} (iv):}
  Since $H(1,\alpha)$ is positive for $\alpha$ positive, eigenfunctions may
  be found by solving the auxiliary equation $y^4-2\alpha y^2-\rho=0$ where
  $\rho>0$. Since
  $\alpha-\sqrt{\alpha^2+\rho}<0<\alpha+\sqrt{\alpha^2+\rho}$, there are four
  distinct roots, two real and two imaginary. Denoting these roots by
  $a,-a,ib,-ib$, we see that $a^2-b^2=2\alpha$ and $a^2b^2=\rho$.

  There exists a combination of the functions $\cosh ax,\sinh ax,\cos bx,\sin 
  bx$ which satisfy the boundary conditions of $H(1,\alpha)$ if and only if
  \[ \det
     \begin{pmatrix}
       1&0&1&0\\0&a&0&b\\
       \cosh a&\sinh a&\cos b&\sin b\\
       a\sinh a&a\cosh a&-b\sin b&b\cos b\\
     \end{pmatrix}
     =0.\]
  Simplifying this determinant, the equation becomes
  \[ 2ab\cosh a\cos b-(a^2-b^2)\sinh a\sin b=2ab, \]
  so $\rho=\beta^2-\alpha^2$ is an eigenvalue of the operator if and only if
  \begin{equation*}
    \cosh\sqrt{\beta+\alpha}\cos\sqrt{\beta-\alpha}-\frac\alpha{\sqrt{\beta
    ^2-\alpha^2}}\sinh\sqrt{\beta+\alpha}\sin\sqrt{\beta-\alpha}=1. \lowearlybox
  \end{equation*}%
\textbf{Proof of theorem~\ref{thm:analysis} (v):}
  The proof of asymptotic formulae for $\rho$ and $\rho'$ as
  $\alpha\rightarrow\infty$ is given in the appendix. It is 
  possible to find the asymptotic formula for $\rho$ as $\alpha\rightarrow0$ by 
  a similar method. The formula for $\rho'$ then follows by differentiation 
  because $\rho$ is analytic at $0$. Here we give a sketch of an alternative 
  proof of the case $\alpha\rightarrow0$ for interest:

  Let $\gamma_n=\cosh c_n-\cos c_n$ and $\delta_n=\sinh c_n-\sin c_n$. Define
  \begin{equation}\label{eqn:beammode}
    g_n(x)=\cosh c_nx-\frac{\gamma_n}{\delta_n}\sinh c_nx-\cos c_nx+\frac{
    \gamma_n}{\delta_n}\sin c_nx.
  \end{equation}
  We claim that $(g_n)_{n\in\natural}$ is an orthonormal sequence of
  eigenfunctions of $H(1,0)$, and so $\rho_n(0)=c_n^4$. These eigenfunctions
  represent the fundamental modes of the clamped beam. We claim
  that, $\|g_n'\|_2^2=d_nc_n^2$, and so from~\eqref{eqn:norm_f'} we see that
  $\rho_n'(0)=2d_nc_n^2$.

  Verification of each of these claims is not trivial; indeed a lengthy 
  calculation is needed even to establish that $\|g_n\|_2=1$ for each $n$. This 
  task is left to the reader.\finbox

\textbf{Proof of theorem~\ref{thm:analysis} (vi):}
  Since $\rho_1$ is concave, $\rho_1'$ is decreasing. The result now follows
  immediately from part (v).\finbox%
\section{Lower Bounds on Eigenvalues}\label{sec:lower}%
We find lower bounds $\lambda_1(h),\lambda_3(h)$ for $\mu_1(h),\mu_3(h)$
respectively, by elementary algebraic manipulation of the biharmonic operator.
This method would be referred to as a finite renormalisation procedure in the
physics literature.

\begin{thm}{Theorem}\label{thm:lowerbd}
  \begin{alignat}{3}
    \lambda_1(h)&:=&&\hspace{0.22ex}\rho_1(\pi^2h^2)h^{-4}+\rho_1(\pi^2h^{-2}
    )-2\pi^4h^{-2}&&\leq\mu_1(h),\label{eqn:l1}\\*
    \lambda_3(h)&:=\;&\min\bigg\{&\frac[0pt]{\rho_1(\pi^2h^2)h^{-4}+\rho_2(
    \pi^2h^{-2})-2\pi^4h^{-2},}{\rho_3(\pi^2h^2)h^{-4}+\rho_1(\pi^2h^{-2})-2
    \pi^4h^{-2}\phantom{,}}\bigg\}&&\leq\mu_3(h).\label{eqn:l3}
  \end{alignat}
  where $\rho_n$ are defined in theorem~\slref{thm:analysis}.
  \end{thm}

\begin{proof}
  In this proof we consider, where relevant, restrictions of operators to
  $C_c^\infty([0,h])$. Let $A_h$ denote the biharmonic operator acting in
  $L^2([0,h])$, and let $B_h$ denote the Dirichlet Laplacian acting in
  $L^2([0,h])$. Then
  \begin{align*}
    \Delta^2&=A_h\otimes1_1+1_h\otimes A_1+2B_h\otimes B_1\\
    &=(A_h+2\pi^2B_h)\otimes1_1+1_h\otimes(A_1+2\pi^2h^{-2}B_1)\\
    &\qquad+2(B_h-\pi^2h^{-2}1_h)\otimes(B_1-\pi^2 1_1)-2\pi^4h^{-2}1_h
    \otimes1_1.
  \end{align*}
  The operator $(B_h-\pi^2h^{-2}1_h)\otimes(B_1-\pi^2 1_1)$ has eigenvalues
  \[ \pi^4h^{-2}(m^2-1)(n^2-1)\qquad m,n\in\natural, \]
  with the corresponding complete orthonormal sequence of eigenfunctions
  \[ 2h^{-1/2}\sin m\pi h^{-1} x\sin n\pi y, \]
  and so is non-negative. Hence
  \[ \bihar\geq(A_h+2\pi^2B_h)\otimes1_1+1_h\otimes(A_1+2\pi^2h^{-2}B_1)-2\pi
     ^4h^{-2}1_h\otimes1_1. \]
  It now follows that for $h$ large enough
  \begin{align*}
    \mu_n(h)&\geq\sigma(h,\pi^2,n)+\sigma(1,\pi^2h^{-2},1)-2\pi^4h^{-2}\\
    &=\rho_n(\pi^2h^2)h^{-4}+\rho_1(\pi^2h^{-2})-2\pi^4h^{-2}\\
    &=\lambda_n(h). \earlybox
  \end{align*}%
\end{proof}%
The functions $\lambda_n$ may be explicitly calculated using
formulae~\eqref{eqn:trans}, \eqref{eqn:l1} and \eqref{eqn:l3}. See
figures~\ref{fig:l1l3}, \ref{fig:l1n1} and table~\ref{tab:lowupp}. 

\begin{note}\label{not:evalveer}
  Despite the fact that the lower bounds $\lambda_n$ cross each other, swapping 
  of eigenvalues is a not actually a genuine feature of the increasing size of 
  the rectangle. In~\cite{BehnMert}, Behnke and Mertins show that the 
  eigenvalues veer away from each other just before the points where one might 
  expect them to cross.\finbox%
\end{note}%
\sngldiagm{%
  \renewcommand\afigwidth{67}%
  \renewcommand\afigheight{42}}{%
  \fig{l1l3}{\label{fig:l1l3}}{%
  \put(50,64){$\lambda_3(h)$}%
  \put(80,36){$\lambda_2(h)$}%
  \put(40,24){$\lambda_1(h)$}%
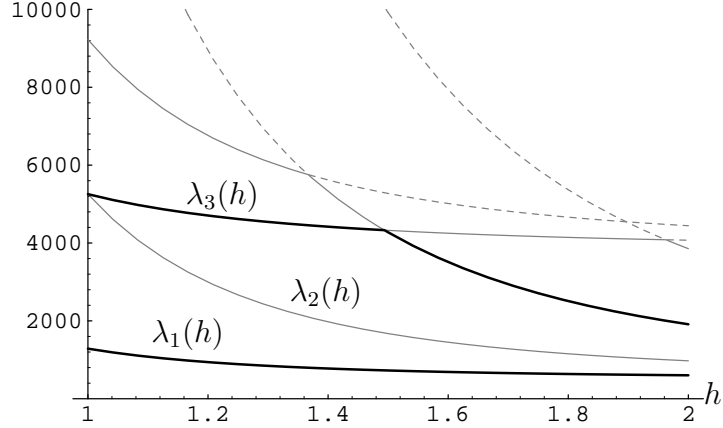
  \put(200,6){$h$}}{Graph of $\lambda_1(h)$ and $\lambda_3(h)$}}%
\section{Upper Bound on the First Eigenvalue}\label{sec:upper}%
We shall find an upper bound $\nu_1(h)$ for $\mu_1(h)$ by approximation of
the groundstate eigenfunction with separable functions. This is called the
Hartree method in the physics/chemistry literature and was used in~\cite{Davi6}
to estimate the groundstate energy of a somewhat similar but simpler problem.
Define the functional
\[ \mathcal E:L^2([0,h])\times L^2([0,1])\rightarrow\real \]
by
\begin{align}
  \mathcal E(f,g)&=Q(f\otimes g)\notag\\
  &=
  \begin{cases}
    {\displaystyle\int_0^h\!\int_0^1|\Delta(f(x)g(y))|^2\dup x\dup y,}&f\in W^{
    2,2}_0([0,h])\textup{and }g\in W^{2,2}_0([0,1]);\\
    \infty,&\textup{otherwise.}
  \end{cases}%
\end{align}%
Let
\begin{equation}\label{eqn:def_of_nu}
  \nu_1(h)=\inf\{\mathcal E(f,g):f\in L^2([0,h]),g\in L^2([0,1]),\|f\|_2=\|g\|
  _2=1\},
\end{equation}
be our upper bound on $\mu_1(h)$, where the two norms are taken in $L^2([0,h])$ 
and $L^2([0,1])$ respectively. The remainder of this section is devoted to 
giving a characterisation of $\nu_1$, and finally an implicit formula.

\begin{thm}{Lemma}
  The infimum in the expression~\sleqref{eqn:def_of_nu} for $\nu_1$ is
  attained.
\end{thm}

\begin{proof}
  Let $A_h$ denote the biharmonic operator acting in $L^2([0,h])$. Since $A_h$
  has compact resolvent, it has a complete orthonormal sequence of
  eigenfunctions whose corresponding eigenvalues form a divergent
  sequence. Using this orthonormal sequence we see that the set
  \[ S_h:=\{f\in L^2([0,h]):\|f\|_2^2+Q_h(f)\leq1\} \]
  is compact. The set
  \begin{align}
    S:=&{\{f\in L^2([0,h]):\|f\|_2^2+Q_h(f)\leq1,\|f\|_2^2\geq1/(2\nu_1+1)\}}
    \notag\\
    &\times{\{g\in L^2([0,1]):\|g\|_2^2+Q_1(g)\leq1,\|g\|_2^2\geq1/(2\nu_1+1)
    \}}
  \end{align}
  is a closed subset of $S_h\times S_1\subseteq L^2(R_h)$, separated from the
  origin, and so is a compact subset of $L^2(R_h)\setminus\{0\}$.
  Since the map
  \begin{equation}\label{eqn:map}
    f\mapsto Q(f)/\|f\|_2^2
  \end{equation}
  is lower semicontinuous on $L^2(R_h)\setminus\{0\}$ it attains its infimum
  when restricted to $S$.

  Suppose that $\mathcal E(f,g)\leq2\nu_1$ where $f$ and $g$ have unit norm.
  We may rescale $f$ and $g$ so that $\|f\|_2^2+Q_h(f)=1$ and
  $\|g\|_2^2+Q_1(g)=1$. Now
  \[ \frac{1-\|f\|_2^2}{\|f\|_2^2}=\frac{Q_h(f)}{\|f\|_2^2}\leq\frac{\mathcal E
     (f,g)}{\|f\|_2^2\|g\|_2^2}\leq2\nu_1, \]
  so $\|f\|_2^2\geq1/(2\nu_1+1)$. A similar argument for $g$ shows that 
  $(f,g)\in S$. It follows that
  \begin{equation*}
    \min\left\{\frac{\mathcal E(f,g)}{\|f\|_2\|g\|_2}:(f,g)\in S\right\}\leq\nu
    _1(h).
  \end{equation*}
  Let $(\tilde f,\tilde g)\in S$ take the minimum value of the
  map~\eqref{eqn:map} restricted to $S$. We may rescale so that
  $\|\tilde f\|_2=\|\tilde g\|_2=1$. Now
  \begin{equation*}
    \nu_1(h)\leq\mathcal E(\tilde f,\tilde g)=\min\left\{\frac{\mathcal E(f,g)}
    {\|f\|_2\|g\|_2}:(f,g)\in S\right\}\leq\nu_1(h). \lowearlybox
  \end{equation*}%
\end{proof}%
From this point onwards, we shall assume that $\mathcal E$ is only applied to
functions of unit norm. We may rewrite the formula for $\mathcal E$ as
\begin{equation}\label{eqn:simplE}
  \mathcal E(f,g)=\bip{\frac{\dup^2f}{\dup x^2}}{\frac{\dup^2f}{\dup x^2}}_h+2
  \bip{\frac{\dup^2f}{\dup x^2}}f_h\bip{\frac{\dup^2g}{\dup y^2}}g_1+\bip{\frac
  {\dup^2g}{\dup y^2}}{\frac{\dup^2g}{\dup y^2}}_1,
\end{equation}
where $\ip.._h$ denotes the inner product on $L^2([0,h])$. We shall identify 
$f$ and $g$ by searching for the critical points of $\mathcal E$.

\begin{thm}{Lemma}\label{thm:2}
  Let $f\in W^{2,2}_0([0,h])$ and $g\in W^{2,2}_0([0,1])$ minimise
  $\mathcal E(f,g)$. Then $f\in C^\infty([0,h])\cap W^{2,2}_0([0,h])$ is the
  groundstate of the operator $H(h,\alpha_g)$, where
  \begin{equation}
    \alpha_g:=\|g'\|_2^2.
  \end{equation}
\end{thm}

\begin{proof}
  Let $\tilde f\in C^\infty_c([0,h])$ be such that $\ip f{\tilde f}_h=0$.
  Define
  \begin{equation}
    f(t):=\frac{f+t\tilde f}{\|f+t\tilde f\|}.
  \end{equation}
  By differentiating we see that
  \[ f(0)=f \text{ and } \frac{\dup f}{\dup t}(0)=\tilde f.\]
  The minimum of $\mathcal E$ will be a critical point so
  \[ 0=\left.\frac{\dup\mathcal E}{\dup t}(f(t),g)\right|_{t=0}=2\re\bip
     {\frac{\dup^4f}{\dup x^4}-2\alpha_g\frac{\dup^2f}{\dup x^2}}{\tilde f}_h \]
  where the fourth derivative has been taken in the distributional sense.
  Replacing $\tilde f$ by $i\tilde f$ we see that
  \begin{equation}\label{eqn:oper_ip}
    \bip{\frac{\dup^4f}{\dup x^4}-2\alpha_g\frac{\dup^2f}{\dup x^2}}{\tilde f}_h
    =0 
  \end{equation}
  for all such $\tilde f$. It now follows that
  \begin{equation}\label{eqn:f_is_efunc}
    \frac{\dup^4f}{\dup x^4}-2\alpha_g\frac{\dup^2f}{\dup x^2}=\mu f\qquad\mu\in
    \real.
  \end{equation}
  For suppose otherwise, then there exist $f_1,f_2\in C_c^\infty$ with
  $\ip f{f_1}_h\not=0$ and $\ip f{f_2}_h\not=0$, and
  $\mu_1\not=\mu_2$ such that
  \[ \bip{\frac{\dup^4f}{\dup x^4}-2\alpha_g\frac{\dup^2f}{\dup x^2}}{f_i}_h
     =\mu_i\ip f{f_i}_h,\qquad i=1,2. \]
  Now let $\tilde f=\ip f{f_2}_h f_1-\ip f{f_1}_h f_2$. Then
  \[ \ip f{\tilde f}_h=0 \]
  and
  \begin{align*}
    \bip{\frac{\dup^4f}{\dup x^4}-2\alpha_g\frac{\dup^2f}{\dup x^2}}{\tilde f}_h
    &= \ip f{f_2}_h\mu_1\ip f{f_1}_h-\ip f{f_1}
    _h\mu_2\ip f{f_2}_h\\
    &=(\mu_1-\mu_2)\ip f{f_1}_h\ip f{f_2}_h\\
    &\not=0,
  \end{align*}
  which contradicts~\eqref{eqn:oper_ip}.

  It follows from equation~\eqref{eqn:f_is_efunc} that $f\in W^{m,2}$ for all
  $m\in\natural$, and consequently $f$ is smooth. From~\eqref{eqn:simplE}
  the value of $\mathcal E$ at the critical point is
  \[ \mu+\bip{\frac{\dup^2g}{\dup y^2}}{\frac{\dup^2g}{\dup y^2}}_1, \]
  In order for this to be the minimum, $f$ must be the groundstate of
  $H(h,\alpha_g)$.\finbox
\end{proof}

The following theorem gives a formula which characterises the abstract 
definition~\eqref{eqn:def_of_nu} of $\nu_1$.

\begin{thm}{Theorem} \label{thm:upperbd}
  Let $\nu_1$ be the upper bound on $\mu_1$ defined 
  by~\sleqref{eqn:def_of_nu}. Then
  \begin{equation}\label{eqn:form_nu} \nu_1(h)=\rho_1(h^2\alpha)h^{-4}+\rho_1
    \left(\frac12\rho_1'(h^2\alpha)h^{-2}\right)-\rho_1'(h^2\alpha)h^{-2}
    \alpha.
  \end{equation}
  where $\alpha\in[\pi^2,dc^2]$ is a solution of the equation
  \begin{equation}\label{eqn:alpha}
    2\alpha=\rho_1'\left(\frac12\rho_1'(h^2\alpha)h^{-2}\right).
  \end{equation}%
\end{thm}%
\begin{proof}
  Let $(f,g)$ be a pair of functions which attain the minimum of $\mathcal E$. 
  Then by lemma~\ref{thm:2}, $f$ is the groundstate of the operator 
  $H(h,\alpha_g)$. Using~\eqref{eqn:norm_f'} we see that
  \begin{equation}\label{eqn:alftoalg}
    \alpha_f=\frac12\rho_1'(h^2\alpha_g)h^{-2}.
  \end{equation}
  By an identical argument,
  \begin{equation}
    \alpha_g=\frac12\rho_1'(\alpha_f),
  \end{equation}
  so
  \[ 2\alpha_g=\rho_1'\left(\frac12\rho_1'(h^2\alpha_g)h^{-2}\right). \]
  Using equations~\eqref{eqn:def_of_nu}, \eqref{eqn:simplE}, 
  \eqref{eqn:wibble}, \eqref{eqn:rhosig}, and~\eqref{eqn:alftoalg} we see that
  \begin{align*}
    \nu_1(h)&=\bip{\frac{\dup^4f}{\dup x^4}}f_h+2\bip{\frac{\dup^2f}{\dup x^2}}
    f_h\bip{\frac{\dup^2g}{\dup y^2}}g_1+\bip{\frac{\dup^4g}{\dup y^4}}g_1\\
    &=\sigma(h,\alpha_g,1)-2\alpha_f\alpha_g+2\alpha_f\alpha_g+\rho_1(\alpha
    _f)-2\alpha_f\alpha_g\\
    &=\rho_1(h^2\alpha_g)h^{-4}-\rho_1'(h^2\alpha_g)h^{-2}\alpha_g+\rho_1
    \left(\frac12\rho_1'(h^2\alpha_g)h^{-2}\right).
  \end{align*}
  Inequality~\eqref{eqn:dr_ineq} implies that
  \[ \rho_1'\left(\frac12\rho_1'(h^2\alpha)h^{-2}\right)\geq2\pi^2>2\alpha \]
  for $\alpha<\pi^2$ and
  \[ \rho_1'\left(\frac12\rho_1'(h^2\alpha)h^{-2}\right)\leq2dc^2<2\alpha \]
  for $\alpha>dc^2$, so there is at least one solution of
  equation~\eqref{eqn:alpha} in the interval $[\pi^2,dc^2]$, and there are no
  solutions outside. Hence $\alpha_g\in[\pi^2,dc^2]$.\finbox%
\end{proof}

\dbledgram{%
  \renewcommand\afigwidth{47}%
  \renewcommand\afigheight{42}%
  \renewcommand\bfigwidth{38}%
  \renewcommand\bfigheight{36}}{%
  \fig{l1n1}{\label{fig:l1n1}}{%
    \put(81,89){$\nu_1(h)$}%
    \put(54,69){$\lambda_1(h)$}%
    \put(201,10){$h$}}{Graph of $\lambda_1(h)$ and $\nu_1(h)$}}{%
  \tab{\label{tab:lowupp}}{c|c|c|c}{%
    \hline%
    $h$ & $\lambda_1$ & $\nu_1$ & \parbox{2.5em}{\vspace{0.5ex}\% err on %
    $\mu_1$\vspace{0.5ex}} \\\hline%
    1.0 & 1286.66 & 1295.93 & 0.720 \\%
    1.2 & 940.070 & 946.421 & 0.676 \\%
    1.4 & 776.088 & 780.618 & 0.584 \\%
    1.6 & 687.796 & 691.129 & 0.485 \\%
    1.8 & 635.529 & 638.044 & 0.396 \\%
    2.0 & 602.282 & 604.221 & 0.322 \\%
    3.0 & 537.444 & 538.111 & 0.124 \\%
    4.0 & 519.496 & 519.794 & 0.058 \\%
    5.0 & 512.080 & 512.237 & 0.031 \\\hline}%
    {Values of $\lambda_1$, $\nu_1$ and a\par guaranteed error estimate on %
     $\mu_1$}}%
\begin{note}
  Numerical evidence strongly suggests that there is a unique solution of
  equation~\eqref{eqn:alpha} for any value of $h$. Plots of
  $\rho_1''(\alpha)$ suggest that
  \[ -0.0603<\rho_1''(\alpha)<0 \]
  for all positive $\alpha$. The fact that $\rho_1''(\alpha)$ is negative is
  clear because we showed that $\rho$ is concave in
  theorem~\ref{thm:analysis} (i). A unique solution of
  equation~\eqref{eqn:alpha} is guaranteed however under the weaker
  requirement that
  \[ -2<\rho_1''(\alpha)<0 \]
  for all positive $\alpha$. For then
  \begin{align*}
    \frac{\dup\phantom\alpha}{\dup\alpha}(2\alpha)&=2<\frac12\rho_1''\left(
    \frac12\rho_1'(h^2\alpha)h^{-2}\right)\rho_1''(h^2\alpha)\\
    &=\frac{\dup\phantom\alpha}{\dup\alpha}\left(\rho_1'\left(\frac12\rho_1'(h^
    2\alpha)h^{-2}\right)\right)
  \end{align*}
  and so
  \[ \rho_1'\left(\frac12\rho_1'(h^2\alpha)h^{-2}\right)-2\alpha \]
  is a strictly increasing function.\finbox
\end{note}

Figure~\ref{fig:l1n1} and table~\ref{tab:lowupp} have been created by 
\mathematica{} using formula~\eqref{eqn:l1} for $\lambda_1$, and 
formulae~\eqref{eqn:form_nu}, \eqref{eqn:alpha} for $\nu_1$. For every value of 
$h$ taken there was, as expected, only one solution of 
equation~\eqref{eqn:alpha}. As in~\cite{Davi6}, the Har\-tree 
approximation $f\otimes g$ gives the worst approximation when there is an extra 
rotational symmetry, as for $h=1$.

In theorem~\ref{thm:asymu} we show that both $\lambda_1$ and $\nu_1$
converge to $c^4\approx500.564$ and the guaranteed percentage error
\[ \left(\frac{\nu_1-\lambda_1}{\lambda_1}\right)\times100 \]
on $\mu_1$ is of order $h^{-3}$. This allows us to prove the asymptotic
formula~\eqref{eqn:asymu1} up to the same order.%
\section{The Groundstate}\label{sec:gstate}%
\begin{thm}{Theorem}\label{thm:7}
  Let $f_1^-$ be the negative part of the groundstate $f_1$ of $\bihar$ acting 
  in $L^2(R_h)$. Then
  \begin{equation}\label{eqn:l2_bound}
    \frac{\fmtwo}{\ftwo}_{\phantom2}<\frac{(\nu_1-\lambda_1)^{1/2}}{(
    \lambda_3-\nu_1)^{1/2}}
  \end{equation}
  where $\lambda_1(h)$, $\lambda_3(h)$ and $\nu_1(h)$ are given by
  theorem~\slref{thm:lowerbd} and
  equation~\sleqref{eqn:def_of_nu}.
\end{thm}

We give an asymptotic expansion of the above bound in
corollary~\ref{thm:negbdexp}. It is also possible to evaluate this bound for 
smaller values of $h$ giving the results in figure~\ref{fig:posit} and 
table~\ref{tab:posit}.

\begin{proof}
  There exists a complete orthonormal sequence of eigenfunctions $f_n$ of the
  biharmonic operator acting in $L^2(R_h)$ with corresponding eigenvalues
  $\mu_n$ listed in increasing order. Let $(f,g)$ minimize $\mathcal{E}$. Then 
  \begin{equation}
    \|f\otimes g\|_2=1\text{ and }\lambda_1\leq\mu_1=\|\Delta f_1\|_2^2<\|
    \Delta(f\otimes g)\|_2^2=\nu_1.
  \end{equation}
  By writing $f\otimes g=\sum_{n=1}^\infty\alpha_nf_n$ with $\alpha_1>0$, we
  see that
  \begin{equation}
    \sum_{n=1}^\infty\alpha_n^2=1,\quad\alpha_2=0\text{ and }\lambda_1\leq\mu
    _1<\sum_{n=1}^\infty\alpha_n^2\mu_n=\nu_1.
  \end{equation}
  Using the above information,
  \begin{align}\label{eqn:1-al^2}
    \nu_1-\mu_1&=\sum_{n=3}^\infty
    \alpha_n^2\mu_n-(1-\alpha_1^2)\mu_1\notag\\
    &\geq\lambda_3\sum_{n=3}^\infty\alpha_n^2-(1-\alpha_1^2)\mu_1\notag\\
    &=(\lambda_3-\mu_1)(1-\alpha_1^2).
  \end{align}
  Rearranging~\eqref{eqn:1-al^2}, we see that
  \begin{equation*}
    1+\alpha_1=2-(1-\alpha_1)\geq2-\frac{\nu_1-\mu_1}{\lambda_3-\mu_1}
    >\frac{2(\lambda_3-\nu_1)}{\lambda_3-\mu_1}.
  \end{equation*}
  Using~\eqref{eqn:1-al^2} again,
  \begin{equation}
    1-\alpha_1<\frac{\nu_1-\lambda_1}{2(\lambda_3-\nu_1)}.
  \end{equation}
  Now
  \begin{align}
    \|f\otimes g-f_1\|_2^2&=(1-\alpha_1)^2+\sum_{n=3}^\infty\alpha_n^2\notag
    \\
    &=\sum_{n=1}^\infty\alpha_n^2+(1-2\alpha_1)\notag\\
    &=2(1-\alpha_1)\notag\\
    &<\frac{\nu_1-\lambda_1}{\lambda_3-\nu_1}.
  \end{align}
  Using theorem~\ref{thm:analysis} (ii) we see that $f\otimes g$ is positive. 
  It follows that
  \begin{equation*}
    |f_1^-|\leq|f\otimes g-f_1|. \earlybox
  \end{equation*}%
\end{proof}%
\dbledgram{%
  \renewcommand\afigwidth{55}%
  \renewcommand\afigheight{35}%
  \renewcommand\bfigwidth{21}%
  \renewcommand\bfigheight{29}}{%
  \fig{negpart}{\label{fig:posit}}{%
    \put(202,5){$h$}}%
    {Graph of $(\nu_1-\lambda_1)^{1/2}/(\lambda_3-\nu_1)^{1/2}$}}{%
  \tab{\label{tab:posit}}{c|c}{%
    \hline%
    $h$ & \parbox{4.1em}{\vspace{0.5ex}$\frac{(\nu_1-\lambda_1)%
    ^{1/2}}{(\lambda_3-\nu_1)^{1/2}}$\vspace{0.5ex}} \\\hline%
    1   & 0.0484 \\%
    10  & 0.0336 \\%
    20  & 0.0249 \\%
    40  & 0.0179 \\%
    60  & 0.0147 \\%
    80  & 0.0128 \\%
    100 & 0.0114 \\\hline}%
    {}}%
The following Sobolev embedding theorem is used to convert the $L^2$ bound
above to an $L^\infty$ bound.

\begin{thm}{Lemma}\label{thm:sobolev}
  Let $f\in W^{2,2}(\real^2)$. Then
  \begin{equation}
    \|f\|^2_\infty\leq\frac14\|f\|_2\|\Delta f\|_2.
  \end{equation}%
\end{thm}%
\begin{proof}
  Let $h\in L^2(\real^2)$ be the function whose Fourier transform is
  \begin{equation}
    \hat h(\xi)=(2\pi)^{-1}(\gamma+\gamma^{-1}|\xi|^4)^{-1/2}.
  \end{equation}
  Then
  \begin{equation}
    \|h\|_2^2=\|\hat h\|_2^2=\frac18.
  \end{equation}
  Now
  \begin{equation}
    (\gamma+\gamma^{-1}\Delta^2)^{-1/2}g=h*g
  \end{equation}
  for all $g\in L^2(\real^2)$, so
  \begin{equation}
    \|(\gamma+\gamma^{-1}\Delta^2)^{-1/2}g\|_\infty^2=\|h*g\|_\infty^2
    \leq\|h\|_2^2\|g\|_2^2=\frac18\|g\|_2^2.
  \end{equation}
  Hence for $f\in W^{2,2}(\real^N)$
  \[ \|f\|_\infty^2\leq\frac18\|(\gamma+\gamma^{-1}\Delta^2)^{1/2}f\|_2^2
     =\frac18(\gamma\|f\|_2^2+\gamma^{-1}\|\Delta f\|_2^2). \]
  The result is obtained by setting
  \begin{equation}
    \gamma=\frac{\|\Delta f\|_2}{\|f\|_2}. \lowearlybox
  \end{equation}%
\end{proof}%
\begin{thm}{Theorem} \label{thm:9}
  Using the notation of theorem~\slref{thm:7},
  \begin{equation}\label{eqn:linf_bound}
    \frac{\fminfty}{\ftwo}_{\phantom\infty}<\frac{{(\nu_1-\lambda_1)}^{1/2}
    \lambda_3^{1/4}}{2(\lambda_3-\nu_1)^{1/2}}.
  \end{equation}%
\end{thm}%
\begin{proof}
  Let $f,g$, and $f_1$ be as in theorem~\ref{thm:7}. Then
  \begin{align*}
    \|\Delta(f\otimes g-f_1)\|_2^2&=\sum_{n=1}^\infty\alpha_n^2\mu_n-\mu_1+2(1-
    \alpha_1)\mu_1 \\
    &<\nu_1-\lambda_1+\left(\frac{\nu_1-\lambda_1}{\lambda_3-\nu_1}\right)\nu_1
    \\
    &=\frac{(\nu_1-\lambda_1)}{(\lambda_3-\nu_1)}\lambda_3.
  \end{align*}
  Lemma~\ref{thm:sobolev} implies that
  \begin{align*}
    \|f\otimes g-f_1\|_\infty^4&\leq\frac1{16}\|f\otimes g-f_1\|_2^2\|\Delta(f
    \otimes g-f_1)\|_2^2 \\
    &<\frac{(\nu_1-\lambda_1)^2\lambda_3}{16(\lambda_3-\nu_1)^2}. \lowearlybox
  \end{align*}%
\end{proof}%
\section{Asymptotic Estimates}%
The functions $\lambda_1$ and $\nu_1$ become increasingly accurate estimates of 
$\mu_1$ as $h\rightarrow\infty$ (see figure~\ref{fig:l1n1} and 
table~\ref{tab:lowupp}). The following theorem is proved using the fact that 
$\lambda_1$ and $\nu_1$ have the same asymptotic formulae up to $O(h^{-3})$.

\begin{thm}{Theorem}\label{thm:asymu}
  The first eigenvalue of the biharmonic operator acting in $L^2(R_h)$ has
  the asymptotic formula
  \begin{equation}\label{eqn:asymu2}
    \mu_1(h)=c^4+2dc^2\pi^2h^{-2}+O(h^{-3})
  \end{equation}
  as $h\rightarrow\infty$.
\end{thm}

\begin{proof}
  Let $\alpha(h)$ be the solution of~\eqref{eqn:alpha} which
  makes~\eqref{eqn:form_nu} valid. It follows from
  inequality~\eqref{eqn:dr_ineq} that
  \[ \frac12\rho_1'(h^2\alpha(h))h^{-2}=O(h^{-2}) \]
  as $h\rightarrow\infty$. Hence by equation~\eqref{eqn:alpha},
  \begin{equation}
    \alpha(h)=dc^2+O(h^{-2})
  \end{equation}
  as $h\rightarrow\infty$. Substituting this expression into the 
  formula~\eqref{eqn:form_nu} for $\nu_1$ we see that
  \begin{align}\label{eqn:formn1}
    \nu_1(h)&=h^{-4}\rho_1(dc^2h^2+O(1))-h^{-2}\rho_1'(dc^2h^2+O(1))(dc^2+O(h
    ^{-2}))\notag\\
    &\qquad+\rho_1(\rho_1'(dch^2+O(1))/2h^2)\notag\\
    &=h^{-4}(2\pi^2dc^2h^2+4\sqrt2\pi^2d^{1/2}ch+O(1))\notag\\
    &\qquad-h^{-2}(2\pi^2+2\sqrt2\pi^2d^{-1/2}c^{-1}h^{-1}+O(h^{-2}))(dc
    ^2+O(h^{-2}))\notag\\
    &\qquad+\rho(\pi^2h^{-2}+\sqrt2\pi^2d^{-1/2}c^{-1}h^{-3}+O(h^{-4}))
    \notag\\
    &=2\pi^2dc^2h^{-2}+4\sqrt2\pi^2d^{1/2}ch^{-3}-2\pi^2dc^2h^{-2}-2\sqrt
    2\pi^2d^{1/2}ch^{-3} \notag\\
    &\qquad+c^4+2dc^2\pi^2h^{-2}+2\sqrt2\pi^2d^{1/2}ch^{-3}+O(h^{-4})
    \notag\\
    &=c^4+2dc^2\pi^2h^{-2}+4\sqrt2\pi^2d^{1/2}ch^{-3}+O(h^{-4}),
  \end{align}
  as $h\rightarrow\infty$.

  Also, by substituting the asymptotic formulae~\eqref{eqn:forms} for $\rho_1$
  into the formula~\eqref{eqn:l1} for $\lambda_1$, we see that
  \begin{align}\label{eqn:forml1}
    \lambda_1(h)&=c^4+2dc^2\pi^2h^{-2}+h^{-4}(2\pi^2h^2\pi^2+4\sqrt2\pi^2h\pi
    +O(1))-2\pi^4h^{-2}+O(h^{-4})\notag\\
    &=c^4+2dc^2\pi^2h^{-2}+4\sqrt2\pi^3h^{-3}+O(h^{-4}),
  \end{align}
  as $h\rightarrow\infty$.\finbox
\end{proof}

\begin{note}
  For long thin rectangles a good approximation to the groundstate of the 
  biharmonic operator is
  \[  \sqrt2h^{-1/2}\sin\left(\frac{\pi x}h\right)g_1(y) \]
  where $g_1$ is defined by formula~\eqref{eqn:beammode}. Intuitively one
  expects this function to be a fairly good approximation because the boundary
  conditions at the ends of the rectangle become less influential on the
  eigenfunction. Note that the above function does not actually lie in the
  quadratic form domain. The energy of this function may be computed however if
  we ignore this fact, and we see that
  \[  Q(f)=c^4+2dc^2\pi^2h^{-2}+\pi^4h^{-4}. \]
  This compares well with the asymptotic formula~\eqref{eqn:asymu2} for $\mu_1$
  as $h\rightarrow\infty$, differing only by terms of order $h^{-3}$.\finbox
\end{note}

\begin{thm}{Corollary}\label{thm:negbdexp}
  The bounds~\textup{\eqref{eqn:l2_bound}}
  and~\textup{\eqref{eqn:linf_bound}} have asymptotic formulae
  \begin{align}\label{eqn:asym22}
    \frac{\fmtwo}{\ftwo}_{\phantom\infty}&<\frac{(\nu_1-\lambda_1)^{1/2}}
    {(\lambda_3-\nu_1)^{1/2}}=\frac{2^{1/4}(d^{1/2}c-\pi)^{1/2}}{2\pi}h^{-1/2}+
    O(h^{-3/2})\\
    \intertext{and}\label{eqn:asyminf2}
    \frac{\fminfty}{\ftwo}_{\phantom\infty} &< \frac{(\nu_1-\lambda_1)^{1/2}
    \lambda_3^{1/4}}{2(\lambda_3-\nu_1)^{1/2}}=\frac{2^{1/4}(d
    ^{1/2}c-\pi)^{1/2}c}{4\pi}h^{-1/2}+O(h^{-3/2})
  \end{align}
  as $h\rightarrow\infty$.
\end{thm}

\begin{proof}
  The asymptotic formula of $\lambda_3$,
  \begin{equation}\label{eqn:forml3}
    \lambda_3(h)=c^4+(2dc^2\pi^2+16\pi^4)h^{-2}+O(h^{-3})
  \end{equation}
  as $h\rightarrow\infty$ is found by using formula~\eqref{eqn:l3} and the
  asymptotic formulae~\eqref{eqn:forms} for $\rho_1$ and $\rho_3$. The corollary
  follows by using the formulae~\eqref{eqn:formn1}, \eqref{eqn:forml1}
  and~\eqref{eqn:forml3}.\finbox
\end{proof}

\begin{note}\label{note}
  It is conjectured that
  \begin{equation}\label{eqn:asy2inf}
  \frac{\|f_1\|\makebox[0cm][l]{$_2$}}{\|f_1\|\makebox[0cm][l]{$_\infty$}}_{
  \phantom\infty}=\frac{\cosh\cbytwo\sinh\cbytwo-\cos\cbytwo\sin\cbytwo}{2
  \sqrt{2}\big[\cosh^2\cbytwo\sin\cbytwo+\sinh\cbytwo\cos^2\cbytwo\big]}h^{
  1/2}+O(h^{-1/2})
  \end{equation}
  as $h\rightarrow\infty$ so
  \begin{align}\label{eqn:asyminfinf}
    \frac{\fminfty}{\finfty}_{\phantom\infty}&=\frac{\fminfty}{\ftwo}_{
    \phantom\infty}\frac{\|f_1\|\makebox[0cm][l]{$_2$}}{\|f_1\|\makebox[0cm]
    [l]{$_\infty$}}_{\phantom\infty}\notag\\
    &<\frac{2^{1/4}(d^{1/2}c-\pi)^{1/2}c\big[\cosh\cbytwo\sinh
    \cbytwo-\cos\cbytwo\sin\cbytwo\big]}{8\sqrt{2}\pi\big[\cosh^2\cbytwo\sin
    \cbytwo+\sinh\cbytwo\cos^2\cbytwo\big]}+O(h^{-1})\notag\\
    &<0.121
  \end{align}
  for $h$ large enough. Comparing this expression with~\eqref{eqn:asym22} we
  see that bound~\eqref{eqn:linf_bound} is of some use, but is a lot weaker
  than~\eqref{eqn:l2_bound}. An improved bound would be desirable.\finbox
\end{note}%
\section{Appendix}%
\textbf{Proof of theorem~\ref{thm:analysis} (v):}

Preliminary calculations show that the roots $\beta_n$ of
equation~\eqref{eqn:trans} are of the form $\alpha+n^2\pi^2+o(1)$ as
$\alpha\rightarrow\infty$. Define
\[ \beta_+(\alpha)=\alpha+n^2\pi^2+2\sqrt2n^2\pi^2\alpha^{-1/2}+6n^2\pi
   ^2\alpha^{-1}+\frac16(-1)^n5\sqrt{2}n^4\pi^2\alpha^{-3/2}. \]
Then
\begin{align*}
  \lefteqn{\cos\sqrt{\beta_+-\alpha}-\frac\alpha{\sqrt{\beta_+^2-\alpha^2}}
  \tanh\sqrt{\beta_++\alpha}\sin\sqrt{\beta_+-\alpha}}
  \quad & \\
  &=\cos(n^2\pi^2+2\sqrt2n^2\pi^2\alpha^{-1/2}+O(\alpha^{-1}))^{1/2
  }\\
  &\quad-\frac{\sqrt2\alpha^{1/2}}{2n\pi}\left(1+2\sqrt2\alpha^{-1/2}+\frac{(12
  +n^2\pi^2)}2\alpha^{-1}+O(\alpha^{-3/2})\right)^{-1/2}\\
  &\qquad\quad\times\tanh(2\alpha+O(1))^{1/2}\\
  &\qquad\quad \times\sin\left(n^2\pi^2+2\sqrt2n^2\pi^2\alpha^{-1/2}+6n
  ^2\pi^2\alpha^{-1}+\frac16(-1)^n5\sqrt2n^4\pi^4\alpha^{-3/2}\right)^{1/2}\\
  &=(-1)^n(1-n^2\pi^2\alpha^{-1}+O(\alpha^{-3/2})) \\
  &\quad- \frac{\sqrt2\alpha^{1/2}}{2n\pi}(1-\sqrt2\alpha^{-1/2}-
  \frac14n^2\pi^2\alpha^{-1}+O(\alpha^{-3/2}))(1+O(\alpha^{-3/2})) \\
  &\qquad\quad\times (-1)^n\bigg(\sqrt2n\pi\alpha^{-1/2}+2n\pi\alpha^{-
  1}-2\sqrt2n\pi\alpha^{-3/2}-\frac13\sqrt2n^3\pi^3\alpha^{-3/2}\\
  &\qquad\qquad\qquad\qquad\qquad\qquad\qquad\quad+\frac1{12}(-1)^n5\sqrt2n
  ^3\pi^3\alpha^{-3/2}+O(\alpha^{-2})\bigg)\\
  &= (-1)^n\left[4-\frac1{12}5n^2\pi^2(1+(-1)^n)\right]\alpha^{-1}+O(\alpha
  ^{-3/2})\\
  &<\operatorname{sech}\sqrt{\beta_++\alpha}
\end{align*}
for $\alpha$ large enough.

Similarly, defining
\[ \beta_-(\alpha)=\alpha+n^2\pi^2+2\sqrt2n^2\pi^2\alpha^{-1/2}+6n^2\pi
   ^2\alpha^{-1}-\frac16(-1)^n5\sqrt2n^4\pi^4\alpha^{-3/2}, \]
we see that
\begin{align*}
  \lefteqn{\cos\sqrt{\beta_--\alpha}-\frac\alpha{\sqrt{\beta_-^2-\alpha^2}}
  \tanh\sqrt{\beta_-+\alpha}\sin\sqrt{\beta_--\alpha}} 
  \qquad\qquad\qquad & \\
  &=(-1)^n\left[4-\frac1{12}5n^2\pi^2(1-(-1)^n)\right]\alpha^{-1}+O(\alpha^
  {-3/2})\\
  &>\operatorname{sech}\sqrt{\beta_-+\alpha}
\end{align*}
for $\alpha$ large enough.

It follows that for $\alpha$ large enough, $\beta_n$ lies between $\beta_-$
and $\beta_+$ so
\begin{equation}\label{eqn:asy_beta}
  \beta_n(\alpha)=\alpha+n^2\pi^2+2\sqrt2n^2\pi^2\alpha^{-1/2}+6n^2\pi^
  2\alpha^{-1}+O(\alpha^{-3/2}).
\end{equation}
Thus
\begin{align*}
  \rho_n(\alpha)&=\beta_n(\alpha)^2-\alpha^2\\
  &=2n^2\pi^2\alpha+4\sqrt2n^2\pi^2\alpha^{1/2}+O(1).
\end{align*}
Let $F(\alpha,\beta)$ be the left hand side of equation~\eqref{eqn:trans}.
Then differentiating~\eqref{eqn:trans}, we see that
\begin{equation}\label{eqn:beta`}
  \beta_n'(\alpha)=-\frac{F_1(\alpha,\beta_n)}{F_2(\alpha,\beta_n)}=1-\left
  (\frac{F_1(\alpha,\beta_n)+F_2(\alpha,\beta_n)}{F_2(\alpha,\beta_n)}
  \right)
\end{equation}
By explicit differentiation of $F$ and substitution of the asymptotic
formula~\eqref{eqn:asy_beta} of $\beta_n$, we see that
\begin{align*}
  \lefteqn{\frac{2(F_1(\alpha,\beta_n)+F_2(\alpha,\beta_n))(\beta_n^2-\alpha^
  2)^{3/2}}{\cosh\sqrt{\beta_n+\alpha}\cos\sqrt{\beta_n-\alpha}}} \qquad
  & \\
  &=2(\beta_n-\alpha)(\beta_n+\alpha)^{1/2}\tanh\sqrt{\beta_n+\alpha}\\
  &\quad-2\alpha(\beta_n-\alpha)^{1/2}\tan\sqrt{\beta_n-\alpha}\\
  &\quad-2\beta_n(\beta_n-\alpha)^{1/2}(\beta_n+\alpha)^{-1/2}\tanh
  \sqrt{\beta_n+\alpha}\tan\sqrt{\beta_n-\alpha}\\
  \ldots&=-2n^2\pi^2+O(\alpha^{-1/2}).
\end{align*}
Also,
\begin{align*}
  \frac{2F_2(\alpha,\beta_n)(\beta_n^2-\alpha^2)^{3/2}}{\cosh\sqrt{
  \beta_n+\alpha}\cos\sqrt{\beta_n-\alpha}}&=(\beta_n-2\alpha)(\beta_n-
  \alpha)^{1/2}(\beta_n+\alpha)\tanh\sqrt{\beta_n+\alpha}\\
  &\quad-(\beta_n+2\alpha)(\beta_n-\alpha)(\beta_n+\alpha)^{1/2}\tan\sqrt{
  \beta_n-\alpha}\\
  &\quad+2\alpha\beta_n\tanh\sqrt{\beta_n+\alpha}\tan\sqrt{\beta_n-\alpha}\\
  \ldots&=-\sqrt2\alpha^{3/2}+O(\alpha).
\end{align*}
Hence using~\eqref{eqn:beta`}
\begin{align}
  \beta_n'(\alpha)&=1-\left(\frac{-2n^2\pi^2+O(\alpha^{-1/2})}{-\sqrt2
  \alpha^{3/2}+O(\alpha)}\right)\notag\\
  &=1-\sqrt2n^2\pi^2\alpha^{-3/2}+O(\alpha^{-2}).
\end{align}
Therefore
\begin{align*}
  \rho_n'(\alpha)&=2\beta_n\beta_n'-2\alpha\\
  &=2n^2\pi^2+2\sqrt2n^2\pi^2\alpha^{-1/2}+O(\alpha^{-1}). \earlybox
\end{align*}

\subsection*{Acknowledgments}

I wish to thank E B Davies for his invaluable guidance in writing this paper
and his considerable time spent supervising my PhD. This research was funded
by an EPSRC studentship.

\bibliographystyle{amsalpha}

\begin{thebibliography}{VFFS66}

\bibitem[Aro51]{Aron}
N.~Aronszajn, \emph{Approximation methods for eigenvalues of completely
  continuous symmetric operators}, Proc. Symp. Spectral Theory Diff. Probl.
  (June-July 1950) (Oklahoma Agricultural and Mechanical College, Stillwater,
  Oklahoma), 1951, pp.~179--202.

\bibitem[BFS67]{BazlFoxStad}
N.~W. Bazley, D.~W. Fox, and J.~T. Stadter, \emph{Upper and lower bounds for
  the frequencies of rectangular clamped plates}, Z. angew. Math. Mech.
  \textbf{47} (1967), 191--198.

\bibitem[BM95]{BehnMert}
H.~Behnke and U.~Mertins, \emph{Eigenwertschranken f{\"u}r das {P}roblem der
  frei schwingenden rechteckigen {P}latte und {U}ntersuchungen zum
  {Aus\-weich\-ph\"anomen}}, Z. angew. Math. Mech. \textbf{75} (1995), no.~5,
  343--363.

\bibitem[BR72]{BaueReis}
L.~Bauer and E.~L. Reiss, \emph{Block five diagonal matrices and the fast
  numerical solution of the biharmonic equation}, Math. Comp. \textbf{26}
  (1972), no.~118, 311--326.

\bibitem[Dav95a]{Davi6}
E.~B. Davies, \emph{Nonlinear {S}chr{\"o}dinger operators and molecular
  structure}, J. Phys. A:Math. Gen. \textbf{28} (1995), 4025--4041.

\bibitem[Dav95b]{Davi5}
E.~B. Davies, \emph{Spectral theory and differential operators}, Cambridge
  University Press, 1995.

\bibitem[Fic78]{Fich}
G.~Fichera, \emph{Numerical and quantitative analysis}, London: Pitman, 1978,
  Trans. by S. Graffi.

\bibitem[Kat80]{Kato}
T.~Kato, \emph{Perturbation theory for linear operators}, second ed.,
  Springer-Verlag, 1980, Reprinted as a softcover edition in 1995.

\bibitem[KKM90]{KozlKondMazy}
V.~A. Kozlov, V.~A. Kondrat'ev, and V.~G. Maz'ya, \emph{On sign variation and
  the absence of ``strong'' zeros of solutions of elliptic equations},
  Translation in Math. USSR-Izv \textbf{34} (1990), no.~2, 337--353.

\bibitem[RS78]{ReedSimo4}
M.~Reed and B.~Simon, \emph{Methods of modern mathematical physics}, vol. IV:
  Analysis of operators, Academic press, 1978.

\bibitem[VFFS66]{DeviFichFuscScha}
L.~De Vito, G.~Fichera, A.~Fusciardi, and M~Sch{\"a}rf, \emph{Sul calcolo degli
  autovalori della piastra quadrata incastrata lungo il bordo}, Rend. Acc. Naz.
  Lincei. Cl. Sci. Fis. Matem. Nat. \textbf{XL} (1966), 725--733.

\bibitem[Wei37]{Wein}
A.~Weinstein, \emph{{\'E}tude des spectres des {\'e}quations aux
  d{\'e}riv{\'e}es partielles de la th{\'e}orie des plaques {\'e}lastiques},
  M{\'e}morial Sci. Math. (1937).

\bibitem[Wie]{Wien3}
C.~Wieners, \emph{Bounds for the {$N$} lowest eigenvalues of fourth-order
  boundary value problems}, Preprint 1995.

\bibitem[Wie96]{Wien2}
C.~Wieners, \emph{A numerical existence proof of nodal lines for the first
  eigenvalue of the plate equation}, Arch. Math. \textbf{66} (1996), 420--427.

\end{thebibliography}
\providecommand{\bysame}{\leavevmode\hbox to3em{\hrulefill}\thinspace}

\end{document}